\documentclass[11pt]{article}
\usepackage{amsmath}
\usepackage{hyperref}
\usepackage{amsmath}
\usepackage{amsfonts}
\usepackage{graphicx}
\usepackage{amssymb}
\usepackage{leftidx}
\usepackage{amsmath}
\usepackage{amsfonts}
\usepackage{graphicx}
\usepackage{amssymb}
\usepackage{leftidx}
\usepackage{mathrsfs}
\newtheorem{theorem}{Theorem}[section]
\newtheorem{lemma}[theorem]{Lemma}
\newtheorem{proposition}[theorem]{Proposition}

\newtheorem{definition}[theorem]{Definition}

\newtheorem{remark}[theorem]{Remark}

\numberwithin{equation}{section}

\begin{document}

\title{ Stochastic Recursive Optimal Control of McKean-Vlasov Type: A
Viscosity Solution Approach }
\author{Liangquan Zhang$\thanks{%
L. Zhang acknowledges the financial support partly by the National Nature
Science Foundation of China (Grant No. 12171053, 11701040, 11871010
\&61871058) and the Fundamental Research Funds for the Central Universities,
and the Research Funds of Renmin University of China (No. 23XNKJ05).}$ \\
{\small School of Mathematics, }\\
{\small Renmin University of China, Beijing 100872, China}}
\maketitle

\begin{abstract}
In this paper, we study a kind of optimal control problem for
forward-backward stochastic differential equations (FBSDEs for short) of
McKean--Vlasov type via the dynamic programming principle (DPP for short)
motivated by studying the infinite dimensional Hamilton--Jacobi--Bellman
(HJB for short) equation derived from the decoupling field of the FBSDEs
posed by Carmona and Delarue (\emph{Ann Probab}, 2015, \cite{cd15}). At the
beginning, the value function is defined by the solution to the controlled
BSDE alluded to earlier. On one hand, we can prove the value function is
deterministic function with respect to the initial random variable; On the
other hand, we can show that the value function is \emph{law-invariant},
i.e., depending on only via its distribution by virtue of BSDE property.
Afterward, utilizing the notion of differentiability with respect to
probability measures introduced by P.L. Lions \cite{Lions2012}, we are able
to establish a DPP for the value function in the Wasserstein space of
probability measures based on the application of BSDE approach,
particularly, employing the notion of stochastic \emph{backward semigroups}
associated with stochastic optimal control problems and It\^{o} formula
along a flow property of the conditional law of the controlled forward state
process. We prove that the value function is the unique viscosity solutions
of the associated generalized HJB equations in some separable Hilbert space.
Finally, as an application, we formulate an optimal control problem for
linear stochastic differential equations with quadratic cost functionals of
McKean-Vlasov type under nonlinear expectation, $g$-expectation introduced
by Peng \cite{Peng04} and derive the optimal feedback control explicitly by
means of several groups of Riccati equations.
\end{abstract}

\noindent \textbf{AMS subject classifications:} 93E20, 60H15, 60H30.

\noindent \textbf{Key words: }Dynamic programming principle,
Forward-backward McKean-Vlasov stochastic differential equations,
Hamilton-Jacobi-Bellman Equation, Wasserstein space, Value function,
Viscosity solutions. 
\newpage
\tableofcontents

\newpage

\section{Introduction\label{sec1}}

Since the fundamental works on McKean-Vlasov equations considered by McKean
Jr. \cite{mck} and Kac \cite{kac1, kac2}, there are huge literature focusing
on uncontrolled SDEs and obtaining the general propagation of chaos results.
Sustained academic attention in this field in the past decade was embodied
in the connection with the so-called mean-field game (MFG for short) theory,
introduced independently and simultaneously by Lasry and Lions in \cite{ll07}
and on Huang, Caines and Malham\'{e} \cite{hcm}. As a matter of fact, the
McKean--Vlasov equation naturally happens whenever one tries to interpret
the mechanism of the behavior of many symmetric agents, all of which
interact via the empirical distribution of their states, to find a Nash
equilibrium (competitive equilibrium) or a Pareto equilibrium (cooperative
equilibrium) (see \cite{bfy13, cdl13}).

As for the optimal control problem in the framework of the McKean--Vlasov
type (e.g. the classical mean--variance portfolio selection problem in
finance, see \cite{YZ1999}) also called mean field type, it is hard to
obtain the DPP for short due to the appearance of the law of the process in
the coefficients and nonlinear dependency with respect to it. Therefore,
problems like this actually belongs to a time inconsistent issue (see Bj\"{o}%
rk, Khapko and Murgoci, \cite{bkm17}\ and T. Bj\"{o}rk, M. Khapko and A.
Murgoci, \cite{bm14}, Hern\'{a}ndez and Possama\"{\i} \cite{HP20} references
therein). It is necessary to point out that, though the problem itself is
time inconsistent, it is possible to capture some form of the DPP by
extending the state space, for instance, Lauri\`{e}re and Pironneau \cite%
{LP14} (see also Bensoussan, Frehse and Yam \cite{bfy13, bfy15, bfy17})
adopted the hypothesis that the existence at all times of a density for the
marginal distribution of the state process and transformed the prime problem
into a density control problem with a family of deterministic controls.
Afterward, they were able to prove a DPP and derive a kind of Bellman's
equation in the space of density functions. Subsequently, Pham and Wei \cite%
{PW18} without the density assumption, obtained the DPP for closed-loop
controls. For open-loop controls, the related topics can be found in Cosso
and Pham \cite{CH2019} for McKean--Vlasov differential games and in
Bayraktar, Cosso and Pham \cite{BCP2018} for the so-called randomised DPP.
If containing the common noise, Pham and Wei \cite{PW2017} proved a DPP
where the control process is adapted to the common noise filtration.
Besides, Bouchard, Djehiche and Kharroubi \cite{BDK2020} investigated a
stochastic McKean-- Vlasov target problem, in which the controlled process
satisfies some target marginal constraints and established a general
geometric dynamic programming (see also \cite{ST2002}). Djete, Possama\"{\i}
and Tan \cite{DPT2022}, however, employed the measurable selection
techniques to prove the DPP. This approach can be used under very general
conditions, avoiding for instance the Markovian property or regularity
assumptions mentioned in the above literature. Alternatively, stochastic
maximum principle (SMP for short) for the control of McKean-Vlasov dynamics
has been considered in \cite{cd15, bdl11, ad10} for state dynamics depending
upon marginal distribution and in \cite{blm17, cz16} for conditional
McKean-Vlasov dynamics. Carmona and Delarue, \cite{cd15} established the SMP
and provided the sufficient conditions for existence of an optimal control.
Both necessary and sufficient parts of the stochastic Pontryagin principle
require the adjoint backward equation, and seeking the optimal control
problem turns de facto to the analysis of FBSDEs in which the marginal
distributions of the solutions appear in the coefficients of the equations.
In the same way, Li \cite{li12} studied a SMP for the mean-field controls.
For mean field linear-quadratic (LQ for short) problem, see Yong \cite{y11}.
A stochastic optimal control problem with delay and of mean-field type was
considered by Shen, Meng and Shi \cite{sms14}.

In history, nonlinear BSDEs of McKean--Vlasov type were investigated first
by Buckdahn, Djehiche Li and Peng \cite{bdlp09}. Since then, the related
results on the theory of FBSDEs of McKean--Vlasov type, as well as that of
the associated partial differential equations (PDEs for short) have been
widely studied. Buckdahn, Li and Peng \cite{blp09} proved not only an
existence and uniqueness theorem, but also a comparison theorem for
expectation type. Employing a BSDE approach, introduced by Peng \cite{pyfw97}%
, the authors also explored a probabilistic interpretation to related
nonlocal PDEs in finite dimensional (comparison of our paper see Remark \ref%
{com}). Subsequently, with the rapid development of the theory of mean field
FBSDEs, many researchers focused on stochastic control problems in the
McKean--Vlasov type. For instance, whenever FBSDEs involving the value
function, with frozen partial initial values, Hao and Li \cite{hl16}
considered an optimal control problem with systems of fully coupled
controlled mean field FBSDEs. Many other works in the mean-field area, see,
e.g., Kotelenez and Kurtz \cite{kk10} for stochastic PDES and Lasry and
Lions \cite{ll07} for Economics, Finance and game theory etc.

The objective of this paper is to study the optimal control of FBSDEs of
McKean--Vlasov type, establish the DPP by virtue of BSDE approach and then
derive the corresponding HJB equation. However, as mentioned before, from a
traditional point of view, the DPP for FBSDEs of McKean--Vlasov type does
not hold true anymore due to the appearance of expectation terms in
coefficients. To overcome this difficulty, one possible direction is
pointing to the generalized state space. In a nutshell, the idea behind this
is that we seize the flow property of the controlled conditional
distribution measurable process and transform the prime BSDE into one
depending on the conditional distribution of forward variable, under which
we are able to announce the DPP via backward semigroup. Nonetheless, there
are serval key points should be stated as follows:

\begin{itemize}
\item At first, unlike the cost functional considered in \cite{PW2017, cd15,
CH2019, DPT2022} defined by taking the expectation both on running function
and terminal function, which derives that the value function is
deterministic no matter of the selection of initial state is random or not
(see Proposition 3.1 in \cite{CH2019}), our cost functional in this paper is
defined by the solution to BSDE, so it can be random whenever the initial
state is measurable to the filtration $\mathcal{F}_{t}$ generated by Brown
motion at time $t\geq 0$. However, if it is measurable to a
\textquotedblleft rich enough\textquotedblright\ (see the interpretation
below) $\sigma $-algebra $\mathcal{G}$ independent of $\mathcal{F}_{t},$ we
are able to prove the value function is deterministic.

\item The value function will be depend on random variable $\xi $, but we
hope it depends on its distribution instead. Therefore, we should prove a
crucial law invariance property of the value function (defined by solution
to BSDE) (see Theorem \ref{inv}), which implies that the value function can
be considered as a function on the Wasserstein space of probability measures
on $\mathcal{P}_{2}\left( \mathbb{R}^{n}\right) $.

\item In an earlier work Barles, Buckdahn and Pardoux \cite{bbp1997},
introduced certain space of continuous functions and proved for in this
space the uniqueness of the viscosity solution of an integro-partial
differential equation$\footnote{%
This wide class of PDEs has some apparent features in contrast to those
derived from forward stochastic controlled systems (see \cite{PW2017,
DPT2022}), for instance, the inhomogeneous term depends not only on the
value function but also its gradient. These new changes make the
investigation more comlicated.}$ associated with decoupled FBSDEs with
jumps. Whereafter, Buckdahn, Li and Peng \cite{blp09} proved the uniqueness
for PDEs associated with a decoupled mean field FBSDEs (see Remark \ref%
{blp09r}) without jumps in the space of continuous functions with at most
polynomial growth. However, in our framework, we have to consider the PDEs
in infinite dimensional set and corresponding result may be borrowed from
\cite{FGS2015}.

\item All the works in \cite{PW2017, cd15, CH2019, DPT2022} consider the
optimal control problem under the classical expectation. In reality, the
objective expectation doesn't reflect people's preferences (for instance,
\cite{All1953, E1961}). There is a general consensus to define a certain
non-linear expectation, such as $g$-expectation (cf. Peng \cite{Peng04}).\
In this paper, we will revisit McKean-Vlasov LQ problem under $g$%
-expectation. We will see, on one hand, the admissible control set needs
more stringent integrability requirements due to the well-posedness for BSDE
regarding the terminal condition. On the other hand, the control problem
under generalized expectation is not necessarily equivalent to the classical
expectation. Hence, this topic is interesting in its own way.
\end{itemize}

The outline of the paper is organized as follows. After recalling briefly
some notations and stochastic framework to present an adequate and precise
definition of the tools that are used throughout the paper, we introduce in
Section \ref{sec2} the controlled FBSDEs of the McKean--Vlasov type of
stochastic control problem and study the properties of corresponding value
function. Next, in Section \ref{sec3}, we derive an infinite dimensional HJB
equation, and prove the viscosity property together with a uniqueness result
for the value function. In Section \ref{sec4}, we revisit a stochastic
McKean-Vlasov control problem under a non-linear expectation, $g$%
-expectation. Section \ref{sec5} states some concluding remarks. At last,
some proofs are displayed in Appendix \ref{Techlemmas}.

Let us introduce the following notations which are needed in what follows:

\begin{itemize}
\item For $r\geq 1,$ $\mathcal{P}_{r}\left( \mathbb{R}^{n}\right) $ denotes
the set probability measure $\mu $ on $\mathbb{R}^{n}$, which satisfies
square integrable, i.e. $\left\Vert \mu \right\Vert _{r}=\left[ \int_{%
\mathbb{R}^{n}}\left\vert x\right\vert ^{r}\mu \left( \mathrm{d}x\right) %
\right] ^{\frac{1}{r}}<\infty .$

\item For any $\mu \in \mathcal{P}_{2}\left( \mathbb{R}^{n}\right) ,$ we
denote by $L_{\mu }^{2}\left( \mathbb{R}^{d}\right) $ the set of measure
functions $\phi :\mathbb{R}^{n}\rightarrow \mathbb{R}^{d},$ which are square
integrable w.r.t. $\mu $ and define $\mu \left( \phi \right) :=\int_{\mathbb{%
R}^{n}}\phi \left( x\right) \mu \left( \mathrm{d}x\right) .$ Similarly, $%
L_{\mu \otimes \mu }^{2}\left( \mathbb{R}^{d}\right) $ denotes the set of
measure functions $\varphi :\mathbb{R}^{n}\times \mathbb{R}^{n}\rightarrow
\mathbb{R}^{d},$ which are square integrable w.r.t. $\mu \otimes \mu .$ We
define $\mu \otimes \mu \left( \varphi \right) :=\int_{\mathbb{R}^{n}\times
\mathbb{R}^{n}}\varphi \left( x,\bar{x}\right) \mu \left( \mathrm{d}x\right)
\mu \left( \mathrm{d}\bar{x}\right) .$

\item $C_{q}\left( \mathbb{R}^{n}\right) $ denotes the set of continuous
functions on $\mathbb{R}^{n}$ with quadratic growth, and for any $\psi \in
C_{q}\left( \mathbb{R}^{n}\right) ,$ we can define a mapping $\Pi _{\psi }:%
\mathcal{P}_{2}\left( \mathbb{R}^{n}\right) \rightarrow \mathbb{R}$ by $\Pi
_{\psi }\left( \mu \right) =\int_{\mathbb{R}^{n}}\psi \left( x\right) \mu
\left( \mathrm{d}x\right) .$

\item $L_{\mu }^{\infty }\left( \mathbb{R}^{d}\right) $ denotes the subset
of elements $\phi \in L_{\mu }^{2}\left( \mathbb{R}^{d}\right) $ which are
bounded w.r.t. $\mu $ with $\left\Vert \phi \right\Vert _{\infty }$ being
their essential supremum. Similarly, $L_{\mu \otimes \mu }^{\infty }\left(
\mathbb{R}^{d}\right) $ denotes the subset of elements $\varphi \in L_{\mu
\otimes \mu }^{2}\left( \mathbb{R}^{d}\right) $ which are bounded w.r.t. $%
\mu \otimes \mu $ with $\left\Vert \varphi \right\Vert _{\infty }$ being
their essential supremum.

\item Let $(\Omega ,\mathcal{F},\{\mathcal{F}_{t}\}_{t\geq 0},\mathbb{P})$
be a complete filtered probability space assume of the form $(\Omega
^{0}\times \Omega ^{1},\mathcal{F}^{0}\otimes \mathcal{G},\mathbb{P}%
^{0}\otimes \mathbb{P}^{1}),$ where $\left( \Omega ^{0},\mathcal{F}^{0},%
\mathbb{P}^{0}\right) $ supports a $d$-dimensional Brownian motion $W.$ Any
element $\omega \in \Omega $ can be expressed as $\omega =\left( \omega
^{0},\omega ^{1}\right) \in \Omega ^{0}\times \Omega ^{1}.$ Moreover, $W$ on
$\Omega $ can be written as $W\left( \omega ^{0},\omega ^{1}\right) =W\left(
\omega ^{0}\right) ,$ etc. Besides, we suppose that $\Omega ^{1}$\ is a
Polish space, $\mathcal{G}$ its Borel $\sigma $-field, $\mathbb{P}^{1}$ an
atomless probability measure on $\left( \Omega ^{1},\mathcal{G}\right) .$ $%
\mathcal{G}\supset \mathcal{N},$ where $\mathcal{N}$ is the set of all $%
\mathbb{P}$-null sets. We denote by $\mathbb{E}^{0}$ (resp. $\mathbb{E}^{1}$%
) the expectation under $\mathbb{P}^{0}$ (resp. $\mathbb{P}^{1}$). Let $%
\mathcal{F}^{0}=\left( \mathcal{F}_{t}^{0}\right) _{0\leq t\leq T}$ is the
natural filtration generated by $W,$ while $\mathbb{F}$ $=\left( \mathcal{F}%
_{t}\right) _{0\leq t\leq T}$ is be the filtration generated by $W$,
completed and augmented by a $\sigma $-field $\mathcal{G}$, i.e.,%
\begin{equation*}
\mathcal{F}_{t}^{0}=\sigma \left\{ W_{s}|s\leq t\right\} ,\text{ }\mathcal{F}%
_{t}=\mathcal{F}_{t+}^{0}\vee \mathcal{G},
\end{equation*}%
where%
\begin{equation*}
\mathcal{F}_{t+}^{0}=\bigcap_{s\geq t}\mathcal{F}_{s}^{0}.
\end{equation*}

\item We assume that the above sub-$\sigma $-field $\mathcal{G}$ of $%
\mathcal{F}$ which is independent of $\mathcal{F}_{\infty }^{0}$ and
\textquotedblleft rich enough\textquotedblright\ in the following sense:%
\begin{equation}
\mathcal{P}_{r}\left( \mathbb{R}^{n}\right) =\left\{ \mathbb{P}_{\xi
}\left\vert \xi \in L^{r}\left( \Omega ^{1},\mathcal{G},\mathbb{P}%
^{1}\right) \right. \right\} ,  \label{rich}
\end{equation}%
where $\mathbb{P}_{\xi }$ or $\mathcal{L}\left( \xi \right) $ denotes the
law of $\xi .$ From Lemma 2.1 in \cite{cg2022}, $\mathcal{G}$ is
\textquotedblleft rich enough\textquotedblright\ if and only if there exists
a $\mathcal{G}$-measurable random variable $U^{\mathcal{G}}:\Omega
^{1}\rightarrow \mathbb{R}$ having uniform distribution on $\left[ 0,1\right]
$. Particularly, if the probability space $\left( \Omega ^{1},\mathcal{G},%
\mathbb{P}^{1}\right) $ is \emph{atomless} (namely, for any $A\in \mathcal{G}
$ such that $\mathbb{P}^{1}\left( A\right) >0$ there exists $B\in \mathcal{G}
$, $B\subset A$, such that $0<\mathbb{P}^{1}\left( B\right) <\mathbb{P}%
^{1}\left( A\right) $), then these two mentioned properties holds. From
Corollary 7.16.1 in \cite{BS1978}, we know that, for each probability space $%
\left( \tilde{\Omega},\mathcal{\tilde{F}},\mathbb{\tilde{P}}\right) ,$
together with $\tilde{\Omega}$ uncountable, separable, complete metric
space, $\mathcal{\tilde{F}}$ its Borel $\sigma $-field, $\mathbb{\tilde{P}}$
an atomless probability, satisfies (\ref{rich}). For instance, (see Theorem
3.19 in \cite{K2002}), the probability space $\left( \tilde{\Omega},\mathcal{%
\tilde{F}},\mathbb{\tilde{P}}\right) $ can be chosen as $\left( \left[ 0,1%
\right] ,B\left( \left[ 0,1\right] \right) ,\lambda \right) $ where $\lambda
$ denotes the Lebesgue measure$.$

\item We define%
\begin{equation*}
\left( \mathcal{F}^{0}\right) _{s}^{t}\triangleq \sigma \left(
W_{r}-W_{t}|t\leq r\leq s\right) \vee \mathcal{N}\text{ },\text{ }s\in \left[
t,T\right] .
\end{equation*}

\item $L^{2}\left( \mathcal{G};\mathbb{R}^{n}\right) $ (resp. $L^{2}\left(
\mathcal{F}_{t}^{0}\vee \mathcal{G};\mathbb{R}^{n}\right) $) denotes the set
of $\mathbb{R}^{n}$-valued square integrable random variables on $\left(
\Omega ^{1},\mathcal{G},\mathbb{P}^{1}\right) $ (resp. $\left( \Omega ,%
\mathcal{F}_{t}^{0}\vee \mathcal{G},\mathbb{P}^{0}\right) $)$.$

\item As we have known, $\mathcal{P}_{2}\left( \mathbb{R}^{n}\right) $ is a
metric space equipped with the $2$-Wasserstein distance%
\begin{equation*}
\begin{array}{lll}
\mathcal{W}_{2}\left( \mu ^{1},\mu ^{2}\right) & \triangleq & \inf \Bigg \{%
\left( \int_{\mathbb{R}^{n}\times \mathbb{R}^{n}}\left\vert x-y\right\vert
^{2}\pi \left( \mathrm{d}x,\mathrm{d}y\right) \right) ^{\frac{1}{2}}:\pi \in
\mathcal{P}_{2}\left( \mathbb{R}^{n}\times \mathbb{R}^{n}\right) \text{ } \\
&  & \text{with marginals }\mu ^{1}\text{ and }\mu ^{2}\Bigg \} \\
&  & =\inf \Bigg \{\left( \mathbb{E}\left\vert X-X^{\prime }\right\vert
^{2}\right) ^{\frac{1}{2}}:X,\text{ }X^{\prime }\in L^{2}\left( \mathcal{O};%
\mathbb{R}^{n}\right) \\
&  & \text{with }\mathbb{P}_{X}=\mu ^{1},\text{ }\mathbb{P}_{X^{\prime
}}=\mu ^{2}\Bigg \}%
\end{array}%
\end{equation*}%
endowed with Borel $\sigma $-field $\mathcal{B}\left( \mathcal{P}_{2}\left(
\mathbb{R}^{n}\right) \right) .$ It follows from Proposition 7.1.5 in \cite%
{AGG2005} that $\left( \mathcal{P}_{2}\left( \mathbb{R}^{n}\right) ,\mathcal{%
W}_{2}\right) $ is a complete separable metric space\footnote{%
By Theorem 7.12 in \cite{V2003}, for $\left( \mu _{n}\right) _{n}$, $\mu
_{n}\in \mathcal{P}_{2}\left( \mathbb{R}^{n}\right) $, we have that $%
\mathcal{W}_{2}(\mu _{n},%
\mu
)\rightarrow 0$ if and
\par
only if, for every $\phi \in $ $C_{q}\left( \mathbb{R}^{n}\right) ,$ $\phi
\left( \mu _{n}\right) \rightarrow \phi \left( \mu \right) .$%
\par
{}}.

\item We need the following spaces:%
\begin{equation*}
\begin{array}{lll}
\mathcal{S}_{\mathcal{F}}^{2}\left( t,T;\mathbb{R}^{n}\right) & \triangleq & %
\Big \{\varphi \Bigg |\varphi :\Omega \times \left[ t,T\right] \rightarrow
\mathbb{R}^{n}\text{ is an }\mathcal{F}_{t}^{0}-\text{adapted process} \\
&  & \text{ satisfying }\mathbb{E}\left[ \sup_{0<s<T}\left\vert \varphi
_{s}\right\vert ^{2}\right] <\infty \Big \}; \\
\mathcal{M}_{\mathcal{F}}^{2}\left( t,T;\mathbb{R}^{n\times d}\right) &  &
\triangleq \Big \{\phi \Bigg |\phi :\Omega \times \left[ t,T\right]
\rightarrow \mathbb{R}^{n\times d}\text{ is }\mathcal{F}_{t}^{0}\text{%
-predictable process } \\
&  & \text{satisfying }\mathbb{E}\left[ \int_{t}^{T}\left\vert \phi
_{s}\right\vert ^{2}\mathrm{d}s\right] <\infty \Big \}.%
\end{array}%
\end{equation*}
\end{itemize}

\begin{remark}
\label{measure}It is straightforward to check that the map $\Pi _{\psi }$ is
$\mathcal{B}\left( \mathcal{P}_{2}\left( \mathbb{R}^{n}\right) \right) $%
-measurable, for any $\psi \in C_{q}\left( \mathbb{R}^{n}\right) $, by
virtue of a monotone class argument since it holds true whenever $\psi \in
C_{q}\left( \mathbb{R}^{n}\right) $. Besides, $\mathcal{B}\left( \mathcal{P}%
_{2}\left( \mathbb{R}^{n}\right) \right) $ coincides with the cylindrical $%
\sigma $-field $\sigma \left( \Pi _{\psi },\psi \in C_{q}\left( \mathbb{R}%
^{n}\right) \right) .$ Therefore, given a measurable space $\left( G,%
\mathcal{G}\right) $ and a map $\rho :G\rightarrow \mathcal{P}_{2}\left(
\mathbb{R}^{n}\right) $, $\rho $ is measurable if and only if the map $\Pi
_{\psi }\circ \rho :G\rightarrow \mathbb{R}$ is measurable, for any $\psi
\in C_{q}\left( \mathbb{R}^{n}\right) $.
\end{remark}



\section{Stochastic Recursive Control Systems of McKean-Vlasov Type}

\label{sec2}

In this section, we consider the following recursive optimal control
problem. For any $\left( t,\xi \right) \in \left[ 0,T\right] \times
L^{2}\left( \mathcal{G};\mathbb{R}^{n}\right) ,$ the dynamics are governed
by McKean-Vlasov controlled FBSDEs in the form:
\begin{equation}
\left\{
\begin{array}{rcl}
\mathrm{d}X_{s}^{t,\xi ;u} & = & b\left( X_{s}^{t,\xi ;u},\mathbb{P}%
_{X_{s}^{t,\xi ;u}}^{W},u_{s}\right) \mathrm{d}s+\sigma \left( X_{s}^{t,\xi
;u},\mathbb{P}_{X_{s}^{t,\xi ;u}}^{W},u_{s}\right) \mathrm{d}W_{s}, \\
\mathrm{d}Y_{s}^{t,\xi ;u} & = & -\mathbb{E}^{\mathcal{F}_{s}^{0}}\left[
f\left( X_{s}^{t,\xi ,u},\mathbb{P}_{X_{s}^{t,\xi ;u}}^{W},\Theta
_{s}^{t,\xi ,u},\mathbb{P}_{\Theta _{s}^{t,\xi ;u}},u_{s}\right) \right]
\mathrm{d}s+Z_{s}^{t,\xi ;u}\mathrm{d}W_{s}, \\
X_{t}^{t,\xi ;u} & = & \xi ,\text{ }Y_{T}^{t,\xi ;u}=\mathbb{E}^{\mathcal{F}%
_{T}^{0}}\left[ \Phi \left( X_{T}^{t,\xi ;u},\mathbb{P}_{X_{s}^{t,\xi
;u}}^{W}\right) \right] ,%
\end{array}%
\right.  \label{FBSDE}
\end{equation}%
where $\Theta _{s}^{t,\xi ,u}=\left( Y_{s}^{t,\xi ,u},Z_{s}^{t,\xi
,u}\right) $ and $\mathbb{E}^{\mathcal{F}_{s}^{0}}$ denotes the conditional
expectation w.r.t. $\mathcal{F}_{s}^{0}.$ Hereafter, $\mathbb{P}%
_{X_{s}^{t,\xi ,u}}^{W}$ denotes the regular conditional distribution of $%
X_{s}^{t,\xi ,u}$ given $\mathcal{F}^{0}$, and its realization at some $%
\omega ^{0}\in \Omega ^{0}$ reads as the law under $\mathbb{P}^{1}$ of the
random variable $X_{s}^{t,\xi ,u}\left( \omega ^{0},\cdot \right) $ on $%
\left( \Omega ^{1},\mathcal{F}^{1}\right) $, i.e. $\mathbb{P}_{X_{s}^{t,\xi
,u}}^{W}\left( \omega ^{0}\right) =\mathbb{P}_{X_{s}^{t,\xi ,u}\left( \omega
^{0},\cdot \right) }^{1}$. As for $\mathbb{P}_{\Theta _{s}^{t,\xi ;u}},$ it
represents as the law under $\mathbb{P}^{0}$ of the random variable $\left(
Y^{t,\xi ,u},Z^{t,\xi ,u}\right) $ on $\left( \Omega ^{0},\mathcal{F}%
^{0}\right) $, i.e. $\mathbb{P}_{\Theta _{s}^{t,\xi ;u}}=\mathbb{P}_{\Theta
_{s}^{t,\xi ;u}}^{0}.$

\begin{remark}
Note that $X^{t,\xi ;u}$ in (\ref{FBSDE}) is $\mathcal{F}$-adapted process,
while $Y^{t,\xi ;u}$ is $\mathcal{F}^{0}$-predictable process since the
conditional expectation $\mathbb{E}^{\mathcal{F}_{s}^{0}}$ appearing in
generator $f$ and $\mathbb{E}^{\mathcal{F}_{T}^{0}}$ in the terminal
function $\Phi $. We will prove that $Y_{t}^{t,\xi ;u}$ is deterministic for
any $\left( t,\xi \right) \in \left[ 0,T\right] \times L^{2}\left( \mathcal{G%
};\mathbb{R}^{n}\right) $ rather than random (see Lemma \ref{vdet} below).
\end{remark}


We explain the following notations:

\begin{itemize}
\item The control process $u,$ satisfying $u_{t}\in U,$ for $0\leq t\leq T$.
Here $U$, the subset of $\mathbb{R}^{k}$, called the control set$\footnote{%
As a matter of fact, the control set can be any Polish set $U$ equipped with
the distance $d_{U}$, satisfying w.l.o.g. $d_{U}<1.$ Recall Theorem 6.1 in
\cite{YZ1999} for classical LQ problem, the optimal control process can be
represented as a linear function in state variable. Inspired by this result,
in Section \ref{sec4}, we will postulate the control set to be a function
space satisfying Lipschitz condition.}$, is a \emph{compact} metric space.
Define the regular admissible control set as follows ($p>0$):%
\begin{equation}
\begin{array}{rcl}
\mathcal{U}_{ad}^{p}\left[ 0,T\right] & \triangleq & \Bigg \{u\left( \cdot
\right) :u_{s}\in U\text{ is an }\mathcal{F}_{s}\text{-progressively
measurable process } \\
&  & \text{satisfying }\mathbb{E}^{0}\left[ \int_{0}^{T}\left\vert
u_{s}\right\vert ^{2}\mathrm{d}s\right] ^{\frac{p}{2}}<\infty \Bigg \}.%
\end{array}%
\text{ }  \label{controlset}
\end{equation}%
Notice that $\mathcal{U}_{ad}^{2}\left[ 0,T\right] $ is separable metric
space under the distance
\begin{equation*}
\left\vert u-v\right\vert ^{2}=\mathbb{E}^{0}\left[ \int_{0}^{T}\left\vert
u_{s}\right\vert ^{2}\mathrm{d}s\right] ^{\frac{1}{2}}.
\end{equation*}

\item For simplicity, denote by $\Lambda _{t,\xi },$ the set of admissible
control under which the controlled systems (\ref{FBSDE}) have the initial
condition $\left( t,\xi \right) \in \left[ 0,T\right] \times L^{2}\left(
\mathcal{G};\mathbb{R}^{n}\right) ,$ being an $\mathcal{F}_{s}$%
-progressively measurable process valued in $U$ such that $u\in \mathcal{U}%
_{ad}^{2}\left[ 0,T\right] .$
\end{itemize}

Now suppose that $\Omega ^{0}=C\left( \mathbb{R}_{+},\mathbb{R}^{n}\right) ,$
the set of continuous functions from $\mathbb{R}_{+}$ into $\mathbb{R}^{n}$,
$\omega ^{0}$ is the canonical process, and $\mathbb{P}^{0}$ is the Wiener
measure. We introduce the shifted control processes (cf. \cite{CTT2016})
constructed by concatenation of paths:

\begin{definition}
\label{sc}Given any $u\in \mathcal{U}_{ad}^{2}\left[ 0,T\right] ,$ $\left( t,%
\tilde{\omega}^{0}\right) \in \left[ 0,T\right] \times \Omega ^{0},$ we
define%
\begin{equation*}
u_{s}^{t,\tilde{\omega}^{0}}\left( \omega ^{0}\right) =u_{s}\left( \tilde{%
\omega}^{0}\otimes _{t}\omega ^{0}\right) ,\text{ for }\left( s,\omega
^{0}\right) \in \left[ 0,T\right] \times \Omega ^{0}
\end{equation*}%
where
\begin{equation*}
\left( \tilde{\omega}^{0}\otimes _{t}\omega ^{0}\right) _{s}=\tilde{\omega}%
_{s}^{0}\mathbb{I}_{s<t}+\left( \tilde{\omega}_{s}^{0}+\omega
_{s}^{0}-\omega _{t}^{0}\right) \mathbb{I}_{s\geq t}.
\end{equation*}
\end{definition}

Notice that $\tilde{\omega}^{0}\otimes _{t}\omega ^{0}\in \Omega ^{0}.$
Moreover, $u^{t,\tilde{\omega}^{0}}\in \mathcal{U}_{t}^{2},$ where $\mathcal{%
U}_{t}^{2}$ denotes the subset of $\mathcal{U}_{ad}^{2}\left[ 0,T\right] $
which are independent of $\mathcal{F}^{0}$ under $\mathbb{P}^{0}.$

On the coefficients $b$ and $\sigma $, and also on the generator function $f$
and terminal function $\Phi $ in (\ref{FBSDE}), we impose the following
assumptions:

\begin{enumerate}
\item[\textbf{(A1)}] The coefficients $b:\mathbb{R}^{n}\times \mathcal{P}%
_{2}\left( \mathbb{R}^{n}\right) \times U\rightarrow \mathbb{R}^{n}$ and $%
\sigma :\mathbb{R}^{n}\times \mathcal{P}_{2}\left( \mathbb{R}^{n}\right)
\times U\rightarrow \mathbb{R}^{n\times d}$ are measurable functions$.$

\item[\textbf{(A2)}] The functions $b,$ $\sigma $ are Lipschitz continuous
in $\left( x,\mu ,u\right) $, i.e., there exists a positive constant $C$
such that, for any $x,$ $x^{\prime }\in \mathbb{R}^{n},$ $\mu ,$ $\mu
^{\prime }\in \mathcal{P}_{2}\left( \mathbb{R}^{n}\right) ,$ and $u,$ $%
u^{\prime }\in U$
\begin{eqnarray*}
&&\left\vert b\left( x,\mu ,u\right) -b\left( x^{\prime },\mu ^{\prime
},u^{\prime }\right) \right\vert +\left\vert \sigma \left( x,\mu ,u\right)
-\sigma \left( x^{\prime },\mu ^{\prime },u^{\prime }\right) \right\vert \\
&\leq &C\left( \left\vert x-x^{\prime }\right\vert +\mathcal{W}_{2}\left(
\mu ,\mu ^{\prime }\right) +\left\vert u-u^{\prime }\right\vert \right)
\end{eqnarray*}%
and
\begin{equation*}
\left\vert b\left( 0,\delta _{0},u\right) \right\vert +\left\vert \sigma
\left( 0,\delta _{0},u\right) \right\vert \leq C.
\end{equation*}
\end{enumerate}

\begin{remark}
It is necessary to point out that in \cite{pyfw97}, Peng (1997) first
studied the optimal control for FBSDEs under Markov framework assuming that
the coefficients $b$ and $\sigma $ have the linear growth with respect to $%
\left( x,u\right) $, are Lipschitz continuous in $x$ but H\"{o}lder
continuous with respect to $u$. Unlike the maximum principle needing some
smoothness of coefficients, we postulate that the coefficients are Lipschitz
continuous in $u.$ In the literature, $b$ and $\sigma $ fulfill the the
Lipchitz condition, uniformly w.r.t. $U$ and $\mathcal{P}\left( U\right) $
(see \cite{CH2019}) here $\delta _{0}$ stands for the Dirac measure with
mass at $0\in \mathbb{R}^{n}.$
\end{remark}

In order to ensure the well-defined of solution to (\ref{FBSDE}), we add the
following conditions:

\begin{enumerate}
\item[\textbf{(A3)}] The generator $f:\mathbb{R}^{n}\times \mathcal{P}%
_{2}\left( \mathbb{R}^{n}\right) \times \mathbb{R\times R}^{d}\times
\mathcal{P}_{2}\left( \mathbb{R}\times \mathbb{R}^{d}\right) \times
U\rightarrow \mathbb{R}$ and the terminal cost function $\Phi :\mathbb{R}%
^{n}\times \mathcal{P}_{2}\left( \mathbb{R}^{n}\right) \rightarrow \mathbb{R}
$ are continuous and measurable functions.

\item[\textbf{(A4)}] The functions $f,$ $\Phi $ satisfy the Lipschitz
condition: there exists a positive constant $C$ such that, for any $%
x,x^{\prime }\in \mathbb{R}^{n},$ $y,$ $y^{\prime }\in \mathbb{R},$ $z,$ $%
z^{\prime }\in \mathbb{R}^{d},$ $\mu ,$ $\mu ^{\prime }\in \mathcal{P}%
_{2}\left( \mathbb{R}^{n}\right) ,$ $\nu ,$ $\nu ^{\prime }\in \mathcal{P}%
_{2}\left( \mathbb{R}\times \mathbb{R}^{d}\right) $ and $u,$ $u^{\prime }\in
U$
\begin{eqnarray*}
&&\left\vert f\left( x,\mu ,y,z,\nu ,u\right) -f\left( x^{\prime },\mu
^{\prime },y^{\prime },z^{\prime },\nu ^{\prime },u^{\prime }\right)
\right\vert +\left\vert \Phi \left( x,\mu \right) -\Phi \left( x^{\prime
},\mu ^{\prime }\right) \right\vert \\
&\leq &C\big (\left\vert x-x^{\prime }\right\vert +\left\vert y-y^{\prime
}\right\vert +\left\vert z-z^{\prime }\right\vert +\mathcal{W}_{2}\left( \mu
,\mu ^{\prime }\right) +\mathcal{W}_{2}\left( \nu ,\nu ^{\prime }\right)
+\left\vert u-u^{\prime }\right\vert \big )
\end{eqnarray*}%
and%
\begin{equation*}
\left\vert f\left( x,\delta _{x},0,0,\delta _{\left( 0,0\right) },u\right)
\right\vert +\left\vert \Phi \left( x,\delta _{0}\right) \right\vert \leq
C\left( 1+\left\vert x\right\vert +\mathcal{W}_{2}\left( \delta _{x},\delta
_{0}\right) \right) .
\end{equation*}

\item[\textbf{(A5)}] Assume that the $\sigma $-algebra $\mathcal{G}$ is
countably generated up to null sets.
\end{enumerate}

\begin{remark}
In order to prove the continuity of value function, Cosso et al. \cite%
{cg2022}, impose a continuity assumption on the payoff functions $f$ and $%
\Phi $ which do not contain $y,z$ and $\nu .$ As a matter of fact, the
continuity of the value functions there is only used to investigate the
viscosity properties of the value functions in \cite{cg2022}. The proof of
the DPP nevertheless is not employed. However, in our paper, notice in
particular the value function is defined via the solution of BSDE, the first
issue is to ensure the existence and uniqueness of solution to FBSDEs (\ref%
{FBSDE}). Thus, the Lipschitz continuity of $f$ and $\Phi $ is generally
required. Besides, in order to prove the law-invariant of value function,
authors in \cite{cg2022} present the following important condition: The
functions $f$ and $\Phi $ are locally uniformly continuous in $\left( x,\mu
\right) $ uniformly with respect to $u$. That is, for each $\varepsilon >0$
and $n\in N,$ there exists $\delta =\delta _{\left( \varepsilon ,n\right) }$
such that, for every $u\in U,$ $\left( x,\mu \right) ,\left( x^{\prime },\mu
^{\prime }\right) \in \mathbb{R}^{n}\times \mathcal{P}_{2}\left( \mathbb{R}%
^{n}\right) ,$ whenever $\left\vert x\right\vert +\mathcal{W}_{2}\left(
\left( \mu ,\nu \right) ,\delta _{\left( 0,0\right) }\right) \leq n$ and $%
\left\vert x^{\prime }\right\vert +\mathcal{W}_{2}\left( \left( \mu ^{\prime
},\nu ^{\prime }\right) ,\delta _{\left( 0,0\right) }\right) \leq n$,%
\begin{equation*}
\left\vert x-x^{\prime }\right\vert +\mathcal{W}_{2}\left( \mu ,\mu ^{\prime
}\right) \leq \delta ,
\end{equation*}%
we have
\begin{eqnarray*}
\left\vert f\left( x,\mu ,u\right) -f\left( x^{\prime },\mu ^{\prime
},u\right) \right\vert &\leq &\varepsilon , \\
\left\vert \Phi \left( x,\mu \right) -\Phi \left( x^{\prime },\mu ^{\prime
}\right) \right\vert &\leq &\varepsilon .
\end{eqnarray*}%
Whilst, Djete et al. in \cite{DPT2022} put the following assumption, which
is weaker than the ones in \cite{cg2022}, Assumption $\left( \mathbf{A}%
_{f,g}\right) _{\text{cont}}$: There exists a constant $C>0$ such that%
\begin{eqnarray*}
\left\vert f\left( x,\mu ,u\right) \right\vert +\left\vert \Phi \left( x,\mu
\right) \right\vert &\leq &C\bigg (1+\left\vert x\right\vert ^{2}+\int_{%
\mathbb{R}^{n}}\left\vert \Theta _{1}\right\vert ^{2}\mu \left( \mathrm{d}%
\Theta _{1}\right) \\
&&+\int_{\mathbb{R}\times \mathbb{R}^{d}}\left\vert \Theta _{2}\right\vert
^{2}\nu \left( \mathrm{d}\Theta _{2}\right) +\left\vert u-u^{\prime
}\right\vert ^{2}\bigg ).
\end{eqnarray*}
\end{remark}

\begin{remark}
The reason why we assume that the $\sigma $-algebra $\mathcal{G}$ is
countably generated up to null set, is to ensures that the Hilbert space $%
L^{2}\left( \Omega ^{1},\mathcal{G},\mathbb{P}^{1};\mathbb{R}^{n}\right) $
is separable (cf. \cite{doob}, Page. 92).
\end{remark}

Under Assumptions (A1)-(A5), for any $\left( t,\xi \right) \in \left[ 0,T%
\right] \times L^{2}\left( \mathcal{G};\mathbb{R}^{n}\right) $, one can
prove the uniqueness and existences of solution to (\ref{FBSDE}) by
classical approach (see Theorem A.1 in \cite{Li2018}), denoted by $\left(
X_{s}^{t,\xi ;u},Y_{s}^{t,\xi ;u},Z_{s}^{t,\xi ;u}\right) \in \mathcal{S}_{%
\mathcal{F}}^{2}\left( t,T;\mathbb{R}^{n}\right) \times \mathcal{S}_{%
\mathcal{F}^{0}}^{2}\left( t,T;\mathbb{R}\right) \times \mathcal{M}_{%
\mathcal{F}^{0}}^{2}\left( t,T;\mathbb{R}^{d}\right) $. Moreover, applying
standard estimates, Burkholder--Davis--Gundy inequality$,$ we have%
\begin{equation}
\mathbb{E}\left[ \sup_{t\leq s\leq T}\left\vert X_{s}^{t,\xi ;u}\right\vert
^{2}\right] \leq C\left( 1+\mathbb{E}\left\vert \xi \right\vert ^{2}\right)
\label{est111}
\end{equation}%
where the constant $C$ is independent on $u,$ $\xi $ and $t$ (see Lemma 2.1
in \cite{CH2019})$.$ Under the Lipschitz condition in (A1)-(A2), employing
the standard arguments, namely, Burkholder--Davis--Gundy and Gronwall lemma,
we get the following estimate similar to the ones for controlled diffusion
processes (cf. Theorem 5.9 and Corollary 5.10, Chapter 2, in \cite{Kry1980}):%
\begin{equation}
\mathbb{E}\left[ \sup_{t\leq \theta \leq T}\left\vert X_{\theta }^{t,\xi
;u}-X_{\theta }^{s,\xi ^{\prime };u}\right\vert ^{2}\right] \leq C\left[
\mathbb{E}\left\vert \xi -\xi ^{\prime }\right\vert ^{2}+\left( 1+\mathbb{E}%
\left\vert \xi \right\vert ^{2}+\mathbb{E}\left\vert \xi ^{\prime
}\right\vert ^{2}\right) \left\vert s-t\right\vert \right] .  \label{sdeest1}
\end{equation}%
In fact, we have a more delicate estimate as follows:

\begin{lemma}
\label{l0}Under the Assumptions \emph{(A1)-(A2)}, for $p\geq 2$, there
exists a positive constant $C_{p}$ depending on the Lipschitz constants of $%
b $, $\sigma $ and $T$, such that for all $t\in \lbrack 0,T]$, $\xi ,\xi
^{\prime }\in L^{2}\left( \mathcal{G};\mathbb{R}^{n}\right) ,$ $u,u^{\prime
}\in \mathcal{U}_{ad}^{2}\left[ t,T\right] ,$ we have%
\begin{equation}
\mathbb{E}\left[ \sup_{t\leq \theta \leq T}\left\vert X_{\theta }^{t,\xi
;u}-X_{\theta }^{t,\xi ^{\prime };u^{\prime }}\right\vert ^{p}\right] \leq
C_{p}\left[ \mathcal{W}_{2}^{p}\left( \mathbb{P}_{\xi },\mathbb{P}_{\xi
^{\prime }}\right) +\Upsilon _{0,T}^{u,u^{\prime };p}\right] ,
\label{estsde1}
\end{equation}%
where
\begin{equation}
\Upsilon _{0,T}^{u,u^{\prime };p}=\mathbb{E}\left[ \int_{0}^{T}\left\vert
u_{s}-u_{s}^{\prime }\right\vert ^{p}\mathrm{d}\theta \right] .  \label{u}
\end{equation}
\end{lemma}


\paragraph{Proof.}

For read's convenience, we display the main steps as follows. To this end,
denote $\hat{X}_{\theta }=X_{\theta }^{t,\xi ;u}-X_{\theta }^{t,\xi ^{\prime
};u^{\prime }},$ $\hat{\xi}=\xi -\xi ^{\prime }$ and Applying
Burkholder-Davis-Gundy inequality, Lipschitz condition (A2) and Gronwall
lemma, we derive that%
\begin{equation}
\mathbb{E}\left[ \sup_{t\leq \theta \leq s}\left\vert \hat{X}_{\theta
}\right\vert ^{p}\right] \leq C_{p}\left\vert \hat{\xi}\right\vert
^{p}+C_{p}\int_{t}^{s}\Big [\left\vert u_{s}-u_{s}^{\prime }\right\vert ^{p}+%
\mathcal{W}_{2}^{p}\left( \mathbb{P}_{X_{\theta }^{t,\xi ;u}}^{W},\mathbb{P}%
_{X_{\theta }^{t,\xi ^{\prime };u^{\prime }}}^{W}\right) \Big ]\mathrm{d}%
\theta .  \label{est01}
\end{equation}%
From the definition of $\mathcal{W}_{2},$ we have $\mathcal{W}_{2}^{2}\left(
\mathbb{P}_{X_{\theta }^{t,\xi ;u}}^{W},\mathbb{P}_{X_{\theta }^{t,\xi
^{\prime };u^{\prime }}}^{W}\right) \leq \mathbb{E}\left[ \left\vert \hat{X}%
_{\theta }\right\vert ^{2}\right] .$ Immediately, Gronwall lemma yields%
\begin{equation*}
\sup_{\theta \in \left[ t,T\right] }\mathcal{W}_{2}^{2}\left( \mathbb{P}%
_{X_{\theta }^{t,\xi ;u}}^{W},\mathbb{P}_{X_{\theta }^{t,\xi ^{\prime
};u^{\prime }}}^{W}\right) \leq C\mathbb{E}\left[ \left\vert \hat{\xi}%
\right\vert ^{2}+\Upsilon _{0,T}^{u,u^{\prime };p}\right]
\end{equation*}%
and taking into account the arbitrariness of $\xi ,\xi ^{\prime }\in
L^{2}\left( \mathcal{G};\mathbb{R}^{n}\right) ,i=1,2,$ we get%
\begin{equation*}
\sup_{\theta \in \left[ t,T\right] }\mathcal{W}_{2}^{2}\left( \mathbb{P}%
_{X_{\theta }^{t,\xi ;u}}^{W},\mathbb{P}_{X_{\theta }^{t,\xi ^{\prime
};u^{\prime }}}^{W}\right) \leq C\left[ \mathcal{W}_{2}^{2}\left( \mathbb{P}%
_{\xi },\mathbb{P}_{\xi ^{\prime }}\right) +\Upsilon _{0,T}^{u,u^{\prime };p}%
\right] .
\end{equation*}%
This together with (\ref{est01}) completes the proof. \hfill $\Box $

\begin{remark}
Immediately, we have
\begin{equation*}
\sup_{\theta \in \left[ t,T\right] }\mathcal{W}_{2}^{2}\left( \mathbb{P}%
_{X_{\theta }^{t,\xi ;u}}^{W},\mathbb{P}_{X_{\theta }^{t,\xi ^{\prime
};u}}^{W}\right) \leq C\mathcal{W}_{2}^{2}\left( \mathbb{P}_{\xi },\mathbb{P}%
_{\xi ^{\prime }}\right) ,
\end{equation*}%
which means for uncontrolled forward system, its solution depends on the law
of $\xi .$
\end{remark}

For any $u\in \Lambda _{t,\xi },$ we know that $\left\{ X_{s}^{t,\xi
;u}\right\} _{t\leq s\leq T},$ which is $\mathcal{F}$-adapted and $W$ is a $%
\left( \Omega ,\mathcal{F}\right) $ Wiener process. Therefore, it derives
that%
\begin{equation*}
\mathbb{P}_{X_{s}^{t,\xi ;u}}^{W}\left( \mathrm{d}x\right) =\mathbb{P}\left(
\left. X_{s}^{t,\xi ;u}\in \mathrm{d}x\right\vert \mathcal{F}^{0}\right) =%
\mathbb{P}\left( \left. X_{s}^{t,\xi ;u}\in \mathrm{d}x\right\vert \mathcal{F%
}_{s}^{0}\right)
\end{equation*}%
thus for any $\psi \in C_{q}\left( \mathbb{R}^{n}\right) $
\begin{equation*}
\mathbb{P}_{X_{s}^{t,\xi ;u}}^{W}\left( \psi \right) =\mathbb{E}\left[
\left. \psi \left( X_{s}^{t,\xi ;u}\right) \right\vert \mathcal{F}^{0}\right]
=\mathbb{E}\left[ \left. \psi \left( X_{s}^{t,\xi ;u}\right) \right\vert
\mathcal{F}_{s}^{0}\right] ,\text{ }t\leq s\leq T,
\end{equation*}%
which indicates that $\left\{ \mathbb{P}_{X_{s}^{t,\xi ;u}}^{W}\left( \psi
\right) \right\} _{t\leq s\leq T}$ is $\left( \mathcal{F}^{0}\right) _{t\leq
s\leq T}$-measurable. We define the map $\mathbb{P}_{X_{s}^{t,\xi ;u}\left(
\omega \right) }^{W}:\omega \rightarrow \mathcal{P}_{2}\left( \mathbb{R}%
^{n}\right) .$ Then, $\mathbb{P}_{X_{s}^{t,\xi ;u}}^{W}$ is $\mathcal{F}%
_{s}^{0}$-measurable if and only if $\mathbb{P}_{X_{s}^{t,\xi ;u}}^{W}\left(
\psi \right) $ is $\mathcal{F}_{s}^{0}$-measurable. Therefore, we report
that $\left\{ \mathbb{P}_{X_{s}^{t,\xi ;u}}^{W}\right\} _{t\leq s\leq T}$ is
$\left( \mathcal{F}_{s}^{0}\right) _{t\leq s\leq T}$-adapted. Moreover, we
have the following flow property for forward SDE in (\ref{FBSDE})
\begin{equation*}
X_{r}^{t,\xi ;u}=X_{r}^{s,X_{s}^{t,\xi ;u};u},\text{ for any }r\in \left[ s,T%
\right] ,\text{ }\mathbb{P}\text{-a.s.,}
\end{equation*}%
which implies that
\begin{equation*}
\mathbb{P}_{X_{r}^{t,\xi ;u}}^{W}=\mathbb{P}_{X_{r}^{s,X_{s}^{t,\xi
;u};u}}^{W},\text{ for any }r\in \left[ s,T\right] ,\text{ }\mathbb{P}\text{%
-a.s.}
\end{equation*}

\begin{definition}
Given any $s\in \left[ t,T\right] ,$ $\xi \in L^{2}\left( \mathcal{G};%
\mathbb{R}^{n}\right) ,$ with $\mathbb{P}_{\xi }^{W}=\mu $ and $u\in
\mathcal{U}_{ad}^{2}\left[ 0,T\right] ,$ we define
\begin{equation}
\rho _{s}^{t,\mu ;u}=\mathbb{P}_{X_{s}^{t,\xi ;u}}^{W}.  \label{law}
\end{equation}
\end{definition}

From Lemma 3.1 in \cite{PW2017}, we know that $\rho _{s}^{t,\mu ;u}$ defines
a square integrable $\mathcal{F}^{0}$-progressive continuous process in $%
\mathcal{P}_{2}\left( \mathbb{R}^{n}\right) $. The map $\left( s,t,\omega
^{0},\mu ,u\right) \in \left[ 0,T\right] \times \left[ 0,T\right] \times
\Omega ^{0}\times \mathcal{P}_{2}\left( \mathbb{R}^{n}\right) \times
\mathcal{U}_{ad}^{2}\left[ 0,T\right] \rightarrow \rho _{s}^{t,\mu ;u}\in
\mathcal{P}_{2}\left( \mathbb{R}^{n}\right) $ is measurable and admits the
following flow property:%
\begin{equation}
\rho _{s}^{t,\mu ;u}=\rho _{s}^{\tau ,\rho _{\tau }^{t,\mu ;u};u^{\tau }},
\label{flow}
\end{equation}%
where $u^{\tau }$ is defined in Definition \ref{sc} for all $\mathcal{F}^{0}$%
-stopping time $\tau $ valued in $\left[ t,T\right] ,$ such a set is denoted
by $\mathbb{T}_{t,T}^{0}.$

We now deal with
\begin{eqnarray*}
&&\mathbb{E}^{\mathcal{F}_{s}^{0}}\left[ f\left( X_{s}^{t,\xi ,u},\mathbb{P}%
_{X_{s}^{t,\xi ;u}}^{W},\Theta _{s}^{t,\xi ,u},\mathbb{P}_{\Theta
_{s}^{t,\xi ;u}},u_{s}\right) \right] \\
&=&\rho _{s}^{t,\mu ;u}\left( f\left( \cdot ,\rho _{s}^{t,\mu ;u},\Theta
_{s}^{t,\xi ,u},\mathbb{P}_{\Theta _{s}^{t,\xi ;u}},u_{s}\right) \right) \\
&=&f^{0}\left( \rho _{s}^{t,\mu ;u},\Theta _{s}^{t,\xi ,u},\mathbb{P}%
_{\Theta _{s}^{t,\xi ;u}},u_{s}\right) ,
\end{eqnarray*}%
and
\begin{eqnarray*}
\mathbb{E}^{\mathcal{F}_{T}^{0}}\left[ \Phi \left( X_{T}^{t,\xi ;u},\mathbb{P%
}_{X_{s}^{t,\xi ;u}}^{W}\right) \right] &=&\rho _{T}^{t,\mu ;u}\left( \Phi
\left( \cdot ,\rho _{T}^{t,\mu ;u}\right) \right) \\
&=&\Phi ^{0}\left( \rho _{T}^{t,\mu ;u}\right)
\end{eqnarray*}%
with the functions $f^{0}:\left( \mu ,y,z,\mathbb{\pi },u\right) \in
\mathcal{P}_{2}\left( \mathbb{R}^{n}\right) \times \mathbb{R}\times \mathbb{R%
}^{d}\times \mathcal{P}_{2}\left( \mathbb{R}\times \mathbb{R}^{d}\right)
\times U\rightarrow \mathbb{R}$ and $\Phi ^{0}:\mu \in \mathcal{P}_{2}\left(
\mathbb{R}^{n}\right) \rightarrow \mathbb{R},$ defined by
\begin{eqnarray*}
f^{0}\left( \mu ,y,z,\mathbb{\pi },u\right) &=&\mu \left( f\left( \cdot ,\mu
,y,z,\mathbb{\pi },u\right) \right) \\
&=&\int_{\mathbb{R}^{n}}f\left( x,\mu ,y,z,\mathbb{\pi },u\right) \mu \left(
\mathrm{d}x\right)
\end{eqnarray*}%
and
\begin{equation*}
\Phi ^{0}\left( \mu \right) =\mu \left( \Phi \left( \cdot ,\mu \right)
\right) =\int_{\mathbb{R}^{n}}\Phi \left( x,\mu \right) \mu \left( \mathrm{d}%
x\right) .
\end{equation*}%
Then FBSDEs (\ref{FBSDE}) can be equivalently expressed as follows:%
\begin{equation}
\left\{
\begin{array}{lll}
\mathrm{d}X_{s}^{t,\xi ;u} & = & b\left( X_{s}^{t,\xi ;u},\mathbb{P}%
_{X_{s}^{t,\xi ;u}}^{W},u_{s}\right) \mathrm{d}s+\sigma \left( X_{s}^{t,\xi
;u},\mathbb{P}_{X_{s}^{t,\xi ;u}}^{W},u_{s}\right) \mathrm{d}W_{s}, \\
X_{t}^{t,\xi ;u} & = & \xi , \\
\mathrm{d}Y_{s}^{t,\xi ;u} & = & -f^{0}\left( \rho _{s}^{t,\mu ;u},\Theta
_{s}^{t,\xi ,u},\mathbb{P}_{\Theta _{s}^{t,\xi ;u}},u_{s}\right) \mathrm{d}%
s+Z_{s}^{t,\xi ;u}\mathrm{d}W_{s}, \\
Y_{T}^{t,\xi ;u} & = & \Phi ^{0}\left( \rho _{T}^{t,\mu ;u}\right) .%
\end{array}%
\right.  \label{fbsde3}
\end{equation}%
However, the new expression of BSDE in (\ref{fbsde3}) means that the
solution $Y^{t,\xi ;u}$ depends on $\xi $ only through its distribution $\mu
=\rho _{t}^{t,\mu ;u}=\mathcal{L}\left( \xi \right) .$ So (\ref{fbsde3}) can
be written as:%
\begin{equation}
\left\{
\begin{array}{rcl}
\mathrm{d}X_{s}^{t,\xi ;u} & = & b\left( X_{s}^{t,\xi ;u},\rho _{s}^{t,\mu
;u},u_{s}\right) \mathrm{d}s+\sigma \left( X_{s}^{t,\xi ;u},\rho _{s}^{t,\mu
;u},u_{s}\right) \mathrm{d}W_{s}, \\
X_{t}^{t,\xi ;u} & = & \xi , \\
\mathrm{d}Y_{s}^{t,\xi ;u} & = & -f^{0}\left( \rho _{s}^{t,\mu ;u},\Theta
_{s}^{t,\mu ;u},\mathbb{P}_{\Theta _{s}^{t,\mu ;u}},u_{s}\right) \mathrm{d}%
s+Z_{s}^{t,\xi ;u}\mathrm{d}W_{s}, \\
Y_{T}^{t,\xi ;u} & = & \Phi ^{0}\left( \rho _{T}^{t,\mu ;u}\right) .%
\end{array}%
\right.  \label{fbsde4}
\end{equation}

\begin{remark}
\label{blp09r} Buckdahn, Li and Peng \cite{blp09} (see also \cite{hl16} for
fully coupled case) considered mean-field FBSDEs in a more general
framework, with general coefficient like%
\begin{equation}
\left\{
\begin{array}{lll}
\mathrm{d}X_{s}^{t,\zeta } & = & \mathbb{E}^{\prime }\left[ b\left( s,\left(
X_{s}^{0,x_{0}}\right) ^{\prime },X_{s}^{t,\zeta }\right) \right] \mathrm{d}%
s+\mathbb{E}^{\prime }\left[ \sigma \left( s,\left( X_{s}^{0,x_{0}}\right)
^{\prime },X_{s}^{t,\zeta }\right) \right] \mathrm{d}W_{s}, \\
X_{t}^{t,\zeta } & = & \zeta \in L^{2}\left( \Omega ,\mathcal{F}_{t}^{0},%
\mathbb{P};\mathbb{R}^{n}\right) , \\
\mathrm{d}Y_{s}^{t,\zeta } & = & -\mathbb{E}^{\prime }\left[ f\left(
s,\left( X_{s}^{0,x_{0}}\right) ^{\prime },X_{s}^{t,\zeta },\left(
Y_{s}^{0,x_{0}}\right) ^{\prime },Y_{s}^{t,\zeta },Z_{s}^{t,\zeta }\right) %
\right] \mathrm{d}s+Z_{s}^{t,\zeta }\mathrm{d}W_{s}, \\
Y_{T}^{t,\zeta } & = & \mathbb{E}^{\prime }\left[ \Phi \left( \left(
X_{s}^{0,x_{0}}\right) ^{\prime },X_{T}^{t,\zeta }\right) \right] ,%
\end{array}%
\right.  \label{blp09}
\end{equation}%
where the driving coefficient of (\ref{blp09}) has to be interpreted as
follows (the other coefficients can be done in the same way):%
\begin{eqnarray*}
&&\mathbb{E}^{\prime }\left[ f\left( s,\left( X_{s}\right) ^{\prime
},X_{s},\left( Y_{s}\right) ^{\prime },Y_{s},Z_{s}\right) \right] \\
&=&\int_{\Omega }f\left( \omega ^{\prime },\omega ,s,X_{s}\left( \omega
^{\prime }\right) ,X_{s}\left( \omega \right) ,Y_{s}\left( \omega ^{\prime
}\right) ,Y_{s}\left( \omega \right) ,Z_{s}\left( \omega \right) \right)
P\left( \mathrm{d}\omega ^{\prime }\right) .
\end{eqnarray*}%
Clearly, $f$ does not dependent on the law of solutions. For the convenience
of readers, let us depict briefly the procedure to solve (\ref{blp09}).
Under certain assumptions (for instance, Lipschitz and linear growth
conditions), SDE in (\ref{blp09}) with $\left( t,\zeta \right) =\left(
0,x_{0}\right) $ admits a unique strong solution. Once we know $%
X_{s}^{0,x_{0}}$, SDE in (\ref{blp09}) with $\left( t,\zeta \right) $
becomes a classical equation. Analogously, we first deal with BSDE in (\ref%
{blp09}) for $\left( t,\zeta \right) =\left( 0,x_{0}\right) ,$ from Theorem
3.1 in \cite{blp09} that there exists a unique solution $\left(
Y_{s}^{0,x_{0}},Z_{s}^{0,x_{0}}\right) _{t\leq s\leq T}.$ Once we have $%
\left( Y_{s}^{0,x_{0}},Z_{s}^{0,x_{0}}\right) _{t\leq s\leq T},$ BSDE in (%
\ref{blp09}) reduces to a classical BSDE of the following type.
\end{remark}

\begin{remark}
The merit to write the dependence on $\mu $ in $Y_{s}^{t,\mu ;u}$ is
embodied in the convenience to derive the DPP by means of the flow property
of $\rho ^{t,\mu ;u}$ (see \ref{p2}). As claimed in Introduction, the DPP
for mean-field FBSDEs (\ref{blp09}) does not hold true anymore because of
the presence of expectation terms in coefficients. Nevertheless, it is
possible to capture some form of the DPP by extending the state space,
namely, the flow of the conditional law of the controlled forward state
process. So this formulation is different from \cite{bdlp09, blp09}. In
addition, if employing the conditional law, one must handle the dependence
of value function on conditional law, i.e. law-invariant (see Theorem \ref%
{inv} below).
\end{remark}

Pardoux and Peng (\cite{PP1}, 1990) introduced the following nonlinear BSDE:%
\begin{equation}
Y_{t}=\xi +\int_{t}^{T}f\left( s,Y_{s},Z_{s}\right) \mathrm{d}%
s-\int_{t}^{T}Z_{s}\mathrm{d}W_{s}  \label{BSDE1992}
\end{equation}%
where $\xi \in L_{\mathcal{F}_{T}}^{2}\left( 0,T;\mathbb{R}\right) ,$ under
certain conditions for the generator $f$, BSDE (\ref{BSDE1992}) admits a
unique strong adapted solution $\left( Y\left( \cdot \right) ,Z\left( \cdot
\right) \right) \in \mathcal{S}_{\mathcal{F}^{0}}^{2}\left( t,T;\mathbb{R}%
\right) \times \mathcal{M}_{\mathcal{F}^{0}}^{2}\left( t,T;\mathbb{R}%
^{d}\right) $. Besides, we need the following estimations for BSDE, whose
proof can be seen in Proposition 3.2 of Briand et al. \cite{Bri2003}.

\begin{lemma}
\label{l1}Let $\left( y^{i},z^{i}\right) ,$ $i=1,2,$ be the solution to the
following
\begin{equation}
y^{i}\left( t\right) =\xi ^{i}+\int_{t}^{T}f^{i}\left(
s,y_{s}^{i},z_{s}^{i}\right) \mathrm{d}s-\int_{t}^{T}z_{s}^{i}\mathrm{d}%
W_{s},  \label{estbdsde}
\end{equation}%
where $\xi ^{i}\in L^{2}\left( \mathcal{F}_{T};\mathbb{R}\right) $
satisfying $\mathbb{E}\left[ \left\vert \xi ^{i}\right\vert ^{\beta }\right]
<\infty ,$ $f^{i}=f^{i}\left( s,y^{i},z^{i}\right) :\left[ 0,T\right] \times
\mathbb{R}\times \mathbb{R}^{d}\rightarrow \mathbb{R}$ satisfies Lipschitz
in $\left( y,z\right) $ and
\begin{equation*}
\mathbb{E}\left[ \left( \int_{t}^{T}\left\vert f^{i}\left(
s,y_{s}^{i},z_{s}^{i}\right) \right\vert \mathrm{d}s\right) ^{\beta }\right]
<\infty .
\end{equation*}%
Then, for some $\beta \geq 2,$ there exists a positive constant $C_{\beta }$
(depending on $\beta ,$ $T$ and Lipschitz constant) such that
\begin{eqnarray}
&&\mathbb{E}\left[ \sup_{0\leq t\leq T}\left\vert
y_{t}^{1}-y_{t}^{2}\right\vert ^{\beta }+\left( \int_{0}^{T}\left\vert
z_{s}^{1}-z_{s}^{2}\right\vert ^{2}\mathrm{d}s\right) ^{\frac{\beta }{2}}%
\right]  \notag \\
&\leq &C_{\beta }\mathbb{E}\Bigg [\left\vert \xi ^{1}-\xi ^{2}\right\vert
^{\beta }+\left( \int_{t}^{T}\left\vert f^{1}\left(
s,y_{s}^{1},z_{s}^{1}\right) -f^{2}\left( s,y_{s}^{1},z_{s}^{1}\right)
\right\vert \mathrm{d}s\right) ^{\beta }\Bigg ].  \label{estbdsde1}
\end{eqnarray}%
Particularly, whenever putting $\xi ^{2}=0,$ $f^{2}=0,$ one has%
\begin{equation}
\mathbb{E}\left[ \sup_{0\leq t\leq T}\left\vert y_{t}^{1}\right\vert ^{\beta
}+\left( \int_{0}^{T}\left\vert z_{s}^{1}\right\vert ^{2}\mathrm{d}s\right)
^{\frac{\beta }{2}}\right] \leq C_{\beta }\mathbb{E}\Bigg [\left\vert \xi
^{1}\right\vert ^{\beta }+\left( \int_{t}^{T}\left\vert f^{1}\left(
s,0,0\right) \right\vert \mathrm{d}s\right) ^{\beta }\Bigg ].
\label{estbdsde2}
\end{equation}
\end{lemma}

\begin{remark}
Notice that Briand et al. \cite{Bri2003} considered the BSDE (\ref{BSDE1992}%
) under rather weak assumptions on the data. They proved the existence and
uniqueness of solutions in $L^{p},$ for $p>1$, extending the condition of
Lipschitz to the case where the monotonicity conditions are satisfied.
Moreover, under an additional assumption, they reported an existence and
uniqueness result for BSDE on a fixed time interval, when the data $\xi $
are only in $L^{1}$. However, in our paper, we only study the case $L^{p},$
for $p\geq 2.$
\end{remark}

\begin{lemma}
\label{estbsde}Suppose that the Assumptions \emph{(A1)-(A4)} are in force.
Then, for $p\geq 2,$ there exists a constant $C_{p}>0$ depending on the
Lipschitz constant of $b,\sigma ,\Phi ,f,$ such that, for $\xi ,\xi ^{\prime
}\in L^{2}\left( \mathcal{F}_{t};\mathbb{R}^{n}\right) $%
\begin{equation}
\mathbb{E}\left[ \sup_{t\leq s\leq T}\left\vert Y_{s}^{t,\xi ;u}\right\vert
^{p}+\left( \int_{t}^{T}\left\vert Z_{s}^{t,\xi ;u}\right\vert ^{2}\mathrm{d}%
s\right) ^{\frac{p}{2}}\right] \leq C_{p},  \label{best1}
\end{equation}%
\begin{eqnarray}
&&\mathbb{E}\left[ \sup_{t\leq s\leq T}\left\vert Y_{s}^{t,\xi
;u}-Y_{s}^{t,\xi ^{\prime };u^{\prime }}\right\vert ^{p}+\left(
\int_{t}^{T}\left\vert Z_{s}^{t,\xi ;u}-Z_{s}^{t,\xi ^{\prime };u^{\prime
}}\right\vert ^{2}\mathrm{d}s\right) ^{\frac{p}{2}}\right]  \notag \\
&\leq &C_{p}\left[ \Upsilon _{0,T}^{u,u^{\prime };p}+\mathcal{W}%
_{2}^{p}\left( \mathbb{P}_{\xi },\mathbb{P}_{\xi ^{\prime }}\right) \right] ,
\label{best2}
\end{eqnarray}%
where $\Upsilon _{0,T}^{u,u^{\prime };p}$ is defined in (\ref{est01}), and%
\begin{equation}
\int_{t}^{T}\mathcal{W}_{2}^{2}\left( \mathbb{P}_{\Theta _{s}^{t,\xi ,u}},%
\mathbb{P}_{\Theta _{s}^{t,\xi ^{\prime },u}}\right) \leq C_{p}\left[
\Upsilon _{0,T}^{u,u^{\prime };2}+\mathcal{W}_{2}^{2}\left( \mathbb{P}_{\xi
},\mathbb{P}_{\xi ^{\prime }}\right) \right] .  \label{best3}
\end{equation}
\end{lemma}

The proof is displayed in Appendix.

We now state a comparison theorem for BSDE which is vital to derive the DPP.

\begin{lemma}
\label{combsde}Let $\left( Y^{i},Z^{i}\right) ,$ $i=1,2,$ be the solution to
the following%
\begin{equation}
Y_{t}^{i}=\zeta ^{i}+\int_{t}^{T}g^{i}\left( Y_{s}^{i},Z_{s}^{i},\mathbb{P}%
_{\left( Y_{s}^{i},Z_{s}^{i}\right) }\right) \mathrm{d}s-%
\int_{t}^{T}Z_{s}^{i}W_{s}\mathrm{d}s,\text{ }i=1,2,  \label{cbsde}
\end{equation}%
where $\zeta ^{i}\in L^{2}\left( \mathcal{F}_{T};\mathbb{R}^{n}\right) $ and
$g^{i}:\mathbb{R\times R}^{d}\times \mathcal{P}_{2}\left( \mathbb{R}%
^{n}\times \mathbb{R}^{n\times d}\right) \rightarrow \mathbb{R}.$ Assume
that $g^{i}$ satisfies the Lipschitz condition: there exists a positive
constant $C$ such that, for any $y,y^{\prime }\in \mathbb{R},z,z^{\prime
}\in \mathbb{R}^{d},$ and $\nu ,\nu ^{\prime }\in \mathcal{P}_{2}\left(
\mathbb{R}\right) $
\begin{equation*}
\left\vert g^{i}\left( y,z,\nu \right) -g^{i}\left( y^{\prime },z^{\prime
},\nu ^{\prime }\right) \right\vert \leq C\big (\left\vert y-y^{\prime
}\right\vert +\left\vert z-z^{\prime }\right\vert +\mathcal{W}_{2}\left( \nu
,\nu ^{\prime }\right) \big ),
\end{equation*}%
Then, BSDE (\ref{cbsde}) (see Theorem A.1 in \cite{Li2018}) admits a unique
adapted solution $\left( Y^{i},Z^{i}\right) \in \mathcal{S}_{\mathcal{F}%
}^{2}\left( t,T;\mathbb{R}\right) \times \mathcal{M}_{\mathcal{F}}^{2}\left(
t,T;\mathbb{R}^{n\times d}\right) $, respectively, for $i=1,2$. Furthermore,
if the following conditions hold, (i) $\zeta ^{1}\geq \zeta ^{2},$ $\mathbb{P%
}$-a.s.$;$ (ii) $g^{1}\left( Y,Z,\gamma \right) \geq g^{2}\left( Y,Z,\gamma
\right) ,$ for any $\left( Y,Z,\gamma \right) \in \mathbb{R\times R}%
^{d}\times \mathcal{P}_{2}\left( \mathbb{R}\right) ;$ (iii) For all $\eta
_{1},\eta _{2}\in L^{2}\left( \mathcal{F}_{t}^{0}\vee \mathcal{G};\mathbb{R}%
\right) $ and all $y\in \mathbb{R},$ $z\in \mathbb{R}^{d}$
\begin{equation*}
g^{i}\left( y,z,\mathbb{P}_{\eta _{1}}\right) -g^{i}\left( y,z,\mathbb{P}%
_{\eta _{2}}\right) \leq C\left( \mathbb{E}\left[ \left( \left( \eta
_{1}-\eta _{2}\right) ^{+}\right) ^{2}\right] \right) ^{\frac{1}{2}}.
\end{equation*}%
Then, we have $Y_{t}^{1}\geq Y_{t}^{2},$ $\mathbb{P}$-a.s.$,$ for all $t\in %
\left[ 0,T\right] .$
\end{lemma}

The proof can be found in Appendix.

In \cite{PW2017}, the authors applied the measurable selection arguments and
shifted control technique to derive DPP, while Cosso et al. \cite{CH2019}
proved the DPP inspired by the classical approach for deterministic
differential games. Under general conditions (Markovian or non-Markovian),
Djete et al. (cf. \cite{DPT2022}) recently developed the classical
measurable selection, conditioning and concatenation arguments for common
noise framework, and established the DPP. In our paper, we will deal with
the control problem driven by FBSDEs of McKean-Vlasov type. The application
of BSDE approach, in particular, the use of the notion of stochastic
backward semigroups introduced by Peng \cite{pyfw97} allows a
straightforward proof of a DPP for value functions associated with
stochastic optimal control problems. We define the family of backward
semigroups associated with the FBSDEs (\ref{fbsde4}).

Given the initial data $\left( t,\xi \right) \in \left[ 0,T\right] \times
L^{2}\left( \mathcal{G};\mathbb{R}^{n}\right) $, a positive number $\delta
\leq T-t$, an admissible control process $u\in \Lambda _{t,\xi }$ and a
real-valued random variable $\eta \in L^{2}\left( \mathcal{F}_{t+\delta
}^{0};\mathbb{R}^{n}\right) $, we set%
\begin{equation*}
\mathbb{G}_{s,t+\delta }^{t,\xi ;u}\left[ \eta \right] :=\mathcal{Y}%
_{s}^{t,\xi ;u},\text{ }t\leq s\leq t+\delta ,
\end{equation*}%
where $\left( \mathcal{Y}_{s}^{t,\xi ;u},\mathcal{Z}_{s}^{t,\xi ;u}\right)
_{t\leq s\leq t+\delta }$ is the solution to the following BSDE of
McKean-Vlasov type with the time horizon $t+\delta :$%
\begin{eqnarray}
\mathcal{Y}_{s}^{t,\xi ;u} &=&\eta +\int_{s}^{t+\delta }f^{0}\left( \rho
_{r}^{t,\xi ;u},\mathcal{Y}_{r}^{t,\xi ;u},\mathcal{Z}_{r}^{t,\xi ;u},%
\mathbb{P}_{\left( \mathcal{Y}_{r}^{t,\xi ;u},\mathcal{Z}_{r}^{t,\xi
;u}\right) },u_{r}\right) \mathrm{d}r  \notag \\
&&-\int_{s}^{t+\delta }\mathcal{Z}_{r}^{t,\xi ;u}\mathrm{d}W_{r},
\label{bsg}
\end{eqnarray}%
where $\rho _{r}^{t,\xi ;u}$ is defined in (\ref{law}).

Apparently, the solution of $\left( Y_{s}^{t,\mu ,u},Z_{s}^{t,\mu ,u}\right)
_{t\leq s\leq T}$ of BSDE in (\ref{fbsde4}) admits%
\begin{equation}
\mathbb{G}_{t,T}^{t,\mu ;u}\left[ \Phi ^{0}\left( \rho _{T}^{t,\mu
;u}\right) \right] =\mathbb{G}_{t,t+\delta }^{t,\mu ;u}\left[ Y_{t+\delta
}^{t,\mu ;u}\right] .  \label{flow3}
\end{equation}%
On the one hand, we are in the Markovian case; On the other hand, Lemma \ref%
{l0} and Lemma \ref{l1} once again indicates that, the processes $\left(
X_{s}^{t,\xi ,u},Y_{s}^{t,\xi ,u},Z_{s}^{t,\xi ,u}\right) $ depend on $\xi $
only through its distribution, which means $\left( X_{s}^{t,\xi
,u},Y_{s}^{t,\xi ,u},Z_{s}^{t,\xi ,u}\right) _{t\leq s\leq T}$ and $\left(
X_{s}^{t,\xi ^{\prime },u},Y_{s}^{t,\xi ^{\prime },u},Z_{s}^{t,\xi ^{\prime
},u}\right) _{t\leq s\leq T}$ are indistinguishable as long as $\xi $ and $%
\xi ^{\prime }$ have the same distribution. Besides, from (\ref{flow3}), we
observe that $Y_{t+\delta }^{t,\mu ;u}$ is $\mathcal{F}_{t+\delta }^{0}\vee
\mathcal{G}$-adapted.

Next we will prove that $Y_{t}^{t,\xi ,u}$ is not random whenever $\left(
t,\xi \right) \in \left[ 0,T\right] \times L^{2}\left( \mathcal{G};\mathbb{R}%
^{n}\right) .$ To this end, we first state a general result.

\begin{lemma}
\label{deter}Assume that $\tilde{f}$ satisfies the Lipschitz condition and $%
\int_{0}^{T}\left\vert \tilde{f}\left( \cdot ,s,0,\delta _{\left( 0,0\right)
}\right) \right\vert \mathrm{d}s\in L^{2}\left( \mathcal{F}_{T};\mathbb{R}%
\right) $. For any fixed $t_{0}\in \left[ 0,T\right] $, if for every $\left(
y,z,\mu \right) \in \mathbb{R\times R}^{d}\times \mathcal{P}_{2}\left(
\mathbb{R}\times \mathbb{R}^{d}\right) $, $\tilde{f}\left( \cdot ,y,z,\mu
\right) $-adapted on the interval $\left[ t_{0},T\right] $ and $\zeta \in
L^{2}\left( \left( \mathcal{F}^{0}\right) _{T}^{t_{0}};\mathbb{R}\right) ,$
then the solution $\left( y_{s},z_{s}\right) _{t_{0}\leq s\leq T}$ to the
following BSDE is also $\mathcal{F}_{s}^{t_{0}}$-adapted on the interval $%
\left[ t_{0},T\right] $:%
\begin{equation}
y_{s}=\zeta +\int_{s}^{T}\tilde{f}\left( s,y_{s},z_{s},\mathbb{P}_{\left(
y_{s},z_{s}\right) }\right) \mathrm{d}s-\int_{s}^{T}z_{s}\mathrm{d}W_{s}.
\label{debsde}
\end{equation}%
Particularly, $y_{t_{0}}$ is a deterministic constant.
\end{lemma}

\paragraph{Proof.}

We set $W_{s}^{t_{0}}=W_{s}-W_{t_{0}}.$ Consider the following BSDE on the
interval $\left[ t_{0},T\right] $:%
\begin{equation}
\tilde{y}_{s}=\zeta +\int_{s}^{T}\tilde{f}\left( s,\tilde{y}_{s},\tilde{z}%
_{s},\mathbb{P}_{\left( \tilde{y}_{s},\tilde{z}_{s}\right) }\right) \mathrm{d%
}s-\int_{s}^{T}\tilde{z}_{s}\mathrm{d}W_{s}^{t_{0}}.  \label{debsde2}
\end{equation}%
It is easy to check that $\int_{t}^{T}z_{s}\mathrm{d}W_{s}=\int_{s}^{T}%
\tilde{z}_{s}\mathrm{d}W_{s}^{t_{0}}.$ Then, according to existence and
uniqueness of (\ref{debsde}), the solution of $\left( y_{s},z_{s}\right)
_{t_{0}\leq s\leq T}$ coincides with $\left( \tilde{y}_{s},\tilde{z}%
_{s}\right) _{t_{0}\leq s\leq T}$ which indicates that $\left(
y_{s},z_{s}\right) _{t_{0}\leq s\leq T}$ is $\mathcal{F}_{s}^{t_{0}}$%
-adapted. \hfill $\Box $

Note that if the initial state $\tilde{\xi}\in L^{2}\left( \mathcal{F}%
_{t}^{0}\vee \mathcal{G};\mathbb{R}^{n}\right) $ with $\mathcal{L}\left(
\tilde{\xi}\right) =\mu ,$ and $u\in \Lambda _{t,\xi },$ then we claim that $%
Y_{t}^{t,\tilde{\xi};u}$ is random.

\begin{lemma}
\label{rdet}For any $\left( t,\tilde{\xi}\right) \in \left[ 0,T\right]
\times L^{2}\left( \mathcal{F}_{t}^{0}\vee \mathcal{G};\mathbb{R}^{n}\right)
$ with $\mathcal{L}\left( \tilde{\xi}\right) =\mu ,$ and $u\in \Lambda _{t,%
\tilde{\xi}},$ $Y_{t}^{t,\tilde{\xi};u}$ is $\mathcal{F}_{t}^{0}\vee
\mathcal{G}$ measurable.
\end{lemma}

The proof can be found in Appendix.

However, when we are in the Markovian case (deterministic coefficients). The
following lemma claims that $Y_{t}^{t,\xi ;u}$ is \emph{deterministic }if%
\emph{\ }$\xi \in L^{2}\left( \mathcal{G};\mathbb{R}^{n}\right) $.

\begin{lemma}
\label{vdet}Under the Assumptions \emph{(A1)-(A5)}, for any $\left( t,\xi
\right) \in \left[ 0,T\right] \times L^{2}\left( \mathcal{G};\mathbb{R}%
^{n}\right) ,$ $\left. Y_{s}^{t,\xi ;u}\right\vert _{s=t}$ is deterministic.
\end{lemma}

\paragraph{Proof.}


Repeating the method in Lemma \ref{rdet}, consider a simple random variable $%
\xi $ of the form: $\xi =\sum_{i=1}^{N}I_{A_{i}}x^{i}$, and a control
process $u$ of the form $u_{\left( \cdot \right)
}=\sum_{i=1}^{N}I_{A_{i}}u_{\cdot }^{i}$, where $\left\{ A_{i}\right\}
_{i=1,\ldots ,N}$ is a finite partition of $\left( \Omega ,\mathcal{G}%
\right) ,$ $x^{i}\in \mathbb{R}^{n}$ and $u_{\cdot }^{i}$ is $\left(
\mathcal{F}^{0}\right) _{\cdot }^{t}$-adapted, for $1\leq i\leq N. $ For
each $i,$ set $\left( X_{s}^{i},Y_{s}^{i},Z_{s}^{i}\right) =\left(
X_{s}^{i},\Theta _{s}^{i}\right) =\left(
X_{s}^{t,x^{i},u^{i}},Y_{s}^{t,x^{i},u^{i}},Z_{s}^{t,x^{i},u^{i}}\right) ,$
for $s\in \left[ t,T\right] $. Then $\left( X^{i},Y^{i},Z^{i}\right) $ is
the solution to the FBSDEs%
\begin{equation}
\left\{
\begin{array}{lll}
\mathrm{d}X_{s}^{i} & = & b\left( X_{s}^{i},\mathbb{P}%
_{X_{s}^{i}}^{W},u_{s}^{i}\right) \mathrm{d}s+\sigma \left( X_{s}^{i},%
\mathbb{P}_{X_{s}^{i}}^{W},u_{s}^{i}\right) \mathrm{d}W_{s}, \\
\mathrm{d}Y_{s}^{i} & = & -f^{0}\left( \mathbb{P}_{X_{s}^{i}}^{W},\Theta
_{s}^{i},\mathbb{P}_{\Theta _{s}^{i}},u_{s}^{i}\right) \mathrm{d}s+Z_{s}^{i}%
\mathrm{d}W_{s}, \\
X_{t}^{i} & = & x^{i},\text{ }Y_{T}^{i}=\Phi ^{0}\left( X_{T}^{i},\mathbb{P}%
_{X_{T}^{i}}^{W}\right) .%
\end{array}%
\right.  \label{df1}
\end{equation}%
Therefore,
\begin{equation}
Y_{t}^{t,\xi ;u}=\sum_{i=1}^{N}I_{A_{i}}Y_{t}^{i}.  \label{df2}
\end{equation}%
Clearly, from Lemma \ref{deter}, each $Y_{t}^{i}$ is deterministic. Now
taking the conditional expectation $\mathbb{E}^{\mathcal{F}_{t}^{0}}$ on
both sides of (\ref{df2}) yields
\begin{eqnarray*}
\mathbb{E}^{\mathcal{F}_{t}^{0}}\left[ Y_{t}^{t,\xi ;u}\right] &=&\mathbb{E}%
^{\mathcal{F}_{t}^{0}}\left[ \sum_{i=1}^{N}I_{A_{i}}Y_{t}^{i}\right] \\
&=&\sum_{i=1}^{N}Y_{t}^{i}\mathbb{E}^{\mathcal{F}_{t}^{0}}\left[ I_{A_{i}}%
\right] \\
&=&\sum_{i=1}^{N}Y_{t}^{i}P\left( A_{i}\right) \text{ since }\mathcal{G}%
\text{ is independent of }\mathcal{F}_{t}^{0}\text{.}
\end{eqnarray*}%
Therefore, for simple random variables, we have the desired result. For a
general $\xi \in L^{2}\left( \mathcal{G};\mathbb{R}^{n}\right) ,$ the rest
proof is analogous in Lemma \ref{rdet}, so we omit it. The proof is
complete. \hfill $\Box $

Now the cost functional corresponding to the control $u\in \Lambda _{t,\xi }$
is defined to be%
\begin{equation}
J\left( t,\xi ;u\right) \triangleq \left. Y_{s}^{t,\xi ;u}\right\vert _{s=t}%
\text{ }  \label{cost1}
\end{equation}%
for $\left( t,\xi \right) \in \left[ 0,T\right] \times L^{2}\left( \mathcal{G%
};\mathbb{R}^{n}\right) $ with $\rho _{t}^{t,\mu ;u}=\mu .$


It follows from the uniqueness of the solution of (\ref{fbsde4}) and flow
property of $\rho ^{t,\mu ;u}$ (see (\ref{flow})),
\begin{equation}
Y_{t+\delta }^{t,\mu ;u}=Y_{t+\delta }^{t+\delta ,\rho _{t+\delta }^{t,\mu
;u};u}=J\left( t+\delta ,\rho _{t+\delta }^{t,\mu ;u};u\right) .
\label{flow2}
\end{equation}

\begin{remark}
When $\xi =x\in \mathbb{R}^{n},$ let $\delta _{x}$ be the Dirac measure with
mass at $x\in \mathbb{R}^{n},$ we instantly have $J\left( t,\delta
_{x};u\right) =Y_{t}^{t,\delta _{x};u}$ since $Y_{t}^{t,\delta _{x};u}$ is
deterministic by virtue of Lemma \ref{deter}, which reduces to the classical
recursive stochastic optimal control problem when $f^{0}$ does not contain
the distributions of $X$ and $\left( Y,Z\right) $ (see \cite{pyfw97,
Zhangsingular}).
\end{remark}

\begin{remark}
A special case is when $f$ doesn't depend on $\left( y,z,\pi \right) \in
\mathbb{R}\times \mathbb{R}^{d}\times \mathcal{P}_{2}\left( \mathbb{R}\times
\mathbb{R}^{d}\right) $, i.e., $f=f\left( x,\mu ,u\right) $. Now if we still
suppose (A1)-(A2) hold then the cost functional (\ref{cost1}) turns into
\begin{eqnarray*}
J\left( t,\xi ;u\right) &=&Y_{t}^{t,\xi ;u} \\
&=&\mathbb{E}\left[ \Phi ^{0}\left( \rho _{T}^{t,\mu ;u}\right)
+\int_{t}^{T}f^{0}\left( \rho _{s}^{t,\mu ;u},u_{s}\right) \mathrm{d}s\right]
,
\end{eqnarray*}%
which has been studied by Djete et al. \cite{DPT2022} with $W\equiv 0,$ Pham
and Wei \cite{PW2017} when $B\left( \cdot \right) \equiv 0,$ Cosso and Pham
\cite{CH2019} with single control process. Note that whenever we define the
cost functional like:
\begin{eqnarray*}
J\left( t,\xi ;u\right) &=&Y_{t}^{t,\xi ;u} \\
&=&\mathbb{E}\left[ \Phi ^{0}\left( \rho _{T}^{t,\mu ;u}\right)
+\int_{t}^{T}c\left( \rho _{s}^{t,\mu ;u}\right) Y_{s}^{t,\xi
;u}+f^{0}\left( \rho _{s}^{t,\mu ;u},u_{s}\right) \mathrm{d}s\right] ,
\end{eqnarray*}%
it actually corresponds to a kind of extended Feynman--Kac formula.
\end{remark}

We now consider the optimal control problem as follows:

\noindent \textbf{Problem MVS.} The McKean-Vlasov recursive control problem
is to derive, over admissible control set $\Lambda _{t,\xi }$, the following
value function%
\begin{equation}
V\left( t,\xi \right) =\inf_{u\in \Lambda _{t,\xi }}J\left( t,\xi ;u\right) ,%
\text{ }  \label{d1}
\end{equation}%
for $\left( t,\xi \right) \in \left[ 0,T\right] \times L^{2}\left( \mathcal{G%
};\mathbb{R}^{n}\right) $ with $\rho _{t}^{t,\mu ;u}=\mathbb{P}%
_{X_{t}^{t,\xi ;u}}^{W}=\mu .$

Before introducing the PDE by virtue of DPP for the value function, we at
present study some properties of the value function, which is useful to
characterize the viscosity solution.

\begin{lemma}
\label{estvalue}Assume that \emph{(A1)-(A5)} hold. There exists a constant $%
C>0$ such that, for any $t,$ $t^{\prime }\in \left[ 0,T\right] ,$ $\xi ,$ $%
\xi ^{\prime }\in L^{2}\left( \mathcal{G};\mathbb{R}^{n}\right) $%
\begin{eqnarray}
\left\vert V\left( t,\xi \right) -V\left( t,\xi ^{\prime }\right)
\right\vert &\leq &C\mathbb{E}\left[ \left\vert \xi -\xi ^{\prime
}\right\vert \right] ,  \label{est1} \\
\left\vert V\left( t,\xi \right) -V\left( t^{\prime },\xi \right)
\right\vert &\leq &C\left( 1+\mathbb{E}\left[ \left\vert \xi \right\vert %
\right] \right) \left\vert t-t^{\prime }\right\vert ^{\frac{1}{2}},
\label{est2} \\
\left\vert V\left( t,\xi \right) \right\vert &\leq &C\left( 1+\mathbb{E}%
\left[ \left\vert \xi \right\vert \right] \right) .  \label{est3}
\end{eqnarray}
\end{lemma}

The proof can be seen in the Appendix.

Apparently, for any $u\in \Lambda _{t,\xi },$ the control processes $u$ may
depend on $\xi $ since it is measurable with respect to $\mathcal{G}$.
Therefore, we can only say that the law $\mathbb{P}_{X_{s}^{t,\xi ;u}}^{W}$
of $X_{s}^{t,\xi ;u}$ depends on the joint law of $\left( \xi
,u,(W_{s}-W_{t})_{s\in \left[ t,T\right] }\right) $. The next lemma will
show that the value function $V$ is \emph{law-invariant}, i.e., depend on $%
\xi $ only via its distribution.

\begin{remark}
For the classical stochastic McKean--Vlasov control problem, for instance,
Cosso and Pham \cite{CH2019}, Cosso et al. \cite{cg2022} they defined the
value function as a function of $\left( t,\xi \right) $ and then proved the
so called \textquotedblleft law invariance property\textquotedblright ,
stating that the value function does not depend on $\xi $ but rather on its
law. In \cite{cg2022}, the authors shown that one can't expect this
invariance property to hold without any regularity conditions on the reward
functionals and terminal function (see \cite{cg2022}, Remark 3.7). Djete et
al. \cite{DPT2022} studied the same topic via measurable selection
techniques and presented weaker conditions than the ones in \cite{cg2022}
(see Remark 6 in \cite{cg2022}).
\end{remark}

\begin{theorem}
\label{inv}We assume that the Assumptions \emph{(A1)-(A5)} are in force. For
any $\xi ,$ $\xi ^{\prime }\in L^{2}\left( \mathcal{F}_{t}^{0}\vee \mathcal{G%
};\mathbb{R}^{n}\right) $ with $\mathbb{P}_{\xi }=\mathbb{P}_{\xi ^{\prime
}},$ we have%
\begin{equation*}
V\left( t,\xi \right) =V\left( t,\xi ^{\prime }\right) ,\text{ }\mathbb{P}%
\text{-a.s.}
\end{equation*}
\end{theorem}

\begin{remark}
Note that in (\ref{cost1}) and (\ref{d1}), the initial condition $\xi $
belongs to $L^{2}\left( \mathcal{G};\mathbb{R}^{n}\right) .$ On one hand, we
can prove that $Y_{t}^{t,\xi ;u}$ is a constant (see Lemma \ref{vdet}) which
is corresponding to a version of deterministic HJB equation (see Section \ref%
{sec3}); As a matter of fact, in \cite{PW2017, DPT2022, CH2019, cg2022,
YZ1999} etc, the reward functional is defined in the following way: No
matter which spaces the initial condition $\xi $ belongs to, the cost
functional is all performed after taking the absolute expectation $\mathbb{E}%
,$ which immediately reduces that the value function is always non-random;
Notwithstanding, the value function in our paper is defined via the solution
to BSDE. Hence very different from the aforementioned works, the current
value function is random whenever $\xi \in L^{2}\left( \mathcal{F}%
_{t}^{0}\vee \mathcal{G};\mathbb{R}^{n}\right) $ (see Lemma \ref{rdet})$.$
On the other hand, as pointed out in \cite{cg2022}, Remark 3.7, the main
reason why we take $\xi \in L^{2}\left( \mathcal{F}_{t}^{0}\vee \mathcal{G};%
\mathbb{R}^{n}\right) $ is to state and prove the DPP (more details see
Theorem 3.4 in \cite{cg2022}). In our paper, for instance, the initial
condition at time $t$ is $\xi \in L^{2}\left( \mathcal{G};\mathbb{R}%
^{n}\right) $, nevertheless, $X_{s}^{t,\xi ;u}$ in general, is not $\mathcal{%
G}$ measurable!
\end{remark}

For reader's convenience, we exhibit two kinds of lemmata (Lemma \ref{inv1}
and Lemma \ref{inv2}) from \cite{cg2022} to proceed our proof, which are
crucial to prove the law-invariant. In \cite{cg2022}, the authors deal with
the state space in a Hilbert space $H$ given by the set of $C\left( \left[
0,T\right] ;H\right) $ of continuous $H$-valued functions on $\left[ 0,T%
\right] .$ Hence, reduced to the finite-dimensional case, i.e., $H=\mathbb{R}%
^{n},$ and in the Markovian case, we have

\begin{lemma}
\label{inv1}For any $u\in \mathcal{U}_{ad}^{2}\left[ 0,T\right] $ and $\xi
\in L^{2}\left( \mathcal{F}_{t}^{0}\vee \mathcal{G};\mathbb{R}^{n}\right) $,
assume that there exists an $\mathcal{F}_{t}^{0}\vee \mathcal{G}$-measurable
random variable $U_{\xi }$, having uniform distribution on $[0,1]$ and being
independent of $\xi $. Then, there exists a measurable function%
\begin{equation*}
\mathbf{a:}\left( \left[ 0,T\right] \times \Omega \times \mathbb{R}%
^{n}\times \left[ 0,1\right] ,\mathcal{B}\left( \left[ 0,T\right] \right)
\otimes \mathcal{F}_{t}^{0}\otimes \mathcal{B}\left( \mathbb{R}^{n}\right)
\otimes \mathcal{B}\left( \left[ 0,1\right] \right) \right) \rightarrow
\left( U,\mathscr{U}\right) ,
\end{equation*}%
where $\mathscr{U}$ denotes the Borel $\sigma $-field of $U$ such that $%
\left\{ \mathbf{a}_{s}\left( \xi ,U_{\xi }\right) \right\} _{0\leq s\leq
T}\in \mathcal{U}_{ad}^{2}\left[ 0,T\right] ,$ $\mathbf{a}_{s}\left( \cdot
,\cdot \right) $ is constant for every $s<t,$ moreover,
\begin{equation*}
\left( \xi ,\left( \mathbf{a}_{s}\left( \xi ,U_{\xi }\right) \right) _{0\leq
s\leq T},\left( W_{s}-W_{t}\right) _{0\leq s\leq T}\right) \overset{%
\mathscr{Law}}{=}\left( \xi ,\left( u_{s}\right) _{0\leq s\leq T},\left(
W_{s}-W_{t}\right) _{0\leq s\leq T}\right) ,
\end{equation*}%
where $\mathscr{Law}$ denotes the equality in law.
\end{lemma}

\begin{lemma}
\label{inv2}For any $\xi \in L^{2}\left( \mathcal{F}_{t}^{0}\vee \mathcal{G};%
\mathbb{R}^{n}\right) $, assume that there exists $\left\{ x_{1},\ldots
,x_{m}\right\} \subset \mathbb{R}^{n},$ with $x_{i}\neq x_{j},$ if $i\neq j,$
such that
\begin{equation*}
\mathbb{P}_{\xi }=\sum_{i=1}^{m}p_{i}\delta _{x_{i}},
\end{equation*}%
where $\delta _{x_{i}}$ is the Dirac measure at $x_{i}$ and $p_{i}>0$, with $%
\sum_{i=1}^{m}p_{i}=1$. Then, there exists an $\mathcal{F}_{t}^{0}\vee
\mathcal{G}$-measurable random variable $U_{\xi }$ having uniform
distribution on $\left[ 0,1\right] $ and being independent of the $\mathcal{F%
}_{t}^{0}\vee \mathcal{G}$-measurable map%
\begin{equation*}
\tilde{\xi}:\Omega \rightarrow \mathbb{R}^{n},\omega \rightarrow \xi \left(
\omega \right) .
\end{equation*}
\end{lemma}

\paragraph{Proof of Theorem \protect\ref{inv}.}

We borrow the idea from \cite{cg2022} and treat the proof in the following
steps.

\noindent \emph{Step 1. }We begin by assuming that there exists two random
variables $U_{\xi }$ and $U_{\xi ^{\prime }}$ which are $\mathcal{F}%
_{t}^{0}\vee \mathcal{G}$-measurable with uniform distribution on $\left[ 0,1%
\right] $, such that $\xi $ and $U_{\xi }$ (resp. $\xi ^{\prime }$ and $%
U_{\xi ^{\prime }}$) are independent. Given any $u\in \mathcal{U}_{ad}^{2}%
\left[ 0,T\right] ,$ from Lemma \ref{inv1}, we are able to find a \ function
$\mathbf{a:}\left[ 0,T\right] \times \Omega \times \mathbb{R}^{n}\times %
\left[ 0,1\right] \rightarrow U$ such that
\begin{equation*}
\left( \xi ,\left( \mathbf{a}_{s}\left( \xi ,U_{\xi }\right) \right) _{0\leq
s\leq T},\left( W_{s}-W_{t}\right) _{0\leq s\leq T}\right) \overset{%
\mathscr{Law}}{=}\left( \xi ,\left( u_{s}\right) _{0\leq s\leq T},\left(
W_{s}-W_{t}\right) _{0\leq s\leq T}\right) .
\end{equation*}%
We notice that $\xi ,$ $\xi ^{\prime }$ have the same law. Therefore,
whenever taking $u^{\prime }=\left( \mathbf{a}_{s}\left( \xi ,U_{\xi
}\right) \right) _{0\leq s\leq T}$ which admits $u^{\prime }\in \mathcal{U}%
_{ad}^{2}\left[ 0,T\right] ,$ we immediately, have
\begin{equation*}
\left( \xi ,\left( \mathbf{a}_{s}\left( \xi ,U_{\xi }\right) \right) _{0\leq
s\leq T},\left( W_{s}-W_{t}\right) _{0\leq s\leq T}\right) \overset{%
\mathscr{Law}}{=}\left( \xi ^{\prime },\left( u_{s}^{\prime }\right) _{0\leq
s\leq T},\left( W_{s}-W_{t}\right) _{0\leq s\leq T}\right)
\end{equation*}%
As we know, under the Assumptions (A1)-(A4), FBSDEs (\ref{fbsde4}) admits a
unique solution. So from inequalities (\ref{estsde1}) and (\ref{best2}), we
imply that
\begin{equation*}
\left( \left( \Upsilon _{s}^{t,\xi ;u}\right) _{0\leq s\leq T},\left(
\mathbf{a}_{s}\left( \xi ,U_{\xi }\right) \right) _{0\leq s\leq T},\right)
\overset{\mathscr{Law}}{=}\left( \left( \Upsilon _{s}^{t,\xi ^{\prime
};u^{\prime }}\right) _{0\leq s\leq T},\left( u_{s}^{\prime }\right) _{0\leq
s\leq T},\right) ,
\end{equation*}%
where $\Upsilon _{s}^{t,\xi ;u}\triangleq \left( X_{s}^{t,\xi
;u},Y_{s}^{t,\xi ;u},Z_{s}^{t,\xi ;u}\right) $ and $\Upsilon _{s}^{t,\xi
^{\prime };u^{\prime }}\triangleq \left( X_{s}^{t,\xi ^{\prime };u^{\prime
}},Y_{s}^{t,\xi ^{\prime };u^{\prime }},Z_{s}^{t,\xi ^{\prime };u^{\prime
}}\right) .$ As a result, from the definition of value function, we get%
\begin{equation}
V\left( t,\xi \right) =\inf_{u\in \Lambda _{t,\xi }}J\left( t,\xi ;u\right)
\leq \inf_{u^{\prime }\in \Lambda _{t,\xi }}J\left( t,\xi ^{\prime
};u^{\prime }\right) =V\left( t,\xi ^{\prime }\right) .  \label{inv3}
\end{equation}%
Now interchanging the roles of $\xi $ and $\xi ^{\prime },$ we derive the
other part, i,e, $V\left( t,\xi \right) \geq V\left( t,\xi ^{\prime }\right)
$.

\noindent \emph{Step 2. }In order to remove the condition ($U_{\xi }$ and $%
U_{\xi ^{\prime }}$ are $\mathcal{F}_{t}^{0}\vee \mathcal{G}$-measurable
with uniform distribution on $\left[ 0,1\right] $, such that $\xi $ and $%
U_{\xi }$ (resp. $\xi ^{\prime }$ and $U_{\xi ^{\prime }}$) are
independent). We first consider the discrete case of $\xi $. Suppose that
there are some $\left\{ x_{1},\ldots ,x_{m}\right\} \subset \mathbb{R}^{n},$
with $x_{i}\neq x_{j},$ if $i\neq j,$ such that
\begin{equation*}
\mathbb{P}_{\xi }=\sum_{i=1}^{m}p_{i}\delta _{x_{i}}.
\end{equation*}%
Instantly, from Lemma \ref{inv2}, there exist two random variables $U_{\xi }$
and $U_{\xi ^{\prime }}$ which are $\mathcal{F}_{t}^{0}\vee \mathcal{G}$%
-measurable with uniform distribution on $\left[ 0,1\right] $, such that $%
\xi $ and $U_{\xi }$ (resp. $\xi ^{\prime }$ and $U_{\xi ^{\prime }}$) are
independent. The desired conclusion is established from Step 1.

\noindent \emph{Step 3. }For the general case, the continuity of value
function will play an important rule as we adopt the technique of
approximation of $\xi $ and $\eta $. Now consider a decomposition of $%
\mathbb{R}^{n}$, namely, $\sum_{i\geq 1}\mathcal{O}_{i}^{k}=\mathbb{R}^{n}$
such that diam$\left( \mathcal{O}_{i}^{k}\right) \leq 2^{-k},$ for each $%
k\geq 1.$ For each $k\in \mathbb{N},$ pick up $\chi _{i}^{k}\in \mathcal{O}%
_{i}^{k}$ and then define%
\begin{equation*}
\check{\xi}_{k}\triangleq \sum_{i\geq 1}\chi _{i}^{k}\mathbb{I}_{\mathcal{O}%
_{i}^{k}}\left( \xi \right) ,\text{ }\check{\eta}_{k}\triangleq \sum_{i\geq
1}\chi _{i}^{k}\mathbb{I}_{\mathcal{O}_{i}^{k}}\left( \eta \right) .
\end{equation*}%
It is easy to check that $\check{\xi}_{k}$, $\check{\eta}_{k}\in L^{2}\left(
\mathcal{F}_{t}^{0}\vee \mathcal{G};\mathbb{R}^{n}\right) $ and $\check{\xi}%
_{k}\rightarrow \xi ,$ $\check{\eta}_{k}\rightarrow \eta $ as $\
k\rightarrow \infty ,$ uniformly w.r.t. $\omega \in \Omega .$ Besides, $%
\check{\xi}_{k}$ and $\check{\eta}_{k}$ have the same law. By virtue of the
diagonal approach, for each $k\in \mathbb{N},$ there exists $K_{k}$ such
that we are able to seek sequences $\xi _{k}$, $\eta _{k},$ denoted by
\begin{equation*}
\xi _{k}\triangleq \sum_{i=1}^{K_{k}}\chi _{i}^{k}\mathbb{I}_{\mathcal{O}%
_{i}^{k}}\left( \xi \right) ,\text{ }\eta _{k}\triangleq
\sum_{i=1}^{K_{k}}\chi _{i}^{k}\mathbb{I}_{\mathcal{O}_{i}^{k}}\left( \eta
\right) ,
\end{equation*}%
which converge respectively to $\xi $ and $\eta .$ Next, from Step 2, we
derive that
\begin{equation*}
V\left( t,\xi _{k}\right) =V\left( t,\eta _{k}\right) ,\text{ }\mathbb{P}%
\text{-a.s. for each }k\in \mathbb{N}.
\end{equation*}%
At last, we get desired result by passing the limit as $k\rightarrow \infty $
by means of the continuity of value function (see Lemma \ref{estvalue}).
\hfill $\Box $

\begin{remark}
In \cite{cg2022}, the authors impose an assumption, so-called $\left(
A_{f,g}\right) \footnote{%
The running function $f$ is locally uniformly continuous in $\left( x,\mu
\right) $ uniformly with respect to $\left( t,u,\nu \right) $. Similarly,
the terminal function $g$ is locally uniformly continuous.}_{\text{cont}}$
is not assumed in our paper since the value function is defined by solution
of BSDE. The assumption (A3)-(A4) can ensure the existence and uniqueness of
the solution.
\end{remark}

\begin{remark}
From now on, the solution of FBSDEs (\ref{fbsde4}) can be written as $\Theta
_{s}^{t,\mu ;u}$ to underline the dependence on the law $\mu =\mathbb{P}%
_{\xi },$ which is important to establish DPP.
\end{remark}

Lemma \ref{inv} shows that the value function $V$ is \emph{law-invariant}.
Therefore, we can define the inverse-lifted function (see Section \ref{sec3}
for the definition) of value function $V.$

\begin{definition}
Given any $\left( t,\mu \right) \in \left[ 0,T\right] \times \mathcal{P}%
_{2}\left( \mathbb{R}^{n}\right) ,$ we define
\begin{equation}
\mathcal{V}\left( t,\mu \right) =V\left( t,\xi \right) ,\text{ }\mathbb{P}%
\text{-a.s.}  \label{d2}
\end{equation}%
where $\xi \in L^{2}\left( \mathcal{F}_{t}^{0}\vee \mathcal{G};\mathbb{R}%
^{n}\right) $ with $\mathbb{P}_{\xi }=\mu .$
\end{definition}

Clearly, under the assumptions (A1)-(A5), the value function $\mathcal{V}%
\left( t,\mu \right) $ is well-defined and a priori it is a bounded $%
\mathcal{F}_{t}^{0}\vee \mathcal{G}$-measurable random variable. But it
turns out that $\mathcal{V}\left( t,\mu \right) $ is even deterministic. We
have the following result.

\begin{proposition}
Suppose that the assumptions \emph{(A1)-(A5)} are in force. The value
function $\mathcal{V}\left( t,\mu \right) $ is a deterministic function.
\end{proposition}

The main idea for this proof is to construct a sequence of admissible
controls used by Peng in \cite{pyfw97} with help of Lemma \ref{vdet}. So we
omit it.

It is necessary to point out that the value function defined in (\ref{d2})
can be taken over the subset $\mathcal{U}_{t}^{2}$ of elements in $\mathcal{U%
}_{ad}^{2}\left[ 0,T\right] ,$ which are independent of $\mathcal{F}^{0}$
under $\mathbb{P}^{0}.$

\begin{lemma}
\label{dpp2}Under the Assumptions \emph{(A1)-(A5)}, for any $\left( t,\mu
\right) \in \left[ 0,T\right] \times \mathcal{P}_{2}\left( \mathbb{R}%
^{n}\right) ,$ we have
\begin{equation}
\mathcal{V}\left( t,\mu \right) =\inf_{u\in \mathcal{U}_{t}^{2}}J\left(
t,\mu ;u\right) .\text{ }  \label{dpp1}
\end{equation}
\end{lemma}

\paragraph{Proof.}

\emph{Step 1.} Given any $\mathcal{F}^{0}$-stopping time $\tau \in \mathbb{T}%
_{t,T}^{0},$ according to the flow property (\ref{flow}), backward semigroup
(\ref{flow3})\ and noting $\rho _{s}^{t,\mu ;u}$ is $\mathcal{F}_{\tau }^{0}$
-measurable, we get for $\mathbb{P}^{0}$-a.s. $\omega ^{0}\in \Omega ^{0},$%
\begin{eqnarray}
J\left( \tau ,\rho _{\tau }^{t,\mu ;u};u^{\tau }\right) \left( \omega
^{0}\right) &=&Y_{\tau }^{t,\mu ;u}\left( \omega ^{0}\right)  \notag \\
&=&\mathbb{G}_{\tau ,T}^{t,\mu ;u}\left[ \Phi ^{0}\left( \rho _{T}^{t,\mu
;u}\right) \right] \left( \omega ^{0}\right)  \notag \\
&=&\left. \mathbb{G}_{r,T}^{r,\pi ;\upsilon }\left[ \Phi ^{0}\left( \rho
_{T}^{r,\pi ;\upsilon }\right) \right] \left( \omega ^{0}\right) \right\vert
_{r=\tau \left( \omega ^{0}\right) ,\pi =\rho _{\tau \left( \omega
^{0}\right) }^{t,\mu ;u}\left( \omega ^{0}\right) ,\upsilon =u^{r,\omega
^{0}}}  \notag \\
&=&\left. J\left( r,\pi ;\upsilon \right) \right\vert _{r=\tau \left( \omega
^{0}\right) ,\pi =\rho _{\tau \left( \omega ^{0}\right) }^{t,\mu ;u}\left(
\omega ^{0}\right) ,\upsilon =u^{r,\omega ^{0}}},\text{ }\mathbb{P}^{0}\text{%
-a.s.}  \label{dp1}
\end{eqnarray}%
in which the process $u^{r,\omega ^{0}}\in \mathcal{U}_{r}^{2}$ is
independent of $\mathcal{F}_{r}^{0}$, which implies that $X_{s}^{t,\vartheta
;u}$ (due to $\mathcal{G}$ \textquotedblleft rich enough\textquotedblright ,
there exists a random variable $\vartheta \in L^{2}\left( \mathcal{G};%
\mathbb{R}^{n}\right) $ such that $\mathcal{L}\left( \vartheta \right) =\pi $%
) is independent of $\mathcal{F}_{r}^{0}$, and consequently $\rho
_{s}^{r,\pi ;\upsilon }$ is also independent of $\mathcal{F}_{r}^{0}$ for $%
r\leq s$. Therefore, we get the relation (\ref{dp1}) by virtue of Lemma \ref%
{deter}.

\emph{Step 2. }For the sake of simplicity, we put
\begin{equation*}
\mathcal{\tilde{V}}\left( t,\mu \right) =\inf_{u\in \mathcal{U}%
_{t}^{2}}J\left( t,\mu ;u\right) .
\end{equation*}%
From the fact $\mathcal{U}_{r}^{2}\subset \mathcal{U}_{ad}^{2}\left[ 0,T%
\right] $, it yields that $\mathcal{V}\left( t,\mu \right) \leq \mathcal{%
\tilde{V}}\left( t,\mu \right) .$ To obtain the opposite relationship, we
let $\tau =t$ in (\ref{dp1})$,$ then get
\begin{equation}
J\left( t,\mu ;u^{t,\omega ^{0}}\right) =\left. J\left( t,\mu ;u\right)
\right\vert _{t,\mu ,u^{t,\omega ^{0}}},\text{ }\mathbb{P}^{0}\text{-a.s.}
\end{equation}%
Immediately,
\begin{equation*}
\int_{\Omega ^{0}}J\left( t,\mu ;u^{t,\omega ^{0}}\right) \mathbb{P}%
^{0}\left( \mathrm{d}\omega ^{0}\right) =J\left( t,\mu ;u\right) .
\end{equation*}%
Besides, noting however $u^{t}\in \mathcal{U}_{t}^{2},$ we have
\begin{equation*}
J\left( t,\mu ;u^{t,\omega ^{0}}\right) \geq \mathcal{\tilde{V}}\left( t,\mu
\right) ,
\end{equation*}%
from which we get the desired result since the arbitrary of $u\in \mathcal{U}%
_{ad}^{2}\left[ 0,T\right] .$ \hfill $\Box $

We can get some properties of the value function $\mathcal{V}\left( t,\mu
\right) $ which is an immediate consequence of its definitions (\ref{d2}), (%
\ref{d1}) and the Lipschitz property (\ref{best2}) of the cost functional.

\begin{lemma}
\label{p1}Under the Assumptions \emph{(A1)-(A5)}, there exists a constant $%
C>0$ such that, for all $t\in \left[ 0,T\right] $, $\mu ,\mu ^{\prime }\in
\mathcal{P}_{2}\left( \mathbb{R}^{n}\right) $%
\begin{eqnarray}
\left\vert \mathcal{V}\left( t,\mu \right) -\mathcal{V}\left( t,\mu ^{\prime
}\right) \right\vert &\leq &C\mathcal{W}_{2}\left( \mu ,\mu ^{\prime
}\right) ,  \label{p11} \\
\left\vert \mathcal{V}\left( t,\mu \right) \right\vert &\leq &C\left[ 1+%
\mathcal{W}_{2}\left( \mu ,\delta _{0}\right) \right] ,  \label{p22}
\end{eqnarray}%
where $\delta _{0}$ is the Dirac measure with mass at $0\in \mathbb{R}^{n}.$
\end{lemma}


\begin{lemma}
\label{p2}Under the Assumptions \emph{(A1)-(A5)}, for any $\left( t,\xi
\right) \in \left[ 0,T\right] \times L^{2}\left( \mathcal{G};\mathbb{R}%
^{n}\right) $ with $\mathbb{P}_{\xi }=\mu $ and all $u\in \mathcal{U}%
_{ad}^{2}\left[ 0,T\right] $ and $\mathcal{F}^{0}$-stopping time $\tau \in
\mathbb{T}_{t,T}^{0},$ we have%
\begin{equation}
\mathcal{V}\left( t,\mu \right) \leq \mathbb{G}_{t,\tau }^{t,\mu ;u}\left[
\mathcal{V}\left( \tau ,\rho _{\tau }^{t,\mu ;u}\right) \right] ,
\label{dpp3}
\end{equation}%
Meanwhile, for any $\varepsilon >0,$ there exists an admissible control $%
u^{\varepsilon }\in \mathcal{U}_{ad}^{2}\left[ 0,T\right] $ for all $\tau
\in \mathbb{T}_{t,T}^{0},$ such that
\begin{equation}
\mathcal{V}\left( t,\mu \right) \geq \mathbb{G}_{t,\tau }^{t,\mu
;u^{\varepsilon }}\left[ \mathcal{V}\left( \tau ,\rho _{\tau }^{t,\mu
;u^{\varepsilon }}\right) \right] -\varepsilon .  \label{dpp4}
\end{equation}
\end{lemma}

\begin{remark}
The inequality in (\ref{dpp4}) indicates that, for any given $\varepsilon
>0, $ there exists a near-optimal control $u^{\varepsilon }.$ It has nice
structure and large availability as well as flexibility. So it is possible
to choose among them which are easier for analysis and implementation (see
\cite{Z2018}). Beyond that, this result will be employed to prove the
viscosity super-solution property of the value function.
\end{remark}

\paragraph{Proof.}

From (\ref{dp1}) and (\ref{dpp1}), we have for all $\tau \in \mathbb{T}%
_{t,T}^{0}$
\begin{eqnarray}
J\left( t,\mu ;u\right) &=&\mathbb{G}_{t,\tau }^{t,\mu ;u}\left[ J\left(
\tau ,\rho _{\tau }^{t,\mu ;u};u^{\tau }\right) \right]  \notag \\
&\geq &\mathbb{G}_{t,\tau }^{t,\mu ;u}\left[ \mathcal{V}\left( \tau ,\rho
_{\tau }^{t,\mu ;u}\right) \right] .  \label{r1}
\end{eqnarray}%
The arbitrary of $u$ implies
\begin{equation}
\mathcal{V}\left( t,\mu \right) \geq \inf_{u\in \mathcal{U}_{ad}^{2}\left[
0,T\right] }\sup_{\tau \in \mathbb{T}_{t,T}^{0}}\mathbb{G}_{t,\tau }^{t,\mu
;u}\left[ \mathcal{V}\left( \tau ,\rho _{\tau }^{t,\mu ;u}\right) \right] .
\label{pp1}
\end{equation}%
Now for any $\varepsilon >0$ and $\omega ^{0}\in \Omega ^{0},$ by Lemma \ref%
{dp1}, there exists certain $u^{\varepsilon }\in \mathcal{U}_{\tau }^{2}$
such that
\begin{equation}
\mathcal{V}\left( \tau \left( \omega ^{0}\right) ,\rho _{\tau \left( \omega
^{0}\right) }^{t,\mu ;u}\right) +\varepsilon \geq J\left( \tau \left( \omega
^{0}\right) ,\rho _{\tau \left( \omega ^{0}\right) }^{t,\mu ;u}\left( \omega
^{0}\right) ;u^{\varepsilon }\left( \omega ^{0}\right) \right) ,\text{ }%
\mathbb{P}^{0}\text{-a.s.}  \label{r2}
\end{equation}%
The measurable selection arguments \cite{W1980} combining to Lemma \ref{p1}
yields that the map $\omega ^{0}\in \Omega ^{0}\rightarrow u^{\varepsilon
}\left( \omega ^{0}\right) $ is measurable. We now define a shifted control
process based on $u^{\varepsilon }$%
\begin{equation*}
\tilde{u}_{s}^{\varepsilon }\left( \omega ^{0}\right) =u_{s}\left( \omega
^{0}\right) \mathbb{I}_{s<\tau \left( \omega ^{0}\right) }+u^{\varepsilon
}\left( \omega ^{0}\right) \mathbb{I}_{s\geq \tau \left( \omega ^{0}\right)
},\text{ }s\in \left[ 0,T\right] .
\end{equation*}%
Clearly, $\rho _{s}^{t,\mu ;\tilde{u}^{\varepsilon }}=\rho _{s}^{t,\mu ;u}$
whenever $s\leq \tau $ and $\rho ^{t,\mu ;u}$ has continuous trajectories.
Notice that $\mathcal{U}_{ad}^{2}\left[ 0,T\right] $ is a separable metric
space, by Lemma 2.1 in \cite{ST2002}, the control process $\tilde{u}%
_{s}^{\varepsilon }$ is $\mathcal{F}_{s}^{0}$-adapted, and thus $\tilde{u}%
^{\varepsilon }\in \mathcal{U}_{ad}^{2}\left[ 0,T\right] .$ Therefore, (\ref%
{r2}) can be expressed as
\begin{equation}
\mathcal{V}\left( \tau ,\rho _{\tau }^{t,\mu ;u}\right) +\varepsilon \geq
J\left( \tau ,\rho _{\tau }^{t,\mu ;u};\tilde{u}^{\varepsilon }\right) ,%
\text{ }\mathbb{P}^{0}\text{-a.s.}  \label{r3}
\end{equation}%
Then%
\begin{eqnarray*}
\mathcal{V}\left( t,\mu \right) &\leq &J\left( t,\mu ;\tilde{u}^{\varepsilon
}\right) \\
&=&\mathbb{G}_{t,\tau }^{t,\mu ;u}\left[ J\left( \tau ,\rho _{\tau }^{t,\mu
;u};\left( \tilde{u}^{\varepsilon }\right) ^{\tau }\right) \right] \\
&\leq &\mathbb{G}_{t,\tau }^{t,\mu ;u}\left[ \mathcal{V}\left( \tau ,\rho
_{\tau }^{t,\mu ;u}\right) +\varepsilon \right] ,
\end{eqnarray*}%
where the third inequality holds by using the comparison theorem (see Lemma %
\ref{combsde}). More preciously, set
\begin{equation*}
\zeta ^{1}=J\left( \tau ,\rho _{\tau }^{t,\mu ;u};\left( \tilde{u}%
^{\varepsilon }\right) ^{\tau }\right) ,\text{ }\zeta ^{2}=\mathcal{V}\left(
\tau ,\rho _{\tau }^{t,\mu ;u}\right) +\varepsilon ,
\end{equation*}%
in Lemma \ref{combsde}, from the fact $\zeta ^{1}\leq \zeta ^{2},$ $\mathbb{P%
}^{0}$-a.s., we get the desired result.

Now since $\tau ,u$ and $\varepsilon $ are arbitrary, it follows that
\begin{equation}
\mathcal{V}\left( t,\mu \right) \leq \inf_{u\in \mathcal{U}_{ad}^{2}\left[
0,T\right] }\inf_{\tau \in \mathbb{T}_{t,T}^{0}}\mathbb{G}_{t,\tau }^{t,\mu
;u}\left[ \mathcal{V}\left( \tau ,\rho _{\tau }^{t,\mu ;u}\right)
+\varepsilon \right] .  \label{pp2}
\end{equation}%
At last, the inequalities (\ref{pp1}) and (\ref{pp2}) yield%
\begin{eqnarray*}
\mathcal{V}\left( t,\mu \right) &\leq &\inf_{u\in \mathcal{U}_{ad}^{2}\left[
0,T\right] }\inf_{\tau \in \mathbb{T}_{t,T}^{0}}\mathbb{G}_{t,\tau }^{t,\mu
;u}\left[ \mathcal{V}\left( \tau ,\rho _{\tau }^{t,\mu ;u}\right) \right] \\
&=&\inf_{u\in \mathcal{U}_{ad}^{2}\left[ 0,T\right] }\sup_{\tau \in \mathbb{T%
}_{t,T}^{0}}\mathbb{G}_{t,\tau }^{t,\mu ;u}\left[ \mathcal{V}\left( \tau
,\rho _{\tau }^{t,\mu ;u}\right) \right] ,
\end{eqnarray*}%
which support the DPP property. We complete the proof. \hfill $\Box $

\begin{remark}
Note that in \cite{CH2019}, the authors consider the initial condition $\xi
\in L^{2}\left( \mathcal{F}_{t}\vee \mathcal{G};\mathbb{R}^{n}\right) .$
Whenever the coefficients in \cite{CH2019} are only dependent on the state
and control rather than the law of state, moreover $\mathcal{G=}\left\{
\varnothing ,\Omega \right\} $, nevertheless the initial state is still a
random variable$,$ then $\xi \in L^{2}\left( \mathcal{F}_{t};\mathbb{R}%
^{n}\right) $. Observe that the cost functional defined in \cite{CH2019} and
\cite{PW2017} by means of taking the absolute expectation which instantly
implies that the value function is deterministic. Thereafter, one can
retrieve the standard dynamic principle in the non-McKean--Vlasov situation
(see Proposition 4.1 in \cite{CH2019}). However, under the similar
assumptions, Problem MVS defined in our paper actually becomes a kind of
standard stochastic recursive optimal control problem (cf. \cite{pyfw97} for
more details). From Lemma \ref{rdet}, one can easily derive that $\left.
Y_{s}^{t,\xi ;u}\right\vert _{s=t}$ is random, which is completely different
from the framework studied in \cite{CH2019} and \cite{PW2017}.
\end{remark}

\begin{lemma}
\label{estv}Under the Assumptions \emph{(A1)-(A5)}, there exists a positive
constant $C$ such that%
\begin{equation*}
\left\vert \mathcal{V}\left( t,\mu \right) -\mathcal{V}\left( t+\delta ,\mu
\right) \right\vert \leq C\left( 1+\mathcal{W}_{2}\left( \mu ,\delta
_{0}\right) \right) \delta ^{\frac{1}{2}},
\end{equation*}%
where $\mu \in \mathcal{P}_{2}\left( \mathbb{R}^{n}\right) ,$ $t\in \left[
0,T\right) $ and $\delta >0$ with $0<\delta +t\leq T.$
\end{lemma}

The proof can be found in Appendix.

\section{Hamilton-Jacobi-Bellman Equation\label{sec3}}

After obtaining the DPP, as usual, we will derive the corresponding Bellman
equation with help of It\^{o}'s formula in Wasserstein space and the
viscosity PDE characterization. The value function will be studied on the
Hilbert space or the Wasserstein space of probability measures via the
associated inverse-lifted identification. To this end, let us briefly recall
the following notions of derivatives w.r.t. a probability measure, as
introduced by Lions (see \cite{Lions2012}).

The following notations are mainly taken from \cite{PW2017, CCD2015}. Let $%
\vartheta $ be a real-valued function defined on $\mathcal{P}_{2}\left(
\mathbb{R}^{n}\right) .$ We say $\upsilon $ is the \textit{lifted} version
of $\vartheta $, namely, the function defined on $L^{2}\left( \Omega ,%
\mathcal{G},\mathbb{P};\mathbb{R}^{n}\right) $ by $\upsilon \left( \xi
\right) =\vartheta \left( \mathbb{P}_{\xi }\right) .$ Conversely, for any
function $\upsilon $ defined on $L^{2}\left( \Omega ,\mathcal{G},\mathbb{P};%
\mathbb{R}^{n}\right) $, we call inverse-lifted function of $\upsilon $ the
function $\vartheta $ defined on $\mathcal{P}_{2}\left( \mathbb{R}%
^{n}\right) $ by $\upsilon \left( \xi \right) =\vartheta \left( \mathbb{\mu }%
\right) $ for $\mu =\mathcal{L}\left( \xi \right) .$ Clearly, $\vartheta $
exists if and only if $\upsilon $ depends only on the distribution of $\xi $
for any $\xi \in L^{2}\left( \Omega ,\mathcal{G},\mathbb{P};\mathbb{R}%
^{n}\right) .$ In the sequel, we set $\upsilon =\vartheta $ often for
simplicity.

In the literatures, $\vartheta $ is said to be $C^{1}$ on $\mathcal{P}%
_{2}\left( \mathbb{R}^{n}\right) $ if the lift $\upsilon $ is Fr\'{e}chet
differentiable with continuous derivatives on $L^{2}\left( \Omega ,\mathcal{G%
},\mathbb{P};\mathbb{R}^{n}\right) .$ By the Riesz representation theorem,
the Fr\'{e}chet derivative $\left[ D\upsilon \right] (\xi )$, regarded as an
element $D\upsilon (\xi )$ of $L^{2}\left( \Omega ,\mathcal{G},\mathbb{P};%
\mathbb{R}^{n}\right) ,$ can be expressed as:
\begin{equation*}
\left[ D\upsilon \right] (\xi )\left( \zeta \right) =\mathbb{E}^{1}\left[
\left\langle D\upsilon (\xi ),\zeta \right\rangle \right]
\end{equation*}%
with $D\upsilon (\xi )=\partial _{\mu }\vartheta \left( \mathcal{L}\left(
\xi \right) \right) \left( \xi \right) $ for $\zeta \in L^{2}\left( \Omega ,%
\mathcal{G},\mathbb{P};\mathbb{R}^{n}\right) ,$ in which $\partial _{\mu
}\vartheta \left( \mathcal{L}\left( \xi \right) \right) :\mathbb{R}%
^{n}\rightarrow \mathbb{R}^{n}$ is called the derivative of $\vartheta $ at $%
\mu =\mathcal{L}\left( \xi \right) .$ Note that $\mathcal{G}$ is assumed to
be \textquotedblleft rich enough\textquotedblright . Therefore, $\mu \in
\mathcal{P}_{2}\left( \mathbb{R}^{n}\right) =\left\{ \left. \mathcal{L}%
\left( \xi \right) \right\vert \xi \in L^{2}\left( \Omega ,\mathcal{G},%
\mathbb{P};\mathbb{R}^{n}\right) \right\} .$ We say that $\vartheta $ is
fully $C^{2}$ if it is $C^{1}$, and one can find, for any $\mu \in \mathcal{P%
}_{2}\left( \mathbb{R}^{n}\right) $, a continuous version of the mapping $%
x\in \mathbb{R}^{n}\rightarrow \partial _{\mu }\vartheta \left( \mathcal{\mu
}\right) \left( x\right) $, such that the mapping $\left( x,\mu \right) \in
\mathbb{R}^{n}\times \mathcal{P}_{2}\left( \mathbb{R}^{n}\right) \rightarrow
\partial _{\mu }\vartheta \left( \mathcal{\mu }\right) \left( x\right) $ is
continuous at any point $\left( x,\mu \right) $ such that $x\in $Supp$\left(
\mu \right) ,$ and (a) for any $\mu \in \mathcal{P}_{2}\left( \mathbb{R}%
^{n}\right) $, the mapping $x\in \mathbb{R}^{n}\rightarrow \partial _{\mu
}\vartheta \left( \mathcal{\mu }\right) \left( x\right) $ is differentiable
in the standard sense, with a gradient denoted by $\partial _{x}\partial
_{\mu }\vartheta \left( \mathcal{\mu }\right) \left( x\right) \in \mathbb{R}%
^{n\times n}$, and s.t. the mapping $\left( x,\mu \right) \in \mathbb{R}%
^{n}\times \mathcal{P}_{2}\left( \mathbb{R}^{n}\right) \rightarrow \partial
_{x}\partial _{\mu }\vartheta \left( \mathcal{\mu }\right) \left( x\right) $
is continuous; (b) for any fixed$\ x\in \mathbb{R}^{n}$, the mapping $\mu
\in \mathcal{P}_{2}\left( \mathbb{R}^{n}\right) \rightarrow \partial _{\mu
}\vartheta \left( \mathcal{\mu }\right) \left( x\right) $ is differentiable
in the above lifted sense. Its derivative, interpreted thus as a mapping $%
\tilde{x}\in \mathbb{R}^{n}\rightarrow \partial _{\mu }\left[ \partial _{\mu
}\vartheta \left( \mathcal{\mu }\right) \left( x\right) \right] \left(
\tilde{x}\right) \in \mathbb{R}^{n\times n}$ in $L_{\mu }^{2}\left( \mathbb{R%
}^{n\times n}\right) $, is denoted by $\tilde{x}\in \mathbb{R}%
^{n}\rightarrow \partial _{\mu }^{2}\vartheta \left( \mathcal{\mu }\right)
\left( x,\tilde{x}\right) $, and s.t. the mapping $(\mu ,x,\tilde{x})\in
\mathcal{P}_{2}\left( \mathbb{R}^{n}\right) \times \mathbb{R}^{n}\times
\mathbb{R}^{n}\rightarrow \partial _{\mu }^{2}\vartheta \left( \mathcal{\mu }%
\right) \left( x,\tilde{x}\right) $ is continuous.

Let $C_{b}^{2}\left( \mathcal{P}_{2}\left( \mathbb{R}^{n}\right) \right) $
denotes the set of functions if it is partially $C^{2},\partial _{x}\partial
_{\mu }\vartheta \left( \mathcal{\mu }\right) \in L_{\mu }^{\infty }\left(
\mathbb{R}^{n\times n}\right) ,$ and for any compact set $\mathcal{K}$ of $%
\mathcal{P}_{2}\left( \mathbb{R}^{n}\right) $, it admits
\begin{equation*}
\sup_{\mu \in \mathcal{K}}\left[ \int_{\mathcal{K}}\left\vert \partial _{\mu
}\vartheta \left( \mathcal{\mu }\right) \left( x\right) \right\vert ^{2}\mu
\left( \mathrm{d}x\right) +\left\Vert \partial _{x}\partial _{\mu }\vartheta
\left( \mathcal{\mu }\right) \right\Vert _{\infty }+\left\Vert \partial
_{\mu }^{2}\vartheta \left( \mathcal{\mu }\right) \right\Vert _{\infty }%
\right] <\infty .
\end{equation*}%
We now ready to state the It\^{o} formula (see \cite{PW2017, CCD2015}).
Consider the following type of It\^{o} process%
\begin{equation}
x_{t}=x_{0}+\int_{0}^{t}b_{r}\mathrm{d}r+\int_{0}^{t}\sigma _{r}\mathrm{d}%
W_{r},\text{ }0\leq t\leq T,  \label{itosde}
\end{equation}%
where $x_{0}$ is independent of $W,$ while $b,\sigma $ are are progressively
measurable processes w.r.t. the natural filtration $\mathcal{F}$ generated
by $\left( x_{0},W\right) $, and satisfying the following condition:
\begin{equation*}
\mathbb{E}^{1}\left[ \int_{0}^{T}\left( \left\vert b_{t}\right\vert
^{2}+\left\vert \sigma _{t}\right\vert ^{2}\right) \mathrm{d}t\right]
<\infty .
\end{equation*}%
As mentioned before, $\mathbb{P}_{x_{t}}^{W}$ denotes the conditional law of
$x_{t}$, $0\leq t\leq T$, given the $\sigma $-field $\mathcal{F}$ generated
by the whole filtration of $W$, namely, $\mathbb{P}_{x_{t}}^{W}\left( \omega
^{0}\right) =\mathbb{P}_{x_{t}\left( \omega ^{0},\cdot \right) }^{1}$. Then,
for any $\vartheta \in C_{b}^{1}\left( \mathcal{P}_{2}\left( \mathbb{R}%
^{n}\right) \right) ,$ we have%
\begin{eqnarray}
\vartheta \left( \mathbb{P}_{x_{t}}^{W}\right) &=&\vartheta \left( \mathbb{P}%
_{x_{0}}\right) +\int_{0}^{t}\Big \{\mathbb{E}^{1}\big [\left\langle
\partial _{\mu }\vartheta \left( \mathbb{P}_{x_{s}}^{W}\right) \left(
x_{s}\right) ,b_{s}\right\rangle  \notag \\
&&+\frac{1}{2}\text{tr}\left( \partial _{x}\partial _{\mu }\vartheta \left(
\mathbb{P}_{x_{s}}^{W}\right) \left( x_{s}\right) \sigma _{s}\sigma
_{s}^{\top }\right) \big ]  \notag \\
&&+\mathbb{E}^{1}\Big [\mathbb{\tilde{E}}^{1}\big [\frac{1}{2}\text{tr}%
\left( \partial _{\mu }^{2}\vartheta \left( \mathbb{P}_{x_{s}}^{W}\right)
\left( x_{s},\tilde{x}_{s}\right) \sigma _{s}\tilde{\sigma}_{s}^{\top
}\right) \big ]\Big \}\mathrm{d}s  \notag \\
&&+\int_{0}^{t}\mathbb{E}^{1}\big [\partial _{\mu }\vartheta \left( \mathbb{P%
}_{x_{s}}^{W}\right) \left( x_{s}\right) ^{\top }\sigma _{s}\big ]\mathrm{d}%
W_{s},  \label{ito1}
\end{eqnarray}%
where $\tilde{x},\tilde{\sigma}$ are copies of $x$ and $\sigma $ on anther
probability space $\tilde{\Omega}=\big (\Omega ^{0}\times \tilde{\Omega}^{1},%
\mathcal{F}^{0}\otimes \mathcal{\tilde{F}}^{1},\mathbb{P}^{0}\otimes \mathbb{%
\tilde{P}}^{1}\big ).$ It will be useful to apply It\^{o} formula for the
lifted function on $L^{2}\left( \tilde{\Omega}^{1},\mathcal{\tilde{F}}^{1},%
\mathbb{\tilde{P}}^{1};\mathbb{R}^{n}\right) .$ To this end, we postulate
that $\upsilon \in C^{2}\left( L^{2}\left( \Omega ^{1},\mathcal{G},\mathbb{P}%
^{1};\mathbb{R}^{n}\right) \right) ,$ then we have the following connection
between derivatives in the Wasserstein space $\mathcal{P}_{2}\left( \mathbb{R%
}^{n}\right) $ and in the Hilbert space $L^{2}\left( \Omega ^{1},\mathcal{G},%
\mathbb{P}^{1};\mathbb{R}^{n}\right) $:%
\begin{eqnarray}
D^{2}\upsilon \left( \xi \right) \left[ \eta ,\eta \right] &=&\mathbb{E}^{1}%
\big [\mathbb{\tilde{E}}^{1}\left[ \text{tr}\left( \partial _{\mu
}^{2}\vartheta \left( \mathcal{L}\left( \xi \right) \right) \left( \xi ,%
\tilde{\xi}\right) \eta \tilde{\eta}^{\top }\right) \right]  \notag \\
&&+\mathbb{\tilde{E}}^{1}\left[ \text{tr}\left( \partial _{x}\partial _{\mu
}\vartheta \left( \mathcal{L}\left( \xi \right) \right) \left( \xi \right)
\eta \eta ^{\top }\right) \right] ,  \label{re1}
\end{eqnarray}%
where $\xi ,\eta \in L^{2}\left( \Omega ^{1},\mathcal{G},\mathbb{P}^{1};%
\mathbb{R}^{n}\right) $ and $\left( \tilde{\xi},\tilde{\eta}\right) $ is a
copy of $\left( \xi ,\eta \right) $ on another Polish and atomless space $%
\left( \tilde{\Omega}^{1},\mathcal{\tilde{F}}^{1},\mathbb{\tilde{P}}%
^{1}\right) .$ Let $\tilde{x}_{0},\tilde{b},\tilde{\sigma}$ be copies of $%
x_{0},b,\sigma $ on $\left( \Omega ^{0}\times \tilde{\Omega}^{1},\mathcal{F}%
^{0}\otimes \mathcal{G},\mathbb{P}^{0}\otimes \mathbb{\tilde{P}}^{1}\right) $
and now consider the following SDE, a copy of (\ref{itosde})%
\begin{equation}
\tilde{x}_{t}=\tilde{x}_{0}+\int_{0}^{t}\tilde{b}_{r}\mathrm{d}r+\int_{0}^{t}%
\tilde{\sigma}_{r}\mathrm{d}W_{r},\text{ }0\leq t\leq T.  \label{itosde2}
\end{equation}

Based on (\ref{itosde2}), we define the process as $\breve{\phi}_{t}\left(
\omega _{0}\right) =\tilde{\phi}_{t}\left( \omega _{0},\cdot \right) ,$ for $%
0\leq t\leq T,\omega _{0}\in \Omega ^{0}.$ Clearly, if $\breve{\phi}=\tilde{b%
},$ $\tilde{\sigma},$ and $\tilde{x},$ $\breve{\phi}_{t}$ is $\mathcal{F}%
^{0} $-progressive, valued in $L^{2}\left( \tilde{\Omega}^{1},\mathcal{%
\tilde{F}}^{1},\mathbb{\tilde{P}}^{1};\mathbb{R}^{n}\right) ,L^{2}\left(
\tilde{\Omega}^{1},\mathcal{\tilde{F}}^{1},\mathbb{\tilde{P}}^{1};\mathbb{R}%
^{n\times d}\right) $ and $L^{2}\left( \tilde{\Omega}^{1},\mathcal{\tilde{F}}%
^{1},\mathbb{\tilde{P}}^{1};\mathbb{R}^{n}\right) ,$ respectively$.$ Based
on (\ref{re1}) and the assumption $\upsilon \in C^{2}\left( \tilde{\Omega}%
^{1},\mathcal{\tilde{F}}^{1},\mathbb{\tilde{P}}^{1};\mathbb{R}^{n}\right) ,$
It\^{o} formula (\ref{ito1}) can be rewritten as (see Proposition 6.3 in
\cite{cd14})
\begin{eqnarray}
\upsilon \left( \breve{x}_{t}\right) &=&\upsilon \left( \breve{x}_{0}\right)
+\int_{0}^{t}\mathbb{E}^{1}\Big [\left\langle D\upsilon \left( \breve{x}%
_{s}\right) ,\breve{b}_{s}\right\rangle +\frac{1}{2}\left\langle
D^{2}\upsilon \left( \breve{x}_{s}\right) \breve{\sigma}_{s},\breve{\sigma}%
_{s}\right\rangle \Big ]\mathrm{d}s  \notag \\
&&+\int_{0}^{t}\mathbb{E}^{1}\left[ D\upsilon \left( \breve{x}_{s}\right)
^{\top }\breve{\sigma}_{s}\right] \mathrm{d}W_{s}.  \label{ito2}
\end{eqnarray}

\begin{remark}
Note that It\^{o} formula (\ref{ito1}) holds even if the lift is not twice
continuously Fr\'{e}chet differentiable as shown in \cite{cdll19}.
\end{remark}

Now suppose that $\mathcal{V}\left( t,\mu \right) $ is sufficiently smooth.
For simplicity of notations, we set%
\begin{eqnarray*}
b_{s}^{t,\mu ,u} &=&b\left( X_{s}^{t,\mu ;u},\mathbb{P}_{X_{s}^{t,\mu
;u}}^{W},u_{s}\right) ,\text{ }\sigma _{s}^{t,\mu ,u}=\sigma \left(
X_{s}^{t,\mu ;u},\mathbb{P}_{X_{s}^{t,\mu ;u}}^{W},u_{s}\right) , \\
\tilde{b}_{s}^{t,\mu ,u} &=&\tilde{b}\left( X_{s}^{t,\mu ;u},\mathbb{P}%
_{X_{s}^{t,\mu ;u}}^{W},u_{s}\right) ,\text{ }\tilde{\sigma}_{s}^{t,\mu ,u}=%
\tilde{\sigma}\left( X_{s}^{t,\mu ;u},\mathbb{P}_{X_{s}^{t,\mu
;u}}^{W},u_{s}\right) .
\end{eqnarray*}%
By virtue of Lemma \ref{p2} (DPP), we have%
\begin{equation*}
\mathcal{V}\left( t,\mu \right) =\inf_{u\in \mathcal{U}_{ad}^{2}\left[ 0,T%
\right] }\mathbb{G}_{t,\tau }^{t,\mu ;u}\left[ \mathcal{V}\left( t+\delta
,\rho _{t+\delta }^{t,\mu ;u}\right) \right] ,
\end{equation*}%
namely%
\begin{eqnarray*}
0 &=&\inf_{u\in \mathcal{U}_{ad}^{2}\left[ 0,T\right] }\Bigg [\mathcal{V}%
\left( t+\delta ,\rho _{t+\delta }^{t,\mu ;u}\right) -\mathcal{V}\left(
t,\mu \right) \\
&&+\int_{t}^{t+\delta }f^{0}\left( \rho _{s}^{t,\mu ;u},\Theta _{s}^{t,\mu
;u},\mathbb{P}_{\Theta _{s}^{t,\mu ;u}},u_{s}\right) \mathrm{d}%
s-\int_{t}^{t+\delta }Z_{s}^{t,\mu ;u}\mathrm{d}W_{s}\Bigg ].
\end{eqnarray*}%
By applying It\^{o} formula (\ref{ito1}) to (\ref{dpp3})-(\ref{dpp4}), we
obtain the HJB equation: all $0\leq t<t+\delta <T$%
\begin{eqnarray}
&&0=\inf_{u\in \mathcal{U}_{ad}^{2}\left[ 0,T\right] }\Bigg \{%
\int_{t}^{t+\delta }f^{0}\left( \rho _{s}^{t,\mu ;u},\Theta _{s}^{t,\mu ;u},%
\mathbb{P}_{\Theta _{s}^{t,\mu ;u}},u_{s}\right) \mathrm{d}%
s-\int_{t}^{t+\delta }Z_{s}^{t,\mu ;u}\mathrm{d}W_{s}  \notag \\
&&+\int_{t}^{t+\delta }\Bigg [\mathbb{E}^{1}\Big [\left\langle \partial
_{\mu }\mathcal{V}\left( s,\mathbb{P}_{X_{s}^{t,\mu ;u}}^{W}\right) \left(
X_{s}^{t,\mu ;u}\right) ,b_{s}^{t,\mu ,u}\right\rangle  \notag \\
&&+\frac{1}{2}\text{tr}\left( \partial _{x}\partial _{\mu }\mathcal{V}\left(
s,\mathbb{P}_{X_{s}^{t,\mu ;u}}^{W}\right) \left( X_{s}^{t,\mu ;u}\right)
\sigma _{s}^{t,\mu ,u}\left( \sigma _{s}^{t,\mu ,u}\right) ^{\top }\right) %
\Big ]  \notag \\
&&+\mathbb{E}^{1}\left[ \mathbb{\tilde{E}}^{1}\left[ \frac{1}{2}\text{tr}%
\left( \partial _{\mu }^{2}\mathcal{V}\left( s,\mathbb{P}_{X_{s}^{t,\mu
;u}}^{W}\right) \left( X_{s}^{t,\mu ;u},\tilde{X}_{s}^{t,\mu ;u}\right)
\sigma _{s}^{t,\mu ,u}\left( \tilde{\sigma}_{s}^{t,\mu ,u}\right) ^{\top
}\right) \right] \right] \Bigg ]\mathrm{d}s  \notag \\
&&+\int_{t}^{t+\delta }\mathbb{E}^{1}\left[ \partial _{\mu }\mathcal{V}%
\left( s,\mathbb{P}_{X_{s}^{t,\mu ;u}}^{W}\right) \left( X_{s}^{t,\mu
;u}\right) ^{\top }\sigma _{s}^{t,\mu ,u}\right] \mathrm{d}W_{s}\Bigg \},
\label{hjb1}
\end{eqnarray}%
from which, we derive that
\begin{equation*}
\left\{
\begin{array}{l}
Y_{s}^{t,\mu ;\bar{u}}=\mathcal{V}\left( s,\mathbb{P}_{X_{s}^{t,\mu
;u}}^{W}\right) ,\text{ }\mathbb{P}\text{-a.s. }t\leq s\leq T, \\
Z_{s}^{t,\mu ;\bar{u}}=\mathbb{E}^{1}\left[ \partial _{\mu }\mathcal{V}%
\left( t,\mathbb{P}_{X_{s}^{t,\mu ;u}}^{W}\right) \left( X_{s}^{t,\mu
;u}\right) ^{\top }\sigma _{s}^{t,\mu ,u}\right] ,\text{ }\mathbb{P}\text{%
-a.s. }t\leq s\leq T.%
\end{array}%
\right.
\end{equation*}%
Arranging the above derivations and combing the fact that $\mathcal{V}$
depends on $\xi $ only through distribution $\mathcal{L}\left( \xi \right) $%
, we announce that HJB equation associated to the value function of
stochastic recursive McKean-Vlasov control problem appearing in the
following type:%
\begin{equation}
\left\{
\begin{array}{rcl}
0 & = & -\partial _{t}\mathcal{V}\left( t,\mu \right) \\
&  & -\inf_{u\in U}\Bigg [f^{0}\Big (\mu ,\mathcal{V}\left( t,\mathbb{\mu }%
\right) ,\mathbb{\mu }\left( \partial _{\mu }\mathcal{V}\left( t,\mathbb{\mu
}\right) x^{\top }\sigma \left( x,\mu ,u\right) \right) , \\
&  & \mathbb{P}_{\left( \mathcal{V}\left( t,\mathbb{\mu }\right) ,\partial
_{\mu }\mathcal{V}\left( t,\mathbb{\mu }\right) x^{\top }\sigma \left( x,\mu
,u\right) \right) },u\Big ) \\
&  & +\mu \left( \mathscr{L}^{u}\mathcal{V}\left( t,\mathbb{\mu }\right)
\right) +\mu \otimes \mu \left( \mathscr{M}^{u}\mathcal{V}\left( t,\mathbb{%
\mu }\right) \right) \Bigg ], \\
\mathcal{V}\left( T,\mathbb{\mu }\right) & = & \Phi ^{0}\left( \mu \right) ,%
\text{ }\left( t,\mu \right) \in \left[ 0,T\right] \times \mathcal{P}%
_{2}\left( \mathbb{R}^{n}\right) ,%
\end{array}%
\right.  \label{hjb}
\end{equation}%
where for $\varphi \in C^{2}\left( \mathcal{P}_{2}\left( \mathbb{R}%
^{n}\right) \right) ,$ the functions $\mathscr{L}^{u}\mathcal{\varphi }%
\left( \mathbb{\mu }\right) :\mathbb{R}^{n}\rightarrow \mathbb{R}$ and $%
\mathscr{M}^{u}\mathcal{\varphi }\left( \mathbb{\mu }\right) :\mathbb{R}%
^{n}\times \mathbb{R}^{n}\rightarrow \mathbb{R}$ are defined by
\begin{eqnarray*}
&&\mathscr{L}^{u}\mathcal{\varphi }\left( \mu \right) \left( x\right)
=\left\langle \partial _{\mu }\mathcal{\varphi }\left( \mu \right) \left(
x\right) ,b\left( x,\mu ,u\right) \right\rangle +\frac{1}{2}\text{tr}\left(
\partial _{x}\partial _{\mu }\mathcal{\varphi }\left( \mu \right) \left(
x\right) \sigma \sigma ^{\top }\left( x,\mu ,u\right) \right) , \\
&&\mathscr{M}^{u}\mathcal{\varphi }\left( \mu \right) \left( x,\tilde{x}%
\right) =\frac{1}{2}\text{tr}\left( \partial _{\mu }^{2}\mathcal{\varphi }%
\left( \mu \right) \left( x,\tilde{x}\right) \sigma \left( x,\mu ,u\right)
\sigma ^{\top }\left( \tilde{x},\mu ,u\right) \right) .
\end{eqnarray*}

\begin{remark}
\label{com} Comparing with the existing results, our paper makes the
following contributions: 1) There appears the term $\mu \otimes \mu \left( %
\mathscr{M}^{u}\mathcal{V}\left( t,\mu \right) \right) $ in (\ref{hjb})
since we apply the It\^{o} formula (\ref{ito1}) with stochastic flow which
is slightly different from the PDE in \cite{Li2018, blp09}. 2) Our
methodology can deal with the situation whenever the generator $f$ can
depends on the value function and its gradient, which actually extends the
framework \cite{CH2019}, \cite{PW2017} and \cite{DPT2022} to a more general
case.
\end{remark}

By means of the lifting identification, HJB equations (\ref{hjb}) can be
expressed in $\left[ 0,T\right] \times L^{2}\left( \Omega ^{1},\mathcal{G},%
\mathbb{P}^{1};\mathbb{R}^{n}\right) .$ To keep the compatibility of
symbols, we still denote $\mathcal{V}\left( t,\mathbb{\xi }\right) =\mathcal{%
V}\left( t,\mathcal{L}\left( \xi \right) \right) $ combing the fact that $%
\mathcal{V}$ depends on $\xi $ only through distribution again. From the
relation (\ref{re1}), it follows, from Proposition 5.6 in \cite{cg2022},
that the HJB equations (\ref{hjb}) can be reported as

\begin{equation}
\left\{
\begin{array}{rcl}
0 & = & -\partial _{t}\mathcal{V}-\mathcal{H}\left( \xi ,\mathcal{V}\left( t,%
\mathbb{\xi }\right) ,D\mathcal{V}\left( t,\mathbb{\xi }\right) ,D^{2}%
\mathcal{V}\left( t,\mathbb{\xi }\right) \right) , \\
&  & \text{ }\left( t,\xi \right) \in \left[ 0,T\right] \times L^{2}\left(
\Omega ^{1},\mathcal{G},\mathbb{P}^{1};\mathbb{R}^{n}\right) , \\
\mathcal{V}\left( T,\mathbb{\xi }\right) & = & \mathbb{E}^{1}\left[ \Phi
\left( \xi ,\mathcal{L}\left( \xi \right) \right) \right] ,\text{ }\xi \in
L^{2}\left( \Omega ^{1},\mathcal{G},\mathbb{P}^{1};\mathbb{R}^{n}\right) ,%
\end{array}%
\right.  \label{hjb2}
\end{equation}%
where $\mathcal{H}:L^{2}\left( \Omega ^{1},\mathcal{G},\mathbb{P}^{1};%
\mathbb{R}^{n}\right) \times L^{2}\left( \Omega ^{1},\mathcal{G},\mathbb{P}%
^{1};\mathbb{R}\right) \times L^{2}\left( \Omega ^{1},\mathcal{G},\mathbb{P}%
^{1};\mathbb{R}^{n}\right) \times L^{2}\left( \Omega ^{1},\mathcal{G},%
\mathbb{P}^{1};\mathbb{S}^{n}\right) \rightarrow \mathbb{R}$ is defined as
\begin{eqnarray}
&&\mathcal{H}\left( \xi ,p,q,Q\right)  \notag \\
&=&\mathbb{E}^{1}\Bigg \{\inf_{u\in U}\Bigg [f\Big (\xi ,\mathcal{L}\left(
\xi \right) ,p,\mathbb{E}^{1}\left[ q^{\top }\sigma \left( \xi ,\mathcal{L}%
\left( \xi \right) ,u\right) \right] ,\mathcal{L}_{\left( p,\mathbb{E}^{1}%
\left[ q^{\top }\sigma \left( \xi ,\mathcal{L}\left( \xi \right) ,u\right) %
\right] \right) },u\Big )  \notag \\
&&+\left\langle q,b\left( \xi ,\mathcal{L}\left( \xi \right) ,u\right)
\right\rangle +\frac{1}{2}Q\left( \sigma \left( \xi ,\mathcal{L}\left( \xi
\right) ,u\right) \right) \sigma \left( \xi ,\mathcal{L}\left( \xi \right)
,u\right) \Bigg ]\Bigg \}.  \label{ham}
\end{eqnarray}

\begin{remark}
It follows, from Proposition 5.6 in \cite{cg2022}, that the HJB equations (%
\ref{hjb2}) can be reported as%
\begin{equation}
\left\{
\begin{array}{rcl}
0 & = & -\partial _{t}\mathcal{V}-\mathcal{\tilde{H}}\left( \xi ,\mathcal{V}%
\left( t,\mathbb{\xi }\right) ,D\mathcal{V}\left( t,\mathbb{\xi }\right)
,D^{2}\mathcal{V}\left( t,\mathbb{\xi }\right) \right) , \\
&  & \text{ }\left( t,\xi \right) \in \left[ 0,T\right] \times L^{2}\left(
\Omega ^{1},\mathcal{G},\mathbb{P}^{1};\mathbb{R}^{n}\right) , \\
\mathcal{V}\left( T,\mathbb{\xi }\right) & = & \mathbb{E}^{1}\left[ \Phi
\left( \xi ,\mathcal{L}\left( \xi \right) \right) \right] ,\text{ }\xi \in
L^{2}\left( \Omega ^{1},\mathcal{G},\mathbb{P}^{1};\mathbb{R}^{n}\right) ,%
\end{array}%
\right.  \label{ham2}
\end{equation}%
where $\mathcal{\tilde{H}}:L^{2}\left( \Omega ^{1},\mathcal{G},\mathbb{P}%
^{1};\mathbb{R}^{n}\right) \times L^{2}\left( \Omega ^{1},\mathcal{G},%
\mathbb{P}^{1};\mathbb{R}\right) \times L^{2}\left( \Omega ^{1},\mathcal{G},%
\mathbb{P}^{1};\mathbb{R}^{n}\right) \times L^{2}\left( \Omega ^{1},\mathcal{%
G},\mathbb{P}^{1};\mathbb{S}^{n}\right) \rightarrow \mathbb{R}$ is defined
as
\begin{eqnarray*}
&&\mathcal{\tilde{H}}\left( \xi ,p,q,Q\right) \\
&=&\mathbb{E}^{1}\Bigg \{\inf_{u\in U}\Bigg [f\Big (\xi ,\mathcal{L}\left(
\xi \right) ,p,\mathbb{E}^{1}\left[ q^{\top }\sigma \left( \xi ,\mathcal{L}%
\left( \xi \right) ,u\right) \right] ,\mathcal{L}_{\left( p,\mathbb{E}^{1}%
\left[ q^{\top }\sigma \left( \xi ,\mathcal{L}\left( \xi \right) ,u\right) %
\right] \right) },u\Big ) \\
&&+\left\langle q,b\left( \xi ,\mathcal{L}\left( \xi \right) ,u\right)
\right\rangle +\frac{1}{2}Q\left( \sigma \left( \xi ,\mathcal{L}\left( \xi
\right) ,u\right) \right) \sigma \left( \xi ,\mathcal{L}\left( \xi \right)
,u\right) \Bigg ]\Bigg \}.
\end{eqnarray*}
\end{remark}

As we have well-known, in the stochastic optimal control theory, the value
function is a solution to the corresponding HJB equation whenever it admits
classical solutions (smooth enough to the order of derivatives involved in
the equation)$\footnote{%
We refer to the paper \cite{GS2015} for existence result of smooth solution
to the HJB equation (\ref{hjb2}) on small time horizon.}$. Unfortunately,
this is not necessarily the case even for some very simple situations.

Normally, whenever the diffusion is possibly degenerate, the HJB equation
may in general have no classical solutions. To remove this obstacle,
Crandall and Lions introduced the so-called viscosity solutions in the early
1980s (cf. \cite{CIL}). Out of question, this idea makes the theory a
powerful tool in studying the optimal control problems.

Similarly, we will employ a notion of viscosity solutions via the lifting
identification in the Hilbert space $L^{2}\left( \Omega ^{1},\mathcal{G},%
\mathbb{P}^{1};\mathbb{R}^{n}\right) $ instead of the Wasserstein space $%
\mathcal{P}_{2}\left( \mathbb{R}^{n}\right) $.

For simplicity, we define, for $\left( \xi ,r,p,P\right) \in L^{2}\left(
\Omega ^{1},\mathcal{G},\mathbb{P}^{1};\mathbb{R}^{n}\right) \times
L^{2}\left( \Omega ^{1},\mathcal{G},\mathbb{P}^{1};\mathbb{R}\right) \times
L^{2}\left( \Omega ^{1},\mathcal{G},\mathbb{P}^{1};\mathbb{R}^{n}\right)
\times L^{2}\left( \Omega ^{1},\mathcal{G},\mathbb{P}^{1};\mathbb{S}%
^{n}\right) ,$
\begin{equation}
F\left( \xi ,r,p,P\right) =-\mathcal{H}\left( \xi ,r,p,P\right) .  \label{F}
\end{equation}%
Then the HJB equation (\ref{hjb2}) can be expressed as%
\begin{equation}
\left\{
\begin{array}{rcl}
0 & = & -\partial _{t}\mathcal{V}\left( t,\mathbb{\xi }\right) -F\left( \xi ,%
\mathcal{V}\left( t,\mathbb{\xi }\right) ,D\mathcal{V}\left( t,\mathbb{\xi }%
\right) ,D^{2}\mathcal{V}\left( t,\mathbb{\xi }\right) \right) , \\
&  & \text{ }\left( t,\xi \right) \in \left[ 0,T\right] \times L^{2}\left(
\Omega ^{1},\mathcal{G},\mathbb{P}^{1};\mathbb{R}^{n}\right) , \\
\mathcal{V}\left( T,\mathbb{\xi }\right) & = & \mathbb{E}^{1}\left[ \Phi
\left( \xi ,\mathcal{L}\left( \xi \right) \right) \right] ,\text{ }\xi \in
L^{2}\left( \Omega ^{1},\mathcal{G},\mathbb{P}^{1};\mathbb{R}^{n}\right) .%
\end{array}%
\right.  \label{hjb3}
\end{equation}

The following definition of viscosity solution for infinite dimension is
similar as Definition 3.40 in \cite{FGS2015} for time-dependent situation
since HJB equation (\ref{hjb3}) does not contain the dissipative operator.

\begin{definition}
\label{dvis1} Let $\mathscr{V}\left( t,\mu \right) \in C\left( \left[ 0,T%
\right] \times \mathcal{P}_{2}\left( \mathbb{R}^{n}\right) \right) $ and $%
\tilde{\mathscr{V}}$ be lifted version of $\mathscr{V}.$ For every $\varphi
\in C^{1,2}\left( \left[ 0,T\right] \times L^{2}\left( \Omega ^{1},\mathcal{G%
},\mathbb{P}^{1};\mathbb{R}^{n}\right) \right) \footnote{%
Here $C^{1,2}\left( \left[ 0,T\right] \times L^{2}\left( \Omega ^{1},%
\mathcal{G},\mathbb{P}^{1};\mathbb{R}^{n}\right) \right) $ denotes the set
of real-valued continuous functions on $\left[ 0,T\right] \times L^{2}\left(
\Omega ^{1},\mathcal{G},\mathbb{P}^{1};\mathbb{R}^{n}\right) $ which are
continuously differentiable in $t\in \left[ 0,T\right) $, and twice
continuously Fr\'{e}chet differentiable on $L^{2}\left( \Omega ^{1},\mathcal{%
G},\mathbb{P}^{1};\mathbb{R}^{n}\right) .$},$

\noindent (1) for each local maximum point $\left( t_{0},\xi _{0}\right) \in %
\left[ 0,T\right] \times L^{2}\left( \Omega ^{1},\mathcal{G},\mathbb{P}^{1};%
\mathbb{R}^{n}\right) $ of $\tilde{\mathscr{V}}-\varphi $ in the interior of
$\left[ 0,T\right] \times L^{2}\left( \Omega ^{1},\mathcal{G},\mathbb{P}^{1};%
\mathbb{R}^{n}\right) ,$ we have%
\begin{equation}
-\frac{\partial }{\partial t}\varphi +F\left( \xi _{0},\varphi ,D\varphi
,D^{2}\varphi \right) \left( t_{0},\xi _{0}\right) \leq 0  \label{vis1}
\end{equation}%
and
\begin{equation}
\tilde{\mathscr{V}}\left( T,\xi _{0}\right) -\Phi ^{0}\left( \mathcal{L}%
\left( \xi _{0}\right) \right) \leq 0,  \label{TT1}
\end{equation}%
i.e., $\mathscr{V}$ is a sub-solution to HJB equation \emph{(\ref{hjb3})};

\noindent (2) for each local minimum point $\left( t_{0},\xi _{0}\right) \in %
\left[ 0,T\right] \times L^{2}\left( \Omega ^{1},\mathcal{G},\mathbb{P}^{1};%
\mathbb{R}^{n}\right) $ of $\tilde{\mathscr{V}}-\varphi $ in the interior of
$\left[ 0,T\right] \times L^{2}\left( \Omega ^{1},\mathcal{G},\mathbb{P}^{1};%
\mathbb{R}^{n}\right) ,$ we have%
\begin{equation}
-\frac{\partial }{\partial t}\varphi +F\left( \xi _{0},\varphi ,D\varphi
,D^{2}\varphi \right) \left( t_{0},\xi _{0}\right) \geq 0  \label{vis2}
\end{equation}%
and
\begin{equation}
\tilde{\mathscr{V}}\left( T,\xi _{0}\right) -\Phi ^{0}\left( \mathcal{L}%
\left( \xi _{0}\right) \right) \geq 0,  \label{TT2}
\end{equation}%
i.e., $\mathscr{V}$ is a super-solution to HJB equation \emph{(\ref{hjb3})}.

\noindent (3) $\mathscr{V}\left( t,\mu \right) \in C\left( \left[ 0,T\right]
\times \mathcal{P}_{2}\left( \mathbb{R}^{n}\right) \right) $ is said to be a
viscosity solution of \emph{(\ref{hjb3})} if it is both a viscosity sub- and
super-solution.
\end{definition}

We now focus on the case, $\nu \in \mathcal{P}_{2}\left( \mathbb{R}\right) ,$
namely, when $f$ does not depend on the law of $Z$ due to Lemma \ref{combsde}%
. We add the following assumption:

\begin{enumerate}
\item[\textbf{(A6)}] For all $\eta _{1},\eta _{2}\in L^{2}\left( \mathcal{F}%
_{t}^{0}\vee \mathcal{G};\mathbb{R}\right) $ and all $x\in \mathbb{R}^{n},$ $%
y\in \mathbb{R},$ $z\in \mathbb{R}^{d},$ $\mu \in \mathcal{P}_{2}\left(
\mathbb{R}^{n}\right) $ and $u\in U$, there exists a positive constant $C$
such that
\begin{equation*}
f\left( x,\mu ,y,z,\mathbb{P}_{\eta _{1}},u\right) -f\left( x,\mu ,y,z,%
\mathbb{P}_{\eta _{2}},u\right) \leq C\left( \mathbb{E}\left[ \left( \left(
\eta _{1}-\eta _{2}\right) ^{+}\right) ^{2}\right] \right) ^{\frac{1}{2}}.
\end{equation*}
\end{enumerate}

Note that the Hamiltonian function $\mathcal{H}$ define in (\ref{ham})
depends not only on $\mathcal{V}$ but its derivative. In order to prove the
uniqueness of solution to HJB equation (\ref{hjb2}), we need the following
hypothesis:

\begin{enumerate}
\item[\textbf{(A7)}] There exists a $\kappa \geq 0$ such that for every $%
\left( x,\mu ,z,\nu ,u\right) \in \mathbb{R}^{n}\times \mathcal{P}_{2}\left(
\mathbb{R}^{n}\right) \times \mathbb{R}^{d}\times \mathcal{P}_{2}\left(
\mathbb{R}\right) \times U$ and $r,$ $s\in \mathbb{R}$ with $r\geq s,$
\begin{equation*}
f\left( x,\mu ,r,z,\nu ,u\right) -f\left( x,\mu ,s,z,\nu ,u\right) \leq
\kappa \left( r-s\right) .
\end{equation*}
\end{enumerate}

\begin{remark}
In stochastic control theory in infinite dimension (see \cite{FGS2015}), we
say that $F$ defined in (\ref{F}) is \emph{proper} whenever $\kappa >0 $ in
(A7), which actually plays an important role to prove the uniqueness of
viscosity solution in Hilbert space (see Theorem 3.50 in \cite{FGS2015}).
\end{remark}

We are at position to assert the main result:

\begin{theorem}
\label{m1}Assume that the Assumptions \emph{(A1)-(A7)} are in force. Then
the value function defined in (\emph{\ref{d2}}) is a unique viscosity
solution to HJB equation (\emph{\ref{hjb2}}).
\end{theorem}

\paragraph{Proof.}

We will employ the DPP of Lemma \ref{p2} for the value function now regarded
as a function on $\left[ 0,T\right] \times L^{2}\left( \Omega ^{1},\mathcal{G%
},\mathbb{P}^{1};\mathbb{R}^{n}\right) .$ To this end, denoting $\mathscr{X}%
_{s}^{t,\xi ;u}\left( \omega ^{0}\right) =X_{s}^{t,\xi ;u}\left( \omega
^{0},\cdot \right) $ for $t\leq s\leq T.$ Apparently, $\left\{ \mathscr{X}%
_{s}^{t,\xi ;u}\right\} _{t\leq s\leq T}$ is $\mathcal{F}^{0}$-adapted
process, valued in $L^{2}\left( \Omega ^{1},\mathcal{G},\mathbb{P}^{1};%
\mathbb{R}^{n}\right) $ and $\rho _{s}^{t,\mu ;u}=\mathbb{P}_{\mathscr{X}%
_{s}^{t,\xi ;u}}^{1}$ for $t\leq s\leq T$ with $\mu =\mathcal{L}\left( \xi
\right) $. Immediately, we derive that the value function satisfying $%
\mathcal{V}\left( s,\mathscr{X}_{s}^{t,\xi ;u}\right) =\mathcal{V}\left(
s,\rho _{s}^{t,\mu ;u}\right) ,$ namely, the lifted valued function on $%
\left[ 0,T\right] \times L^{2}\left( \Omega ^{1},\mathcal{G},\mathbb{P}^{1};%
\mathbb{R}^{n}\right) $ equals to the value function on $\left[ 0,T\right]
\times \mathcal{P}_{2}\left( \mathbb{R}^{n}\right) .$ Therefore, DPP in
Lemma \ref{p2} can be written as, for any $\left( t,\xi \right) \in \left[
0,T\right] \times L^{2}\left( \Omega ^{1},\mathcal{G},\mathbb{P}^{1};\mathbb{%
R}^{n}\right) $ with $\mu =\mathcal{L}\left( \xi \right) ,$%
\begin{eqnarray*}
\mathcal{V}\left( t,\xi \right) &=&\inf_{u\in \mathcal{U}_{ad}^{2}\left[ 0,T%
\right] }\inf_{\tau \in \mathbb{T}_{t,T}^{0}}\mathbb{G}_{t,\tau }^{t,\xi ;u}%
\left[ \mathcal{V}\left( \tau ,\mathscr{X}_{\tau }^{t,\xi ;u}\right) \right]
\\
&=&\inf_{u\in \mathcal{U}_{ad}^{2}\left[ 0,T\right] }\sup_{\tau \in \mathbb{T%
}_{t,T}^{0}}\mathbb{G}_{t,\tau }^{t,\xi ;u}\left[ \mathcal{V}\left( \tau ,%
\mathscr{X}_{\tau }^{t,\xi ;u}\right) \right] ,
\end{eqnarray*}

From Lemma \ref{p1} and Lemma \ref{estv}, we know that the value function is
continuous on $\left[ 0,T\right] \times L^{2}\left( \Omega ^{1},\mathcal{G},%
\mathbb{P}^{1};\mathbb{R}^{n}\right) .$ Particularly, at the terminal time,
it requires that $\mathcal{V}\left( T,\xi \right) =\mathbb{E}^{1}\left[ \Phi
\left( \xi ,\mathcal{L}\left( \xi \right) \right) \right] .$ For the proof
of this theorem, we take the following three steps.

\emph{Step 1.}\textit{\ }For arbitrarily chosen but fixed $\varphi \in
C_{l,b}^{2}\left( \left[ 0,T\right] \times L^{2}\left( \Omega ^{1},\mathcal{G%
},\mathbb{P}^{1};\mathbb{R}^{n}\right) \right) \footnote{%
Here $C_{l,b}^{2}\left( \left[ 0,T\right] \times L^{2}\left( \Omega ^{1},%
\mathcal{G},\mathbb{P}^{1};\mathbb{R}^{n}\right) \right) $ denotes the set
of class $C^{2}$ (Lions's sense) whose partial derivatives of order less
than or equal to $2$ are bounded (and hence the function itself grows at
most linearly at infinity).},$ we define
\begin{eqnarray*}
\mathscr{F}\left( t,\xi ,y,z,u\right) &=&\partial _{t}\varphi \left( t,\xi
\right) +\mathbb{E}^{1}\Bigg \{\left\langle D\varphi \left( t,\xi \right)
,b\left( \xi ,\mathcal{L}\left( \xi \right) ,u\right) \right\rangle \\
&&+\frac{1}{2}D^{2}\varphi \left( t,\xi \right) \left( \sigma \left( \xi ,%
\mathcal{L}\left( \xi \right) ,u\right) \right) \sigma \left( \xi ,\mathcal{L%
}\left( \xi \right) ,u\right) \\
&&+f\Big (\xi ,\mathcal{L}\left( \xi \right) ,y+\varphi \left( t,\xi \right)
,z+\mathbb{E}^{1}\left[ D\varphi \left( t,\xi \right) ^{\top }\sigma \left(
\xi ,\mathcal{L}\left( \xi \right) ,u\right) \right] , \\
&&\mathcal{L}_{\left( y+\varphi \left( t,\xi \right) \right) },u\Big )\Bigg
\}
\end{eqnarray*}%
$\left( t,\xi ,y,z,u\right) \in \left[ 0,T\right] \times L^{2}\left( \Omega
^{1},\mathcal{G},\mathbb{P}^{1};\mathbb{R}^{n}\right) \times \mathbb{R}%
\times \mathbb{R}^{d}\times U,$ and we introduce the following BSDE defined
on the interval $\left[ t,t+\delta \right] $ for $0<\delta \leq T-t,$%
\begin{equation}
\left\{
\begin{array}{rcl}
-\mathrm{d}\mathscr{Y}_{s}^{1,u} & = & \mathscr{F}\left( s,\mathscr{X}%
_{s}^{t,\xi ;u},\mathscr{Y}_{s}^{1,u},\mathscr{Z}_{s}^{1,u},u_{s}\right)
\mathrm{d}s-\mathscr{Z}_{s}^{1,u}\mathrm{d}W_{s},\text{ }t\leq s\leq
t+\delta , \\
\mathscr{Y}_{t+\delta }^{1,u} & = & 0,%
\end{array}%
\right.  \label{bs1}
\end{equation}%
where $u\in \mathcal{U}_{ad}^{2}\left[ 0,T\right] $ and the process $%
\mathscr{X}_{s}^{t,\xi ;u}$ is defined by $X_{s}^{t,\xi ;u}\left( \omega
^{0},\cdot \right) .$ It is easy to check that $\mathscr{F}\left( t,\xi
,y,z,u\right) $ fulfills assumptions (H2.1)-(H2.2) in \cite{pyfw97}. Thus
BSDE (\ref{bs1}) admits a unique solution. As a matter of fact, we are able
to establish the relationship between $\mathscr{Y}_{s}^{1,u}$ and $%
\mathscr{X}_{s}^{t,\xi ;u}$ by virtue of It\^{o}'s formula as follows:%
\begin{equation}
\mathscr{Y}_{s}^{1,u}=\mathbb{G}_{s,t+\delta }^{t,\xi ;u}\left[ \varphi
\left( t+\delta ,\mathscr{X}_{t+\delta }^{t,\xi ;u}\right) \right] -\varphi
\left( s,\mathscr{X}_{s}^{t,\xi ;u}\right)  \label{rebs1}
\end{equation}%
Recall that $\mathbb{G}_{s,t+\delta }^{t,\xi ;u}\left[ \varphi \left(
t+\delta ,\mathscr{X}_{t+\delta }^{t,\xi ;u}\right) \right] $ is
corresponding to the solution of the following BSDE:%
\begin{equation}
\left\{
\begin{array}{rcl}
-\mathrm{d}\mathscr{Y}_{s}^{u} & = & f\left( s,\mathscr{X}_{s}^{t,\xi ;u},%
\mathscr{Y}_{s}^{u},\mathscr{Z}_{s}^{u},u_{s}\right) \mathrm{d}s-\mathscr{Z}%
_{s}^{u}\mathrm{d}W_{s}, \\
\mathscr{Y}_{T}^{u} & = & \varphi \left( t+\delta ,\mathscr{X}_{t+\delta
}^{t,\xi ;u}\right)%
\end{array}%
\right.  \label{bs2}
\end{equation}%
by
\begin{equation*}
\mathbb{G}_{s,t+\delta }^{t,\xi ;u}\left[ \varphi \left( t+\delta ,%
\mathscr{X}_{t+\delta }^{t,\xi ;u}\right) \right] =\mathscr{Y}_{s}^{u},\text{
}t\leq s\leq t+\delta .
\end{equation*}%
The relation (\ref{rebs1}) can be verified by applying It\^{o}'s formula to $%
\varphi \left( s,\mathscr{X}_{s}^{t,\xi ;u}\right) .$ Afterward, we compare
the coefficients of the stochastic differentials of $\mathscr{Y}%
_{s}^{u}-\varphi \left( s,\mathscr{X}_{s}^{t,\xi ;u}\right) $ and $%
\mathscr{Y}_{s}^{1,u},$ meanwhile we notice that at the terminal time $%
t+\delta ,$ the equality
\begin{equation*}
\mathscr{Y}_{t+\delta }^{u}-\varphi \left( s,\mathscr{X}_{t+\delta }^{t,\xi
;u}\right) =0=\mathscr{Y}_{t+\delta }^{1,u}
\end{equation*}%
holds.

\emph{Step 2. }Consider the following BSDE
\begin{equation}
\left\{
\begin{array}{rcl}
-\mathrm{d}\mathscr{Y}_{s}^{2,u} & = & \mathscr{F}\left( s,\xi ,\mathscr{Y}%
_{s}^{2,u},\mathscr{Z}_{s}^{2,u},u_{s}\right) \mathrm{d}s-\mathscr{Z}%
_{s}^{2,u}\mathrm{d}W_{s},\text{ }t\leq s\leq t+\delta , \\
\mathscr{Y}_{t+\delta }^{2,u} & = & 0,%
\end{array}%
\right.  \label{bs3}
\end{equation}%
where the driving process $\mathscr{X}_{s}^{t,\xi ;u}$ in (\ref{bs1}) are
replaced by the initial state $\xi .$

\begin{lemma}
\label{estvis}For any $u\in \mathcal{U}_{ad}^{2}\left[ 0,T\right] ,$ there
exists a constant $C$ independent of $u\in \mathcal{U}_{ad}^{2}\left[ 0,T%
\right] $ such that%
\begin{equation*}
\left\vert \mathscr{Y}_{t}^{1,u}-\mathscr{Y}_{t}^{2,u}\right\vert \leq
C\delta \rho _{1}\left( \delta \right) ,
\end{equation*}%
where $\rho _{1}\left( \delta \right) :\mathbb{R}_{+}\rightarrow \mathbb{R}.$
When $\delta \rightarrow 0,$ it has $\rho _{1}\left( \delta \right)
\rightarrow 0$ and $\rho _{1}\left( \delta \right) $ is independent of
control $u\in \mathcal{U}_{ad}^{2}\left[ 0,T\right] .$
\end{lemma}

The proof can be seen in the Appendix.

\begin{lemma}
\label{ode}Consider the following ordinary differential equation with law
(ODEL for short) :%
\begin{equation}
\left\{
\begin{array}{rcl}
-\dot{\mathscr{Y}}_{s} & = & \Upsilon \left( s,\xi ,\mathscr{Y}_{s},0\right)
, \\
\mathscr{Y}_{t+\delta } & = & 0,%
\end{array}%
\right.  \label{ode1}
\end{equation}%
where the function $\Upsilon $ is defined by
\begin{equation}
\Upsilon \left( t,\xi ,y,z\right) \triangleq \inf_{u\in U}\mathscr{F}\left(
t,\xi ,y,z,u\right) ,  \label{deode}
\end{equation}%
in which $\left( t,\xi ,y,z\right) \in \left[ 0,T\right] \times L^{2}\left(
\Omega ^{1},\mathcal{G},\mathbb{P}^{1};\mathbb{R}^{n}\right) \times \mathbb{R%
}\times \mathbb{R}^{d}.$ Then%
\begin{equation*}
\underset{u\in \mathcal{U}_{ad}^{2}\left[ 0,T\right] }{\text{essinf}}%
\mathscr{Y}_{t}^{2,u}=\mathscr{Y}_{t},\text{ a.s.}
\end{equation*}
\end{lemma}

The proof also can be seen in the Appendix.

Apparently, $\Upsilon \left( t,\xi ,y,z\right) $ is Lipschitz in $\left(
y,z\right) $, uniformly with respect to $\left( t,\xi \right) $. This
guarantees the existence and uniqueness of ODEL (\ref{ode1}). \newline
Consider the BSDE%
\begin{equation}
\left\{
\begin{array}{rcl}
-\dot{\mathscr{Y}}_{s} & = & \Upsilon \left( s,\xi ,\mathscr{Y}_{s},%
\mathscr{Z}_{s}\right) \mathrm{d}s-\mathscr{Z}_{s}\mathrm{d}W_{s} \\
\mathscr{Y}_{t+\delta } & = & 0,\text{ }t\leq s\leq t+\delta .%
\end{array}%
\right.  \label{bsdec1}
\end{equation}%
From the definition of $\Upsilon $ (see (\ref{deode})) we know
\begin{equation*}
\Upsilon \left( t,\xi ,y,z\right) \leq \mathscr{F}\left( t,\xi
,y,z,u_{s}\right) ,\text{ for all }u\in \mathcal{U}_{ad}^{2}\left[ 0,T\right]
.
\end{equation*}%
From Lemma \ref{combsde} (comparison theorem), we have%
\begin{equation}
\mathscr{Y}_{t}\leq \underset{u\in \mathcal{U}_{ad}^{2}\left[ 0,T\right] }{%
\text{essinf}}\mathscr{Y}_{t}^{2,u}.  \label{invereq}
\end{equation}%
Observe that $\Upsilon $ is a deterministic function of $\left( t,\xi
,y,z\right) $ (see the definition of $F$), therefore $\left( \mathscr{Y}_{s},%
\mathscr{Z}_{s}\right) \equiv $ $\left( \mathscr{Y}_{s},0\right) $ where $%
\mathscr{Y}$ is the solution of the ODEL (\ref{ode1}). Now to prove the
reverse inequality in (\ref{invereq}), we introduce the following ODEL%
\begin{equation}
\left\{
\begin{array}{rcl}
-\dot{\mathscr{Y}_{s}^{2,\upsilon }} & = & \mathscr{F}\left( s,\xi ,%
\mathscr{Y}_{s}^{2,\upsilon },0,\upsilon _{s}\right) , \\
\mathscr{Y}_{t+\delta }^{2,\upsilon } & = & 0,\text{ }t\leq s\leq t+\delta ,%
\end{array}%
\right.
\end{equation}%
where $\upsilon \in \mathscr{U}\left[ 0,T\right] $ (Here $\mathscr{U}\left[
0,T\right] $ denotes the set of all deterministic admissible controls).
Immediately, one can show%
\begin{equation*}
\mathscr{Y}_{t}=\inf_{\upsilon \in \mathscr{U}\left[ 0,T\right] }\mathscr{Y}%
_{t}^{2,\upsilon }.
\end{equation*}

\emph{Step 3. }We will prove that $\mathcal{V}$ is a viscosity solution of (%
\ref{hjb2}). At the beginning, we have $\mathcal{V}\left( T,\xi \right) =%
\mathbb{E}^{1}\left[ \Phi \left( \xi ,\mathcal{L}\left( \xi \right) \right) %
\right] .$ We only prove that $\mathcal{V}$ is a viscosity sub-solution
since the proof is analogous for the viscosity super-solution. For this we
suppose that $\varphi \in C^{1,2}\left( \left[ 0,T\right] \times L^{2}\left(
\Omega ^{1},\mathcal{G},\mathbb{P}^{1};\mathbb{R}^{n}\right) \right) ,$ and $%
\left( t,\xi \right) \in \left[ 0,T\right] \times L^{2}\left( \Omega ^{1},%
\mathcal{G},\mathbb{P}^{1};\mathbb{R}^{n}\right) $ such that $\mathcal{V}%
-\varphi $ attains a maximum at $\left( t,\xi \right) .$ Without loss of
generality, we may also suppose that $\mathcal{V}\left( t,\xi \right)
=\varphi \left( t,\xi \right) $. Then, by virtue of the DPP (see Lemma \ref%
{p2} and Lemma \ref{inv2}),
\begin{equation*}
\varphi \left( t,\xi \right) =\mathcal{V}\left( t,\xi \right) =\inf_{u\in
\mathcal{U}_{ad}^{2}\left[ 0,T\right] }\mathbb{G}_{t,t+\delta }^{t,\mu ;u}%
\left[ \mathcal{V}\left( t+\delta ,X_{t+\delta }^{t,\mu ;u}\right) \right] ,%
\text{ }0\leq \delta \leq T-t,
\end{equation*}%
and from $\mathcal{V}\leq \varphi $ and the monotonicity property of $%
\mathbb{G}_{t,t+\delta }^{t,\mu ;u}\left[ \mathcal{\cdot }\right] $ (see
Lemma \ref{combsde}), we get%
\begin{equation*}
\inf_{u\in \mathcal{U}_{ad}^{2}\left[ 0,T\right] }\left\{ \mathbb{G}%
_{t,t+\delta }^{t,\mu ;u}\left[ \varphi \left( t+\delta ,X_{t+\delta
}^{t,\mu ;u}\right) \right] -\varphi \left( t,\xi \right) \right\} \geq 0,%
\text{ a.s.}
\end{equation*}%
From the relation (\ref{rebs1}), we have $\inf_{u\in \mathcal{U}_{ad}^{2}%
\left[ 0,T\right] }\mathscr{Y}_{t}^{1,u}\geq 0,$ a.s., and furthermore, from
Lemma \ref{estvis}, we have $\inf_{u\in \mathcal{U}_{ad}^{2}\left[ 0,T\right]
}\mathscr{Y}_{t}^{2,u}\geq -C\delta \rho _{1}\left( \delta \right) .$ By
means of Lemma \ref{ode}, it follows%
\begin{equation*}
\mathscr{Y}_{t}\geq -C\delta \rho _{1}\left( \delta \right) .
\end{equation*}%
Therefore, we deduce that%
\begin{equation*}
\Upsilon \left( t,\xi ,0,0\right) =\inf_{u\in U}\mathscr{F}\left( t,\xi
,0,0,u\right) \geq 0.
\end{equation*}%
From the definition of $F$ we report that $\mathcal{V}$ is the viscosity
subsolution (resp., supersolution) of equation (\ref{hjb2}) with the
terminal condition (\ref{TT1}), i.e., $\mathcal{V}$ is the viscosity
solution.

Next we proceed the uniqueness of viscosity solution via a comparison
principle for viscosity solutions in infinite dimension to the lifted HJB
equation (\ref{hjb3}). Just like \cite{PW2017}, we need check that all the
assumptions of Theorem 3.50 in \cite{FGS2015} are fulfilled with the
framework of the lifted Hamiltonian $F$ defined in (\ref{F}).

Observe first that there is no appearance of linear dissipative operator in
the equation, thus the HJB equation (\ref{hjb3}) is a bounded equation in
the terminology of \cite{FGS2015}. Therefore, the notion of $B$-continuity
reduces to the standard notion of continuity in $L^{2}\left( \Omega ^{1},%
\mathcal{G},\mathbb{P}^{1};\mathbb{R}^{n}\right) $, namely, $B=I$ (the
identity operator)$\footnote{%
Let $A$ be a linear, densely defined, closed operator in $H.$ Whenever $A=0,$
and if one sets $B=I$, it follows immediately that $B$ satisfies the strong $%
B$-condition (Defintion 3.10 in \cite{FGS2015}). \ }$. The norm define in
\cite{FGS2015}
\begin{equation*}
\left\vert x\right\vert _{-1}^{2}=\left\vert x\right\vert ^{2}=\left\langle
x,x\right\rangle ,\text{ }x\in H,
\end{equation*}%
where $H$ is a real, separable Hilbert space with inner product $%
\left\langle \cdot ,\cdot \right\rangle .$ In particular, we set $\left\{
e_{1},e_{2},\cdots \right\} $ is an orthogonal basis of $H$, and for $N>2$, $%
H_{N}=$span$\left\{ e_{1},e_{2},\cdots e_{N}\right\} $, $P_{N}$ is the
orthogonal projection in $H$ onto $H_{N}$. We denote $S\left( H\right) $ the
set of bounded, self-adjoint operators on $H.$

We will check one by one the conditions in Theorem 3.50 in \cite{FGS2015}.
The Hypothesis 3.44 holds from the uniform continuity of $b$, $\sigma $ and $%
f$ in (A1)-(A2). As for Hypothesis 3.45, invoking (A6), Hypothesis 3.45 is
immediately satisfied. It is rather easy to check that the monotonicity
condition in $P\in L^{2}\left( \Omega ^{1},\mathcal{G},\mathbb{P}^{1};%
\mathbb{S}^{n}\right) $) of $F$ in Hypothesis 3.46 is satisfied since the
generator $f$ is independent of $P.$ Hypothesis 3.47 holds directly when
dealing with bounded equations. We now deal with Hypothesis 3.48. Let $%
X,Y\in S\left( H\right) $ such that $X=P_{N}^{\ast }XP_{N},$ $Y=P_{N}^{\ast
}YP_{N}$ for some $N$ and satisfying
\begin{equation}
-\frac{3}{\varepsilon }\left(
\begin{array}{cc}
P_{N} & 0 \\
0 & P_{N}%
\end{array}%
\right) \leq \left(
\begin{array}{cc}
X & 0 \\
0 & -Y%
\end{array}%
\right) \leq \frac{1}{\varepsilon }\left(
\begin{array}{cc}
P_{N} & -P_{N} \\
-P_{N} & P_{N}%
\end{array}%
\right) .  \label{XY}
\end{equation}%
We proceed the following term%
\begin{eqnarray*}
&&F\left( x,r,\frac{x-y}{\varepsilon },X\right) -F\left( y,r,\frac{x-y}{%
\varepsilon },Y\right) \\
&=&\inf_{u\in U}\mathbb{E}^{1}\Bigg [f\Big (x,\mathcal{L}\left( x\right) ,r,%
\mathbb{E}^{1}\left[ \frac{\left( x-y\right) ^{\top }}{\varepsilon }\sigma
\left( x,\mathcal{L}\left( \xi \right) ,u\right) \right] ,\mathcal{L}_{r},u%
\Big ) \\
&&+\left\langle \frac{\left( x-y\right) }{\varepsilon },b\left( x,\mathcal{L}%
\left( x\right) ,u\right) \right\rangle +\frac{1}{2}\left\langle X\left(
\sigma \left( x,\mathcal{L}\left( x\right) ,u\right) \right) ,\sigma \left(
x,\mathcal{L}\left( x\right) ,u\right) \right\rangle \Bigg ] \\
&&-\inf_{u\in U}\mathbb{E}^{1}\Bigg [f\Big (y,\mathcal{L}\left( y\right) ,r,%
\mathbb{E}^{1}\left[ \frac{\left( x-y\right) ^{\top }}{\varepsilon }\sigma
\left( y,\mathcal{L}\left( y\right) ,u\right) \right] ,\mathcal{L}_{r},u\Big
) \\
&&+\left\langle \frac{\left( x-y\right) ^{\top }}{\varepsilon },b\left( y,%
\mathcal{L}\left( y\right) ,u\right) \right\rangle +\frac{1}{2}\left\langle
Y\left( \sigma \left( y,\mathcal{L}\left( y\right) ,u\right) \right) ,\sigma
\left( y,\mathcal{L}\left( y\right) ,u\right) \right\rangle \Bigg ],
\end{eqnarray*}%
For any $\epsilon >0,$ there exists a measurable functions $u^{\varepsilon }:%
\mathbb{R}^{n}\rightarrow U,$ such that
\begin{eqnarray*}
&&\mathbb{E}^{1}\Bigg [f\Big (x,\mathcal{L}\left( x\right) ,r,\mathbb{E}^{1}%
\left[ \frac{\left( x-y\right) ^{\top }}{\varepsilon }\sigma \left( x,%
\mathcal{L}\left( \xi \right) ,u^{\varepsilon }\right) \right] ,\mathcal{L}%
_{r},u^{\varepsilon }\Big ) \\
&&+\left\langle \frac{\left( x-y\right) }{\varepsilon },b\left( x,\mathcal{L}%
\left( x\right) ,u^{\varepsilon }\right) \right\rangle +\frac{1}{2}%
\left\langle P\left( \sigma \left( x,\mathcal{L}\left( x\right)
,u^{\varepsilon }\right) \right) ,\sigma \left( x,\mathcal{L}\left( x\right)
,u^{\varepsilon }\right) \right\rangle \Bigg ]-\epsilon \\
&\leq &\inf_{u\in U}\mathbb{E}^{1}\Bigg [f\Big (x,\mathcal{L}\left( x\right)
,r,\mathbb{E}^{1}\left[ \frac{\left( x-y\right) ^{\top }}{\varepsilon }%
\sigma \left( x,\mathcal{L}\left( \xi \right) ,u\right) \right] ,\mathcal{L}%
_{r},u\Big ) \\
&&+\left\langle \frac{\left( x-y\right) }{\varepsilon },b\left( x,\mathcal{L}%
\left( x\right) ,u\right) \right\rangle +\frac{1}{2}\left\langle P\left(
\sigma \left( x,\mathcal{L}\left( x\right) ,u\right) \right) ,\sigma \left(
x,\mathcal{L}\left( x\right) ,u\right) \right\rangle \Bigg ].
\end{eqnarray*}%
Immediately, we have, from the condition of $b$ and $\sigma $ in (A1)-(A2),
and the uniform continuity condition on $f$ in (A3),
\begin{eqnarray*}
&&F\left( x,r,\frac{x-y}{\varepsilon },P\right) -F\left( y,r,\frac{x-y}{%
\varepsilon },\tilde{P}\right) \\
&\geq &\mathbb{E}^{1}\Bigg [f\Big (x,\mathcal{L}\left( x\right) ,r,\mathbb{E}%
^{1}\left[ \frac{\left( x-y\right) ^{\top }}{\varepsilon }\sigma \left( x,%
\mathcal{L}\left( \xi \right) ,u^{\varepsilon }\right) \right] ,\mathcal{L}%
_{r},u^{\varepsilon }\Big )
\end{eqnarray*}%
\begin{eqnarray*}
&&+\left\langle \frac{\left( x-y\right) }{\varepsilon },b\left( x,\mathcal{L}%
\left( x\right) ,u^{\varepsilon }\right) \right\rangle +\frac{1}{2}%
\left\langle P\left( \sigma \left( x,\mathcal{L}\left( x\right)
,u^{\varepsilon }\right) \right) ,\sigma \left( x,\mathcal{L}\left( x\right)
,u^{\varepsilon }\right) \right\rangle \Bigg ]-\varepsilon \\
&&-\mathbb{E}^{1}\Bigg [f\Big (y,\mathcal{L}\left( y\right) ,r,\mathbb{E}^{1}%
\left[ \frac{\left( x-y\right) ^{\top }}{\varepsilon }\sigma \left( y,%
\mathcal{L}\left( y\right) ,u^{\varepsilon }\right) \right] ,\mathcal{L}%
_{p},u^{\varepsilon }\Big ) \\
&&-\left\langle \frac{\left( x-y\right) ^{\top }}{\varepsilon },b\left( y,%
\mathcal{L}\left( y\right) ,u^{\varepsilon }\right) \right\rangle -\frac{1}{2%
}\left\langle \tilde{P}\left( \sigma \left( y,\mathcal{L}\left( y\right)
,u^{\varepsilon }\right) \right) ,\sigma \left( y,\mathcal{L}\left( y\right)
,u^{\varepsilon }\right) \right\rangle \Bigg ] \\
&\geq &-C\left( \left\vert x-y\right\vert \left( 1+\frac{\left\vert
x-y\right\vert }{\varepsilon }\right) \right) +\frac{1}{2}\left\langle
P\left( \sigma \left( x,\mathcal{L}\left( x\right) ,u^{\varepsilon }\right)
\right) ,\sigma \left( x,\mathcal{L}\left( x\right) ,u^{\varepsilon }\right)
\right\rangle \\
&&-\frac{1}{2}\left\langle \tilde{P}\left( \sigma \left( y,\mathcal{L}\left(
y\right) ,u^{\varepsilon }\right) \right) ,\sigma \left( y,\mathcal{L}\left(
y\right) ,u^{\varepsilon }\right) \right\rangle .
\end{eqnarray*}%
We proceed the term (with $\sigma ^{\varepsilon }\left( x\right) =\sigma
\left( x,\mathcal{L}\left( x\right) ,u^{\varepsilon }\right) ,$ $\sigma
^{\varepsilon }\left( y\right) =\sigma \left( y,\mathcal{L}\left( y\right)
,u^{\varepsilon }\right) $)%
\begin{eqnarray*}
&&\left\vert \text{tr}\left( \sigma ^{\varepsilon }\left( x\right) \sigma
^{\varepsilon }\left( x\right) ^{\top }X-\sigma ^{\varepsilon }\left(
y\right) \sigma ^{\varepsilon }\left( y\right) ^{\top }Y\right) \right\vert
\\
&=&\left\vert \text{tr}\left( \left(
\begin{array}{cc}
\sigma ^{\varepsilon }\left( x\right) \sigma ^{\varepsilon }\left( x\right)
^{\top } & \sigma ^{\varepsilon }\left( x\right) \sigma ^{\varepsilon
}\left( y\right) ^{\top } \\
\sigma ^{\varepsilon }\left( y\right) \sigma ^{\varepsilon }\left( x\right)
^{\top } & \sigma ^{\varepsilon }\left( y\right) \sigma ^{\varepsilon
}\left( y\right) ^{\top }%
\end{array}%
\right) \right) \right\vert \\
&\leq &\frac{3}{\varepsilon }\left\vert \text{tr}\left( \left(
\begin{array}{cc}
\sigma ^{\varepsilon }\left( x\right) \sigma ^{\varepsilon }\left( x\right)
^{\top } & \sigma ^{\varepsilon }\left( x\right) \sigma ^{\varepsilon
}\left( y\right) ^{\top } \\
\sigma ^{\varepsilon }\left( y\right) \sigma ^{\varepsilon }\left( x\right)
^{\top } & \sigma ^{\varepsilon }\left( y\right) \sigma ^{\varepsilon
}\left( y\right) ^{\top }%
\end{array}%
\right) \left(
\begin{array}{cc}
P_{N} & -P_{N} \\
-P_{N} & P_{N}%
\end{array}%
\right) \right) \right\vert \\
&=&\frac{3}{\varepsilon }\left\vert \text{tr}\left( \sigma ^{\varepsilon
}\left( x\right) \sigma ^{\varepsilon }\left( x\right) ^{\top }P_{N}-\sigma
^{\varepsilon }\left( x\right) \sigma ^{\varepsilon }\left( y\right) ^{\top
}P_{N}-\sigma ^{\varepsilon }\left( y\right) \sigma ^{\varepsilon }\left(
x\right) ^{\top }P_{N}+\sigma ^{\varepsilon }\left( y\right) \sigma
^{\varepsilon }\left( y\right) ^{\top }P_{N}\right) \right\vert \\
&\leq &C\frac{3}{\varepsilon }\left\vert x-y\right\vert ^{2}.
\end{eqnarray*}%
So
\begin{equation*}
F\left( x,r,\frac{x-y}{\varepsilon },X\right) -F\left( y,r,\frac{x-y}{%
\varepsilon },Y\right) \geq -C\left( \left\vert x-y\right\vert \left( 1+%
\frac{\left\vert x-y\right\vert }{\varepsilon }\right) \right) .
\end{equation*}

Hypothesis 3.49 follows from the growth condition of $b$ and $\sigma $ in
(H1). Consequently, from Theorem 3.50 in \cite{FGS2015}, we assert that
comparison principle holds for the HJB (\ref{hjb3}). The proof is thus
complete. \hfill $\Box $

\begin{remark}
In 1997, Barles, Buckdahn and Pardoux \cite{bbp1997} introduced the
following space of continuous functions%
\begin{equation*}
\begin{array}{ll}
\Theta & =\Bigg \{\varphi \in C\left( \left[ 0,T\right] \times \mathbb{R}%
^{n}\right) \bigg |\exists \tilde{A}>0\text{ such that} \\
& \lim_{\left\vert x\right\vert \rightarrow \infty }\left\vert \varphi
\left( t,x\right) \right\vert \exp \left\{ -\tilde{A}\left[ \log \left(
\left( \left\vert x\right\vert ^{2}+1\right) ^{\frac{1}{2}}\right) \right]
^{2}\right\} =0,\text{ uniformly in }t\in \left[ 0,T\right] \Bigg \}. \\
&
\end{array}%
\end{equation*}%
to prove the uniqueness of the viscosity solution of an integro-partial
differential equation associated with a decoupled FBSDE with jumps. As
pointed out in \cite{bbp1997}, its growth condition is slightly weaker than
the assumption of polynomial growth but more restrictive than that of
exponential growth. In addition, they asserted in \cite{bbp1997} that this
kind of growth condition is optimal for the uniqueness. Whereafter,
Buckdahn, Li and Peng \cite{blp09}, proved the uniqueness in the space of
continuous functions with at most polynomial growth with the help of the
comparison principle (see, for instance, \cite{CIL}) for PDEs whose
solutions can be stochastically interpreted in terms of mean-field BSDEs
(see Remark \ref{blp09r} for more details) with standard assumptions on the
coefficients. Note that, however, PDEs mentioned above are considered in
finite dimensional case. Besides, the structure of them are very different
from ours.
\end{remark}

As we have well known, the target of the optimal control problem is to seek
an optimal control and the corresponding state trajectory. From Theorem \ref%
{m1}, one of the potential approach is that one might be able to construct
an optimal feedback control via the value function. Due to the technical
restriction, we present the following verification theorem to gives a way of
testing whether a given admissible pair is optimal or not under smoothness
of value function.

\begin{theorem}[Verification theorem]
\label{m2}Assume that \emph{(A1)-(A6)} are in force. Suppose that $v\in
C_{b}^{1,2}\left( \left[ 0,T\right] \times \mathcal{P}_{2}\left( \mathbb{R}%
^{n}\right) \right) $ satisfying a quadratic growth condition as in (\ref%
{p22}), together with a linear growth condition for its derivative%
\begin{equation*}
|\partial _{\mu }v(t,%
\mu
)(x)|\leq C(1+|x|+\left\Vert \mu \right\Vert _{2}),\text{ }\forall (t,x,%
\mu
)\in \left[ 0,T\right] \times \mathbb{R}\times \mathcal{P}_{2}\left( \mathbb{%
R}^{n}\right) ,
\end{equation*}%
for some positive constant $C$. Let $v$ be the solution to the HJB equation (%
\ref{ham2}). Then%
\begin{equation*}
v(t,%
\mu
)\leq J\left( t,\mu ;u\right) ,\text{ }u\in \mathcal{U}_{ad}^{2}\left[ 0,T%
\right] ,\text{ }(t,%
\mu
)\in \left[ 0,T\right] \times \mathcal{P}_{2}\left( \mathbb{R}^{n}\right) .
\end{equation*}%
Furthermore, an admissible pair $\left( \bar{X}^{t,\xi ,\bar{u}},\bar{Y}%
^{t,\xi ,\bar{u}},\bar{Z}^{t,\xi ,\bar{u}},\bar{u}\right) _{0\leq s\leq T}$
is optimal for Problem MVS if and only if and there exists an element $\hat{u%
}\left( t,x,%
\mu
\right) \in U$ attaining the infimumin (\ref{hjb}) s.t. the map $(t,x,%
\mu
)\rightarrow \hat{u}\left( t,x,%
\mu
\right) $ is measurable, and
\begin{equation}
\left\{
\begin{array}{rcl}
\mathrm{d}\bar{X}_{s} & = & b\left( \bar{X}_{s},\mathbb{P}_{\bar{X}_{s}}^{W},%
\hat{u}\left( t,\bar{X}_{s},\mathbb{P}_{\bar{X}_{s}}^{W}\right) \right)
\mathrm{d}s+\sigma \left( \bar{X}_{s},\mathbb{P}_{\bar{X}_{s}}^{W},\hat{u}%
\left( s,\bar{X}_{s},\mathbb{P}_{\bar{X}_{s}}^{W}\right) \right) \mathrm{d}%
W_{s}, \\
\bar{X}_{t} & = & \xi , \\
\mathrm{d}\bar{Y}_{s} & = & -f^{0}\left( \mathbb{P}_{\bar{X}_{s}}^{W},\bar{%
\Theta}_{s},\mathbb{P}_{\bar{\Theta}_{s}},\hat{u}\left( s,\bar{X}_{s},%
\mathbb{P}_{\bar{X}_{s}}^{W}\right) \right) \mathrm{d}s+\bar{Z}_{s}\mathrm{d}%
W_{s}, \\
\bar{Y}_{T} & = & \Phi ^{0}\left( \mathbb{P}_{\bar{X}_{T}}^{W}\right) .%
\end{array}%
\right.
\end{equation}%
admits a unique solution denoted $\left( \bar{X}_{s},\bar{Y}_{s},\bar{Z}_{s},%
\hat{u}_{s}\right) _{0\leq s\leq T}$, for any $(t,\xi )\in \lbrack
0,T]\times L^{2}\left( \mathcal{G};\mathbb{R}^{n}\right) $. Then, the
admissible feed back control $\bar{u}\in \mathcal{U}_{ad}^{2}\left[ 0,T%
\right] $ defined by
\begin{equation}
\bar{u}\left( s\right) =\hat{u}\left( t,\bar{X}_{s},\mathbb{P}_{\bar{X}%
_{s}}^{W}\right) ,\text{ }t\leq s<T  \label{feed}
\end{equation}%
is an optimal control in the sense of $v(t,%
\mu
)=J\left( t,\mu ;\hat{u}\right) $ for $(t,%
\mu
)\in \left[ 0,T\right] \times \mathcal{P}_{2}\left( \mathbb{R}^{n}\right) $.
\end{theorem}

\paragraph{Proof.}

For any $(t,%
\mu
)\in \left[ 0,T\right] \times \mathcal{P}_{2}\left( \mathbb{R}^{n}\right) $
and $u\in U,$ we have%
\begin{eqnarray}
0 &=&\partial _{t}v\left( t,\mu \right) +\inf_{u\in U}\Bigg [f^{0}\Big (\mu
,v\left( t,\mathbb{\mu }\right) ,\mu \left( \partial _{\mu }v\left( t,%
\mathbb{\mu }\right) x^{\top }\sigma \left( x,\mu ,u\right) \right) ,\mathbb{%
P}_{v\left( t,\mu \right) },u\Big )  \notag \\
&&+\mu \left( \mathscr{L}^{u}v\left( t,\mu \right) \right) +\mu \otimes \mu
\left( \mathscr{M}^{u}v\left( t,\mathbb{\mu }\right) \right) \Bigg ]  \notag
\\
&\leq &\partial _{t}v\left( t,\mu \right) +\mu \left( \mathscr{L}^{u}v\left(
t,\mathbb{\mu }\right) \right) +\mu \otimes \mu \left( \mathscr{M}%
^{u}v\left( t,\mathbb{\mu }\right) \right)  \notag \\
&&+f^{0}\Big (\mu ,v\left( t,\mu \right) ,\mu \left( \partial _{\mu }v\left(
t,\mu \right) x^{\top }\sigma \left( x,\mu ,u\right) \right) ,\mathbb{P}%
_{v\left( t,\mathbb{\mu }\right) },u\Big )  \label{w1}
\end{eqnarray}

Next, for any $u\in \mathcal{U}_{ad}^{2}\left[ 0,T\right] ,$ applying It\^{o}%
's formula (\ref{ito1}) to $v(t,%
\mu
)$, we get (with $b_{s}^{u}=b_{s}\left( X_{s}^{u},u_{s}\right) ,$ $\sigma
_{s}^{u}=\sigma _{s}\left( X_{s}^{u},u_{s}\right) ,$ \textit{etc.})%
\begin{eqnarray}
v(t,%
\mu
) &=&\Phi ^{0}\left( \mathbb{P}_{X_{T}^{u}}^{W}\right) -\int_{t}^{T}\Bigg [%
\mathbb{E}^{1}\Big [\left\langle \partial _{\mu }v\left( s,\mathbb{P}%
_{X_{s}^{u}}^{W}\right) \left( X_{s}^{u}\right) ,b_{s}^{u}\right\rangle
\notag \\
&&+\frac{1}{2}\text{tr}\left( \partial _{x}\partial _{\mu }v\left( s,\mathbb{%
P}_{X_{s}^{u}}^{W}\right) \left( X_{s}^{u}\right) \sigma _{s}^{u}\left(
\sigma _{s}^{u}\right) ^{\top }\right) \Big ]  \notag \\
&&+\mathbb{E}^{1}\left[ \mathbb{\tilde{E}}^{1}\left[ \frac{1}{2}\text{tr}%
\left( \partial _{\mu }^{2}v\left( s,\mathbb{P}_{X_{s}^{u}}^{W}\right)
\left( X_{s}^{u},\tilde{X}_{s}^{u}\right) \sigma _{s}^{t,\mu ,u}\left(
\tilde{\sigma}_{s}^{u}\right) ^{\top }\right) \right] \right] \Bigg ]\mathrm{%
d}s  \notag \\
&&-\int_{t}^{T}\mathbb{E}^{1}\left[ \partial _{\mu }v\left( s,\mathbb{P}%
_{X_{s}^{u}}^{W}\right) \left( X_{s}^{u}\right) ^{\top }\sigma _{s}^{u}%
\right] \mathrm{d}W_{s}  \notag \\
&\leq &\Phi ^{0}\left( \mathbb{P}_{X_{T}^{u}}^{W}\right) +\int_{t}^{T}f^{0}%
\Big (\mathbb{P}_{X_{s}^{u}}^{W},v\left( s,\mathbb{P}_{X_{s}^{u}}^{W}\right)
,\mathbb{E}^{1}\left[ \partial _{\mu }v\left( s,\mathbb{P}%
_{X_{s}^{u}}^{W}\right) \left( X_{s}^{u}\right) ^{\top }\sigma _{s}^{u}%
\right] ,  \notag \\
&&\mathbb{P}_{v\left( s,\mathbb{P}_{X_{s}^{u}}^{W}\right) },u_{s}\Big )
\notag \\
&&-\int_{t}^{T}\mathbb{E}^{1}\left[ \partial _{\mu }v\left( s,\mathbb{P}%
_{X_{s}^{u}}^{W}\right) \left( X_{s}^{u}\right) ^{\top }\sigma _{s}^{u}%
\right] \mathrm{d}W_{s}\text{ by (\ref{w1}).}  \label{w2}
\end{eqnarray}%
For simplicity, we set
\begin{equation*}
\left\{
\begin{array}{l}
\bar{Y}_{s}^{u}=v\left( s,\mathbb{P}_{X_{s}^{u}}^{W}\right) ,\text{ }\mathbb{%
P}\text{-a.s. }t\leq s\leq T, \\
\bar{Z}_{s}^{u}=\mathbb{E}^{1}\left[ \partial _{\mu }v\left( t,\mathbb{P}%
_{X_{s}^{u}}^{W}\right) \left( X_{s}^{u}\right) ^{\top }\sigma _{s}^{u}%
\right] ,\text{ }\mathbb{P}\text{-a.s. }t\leq s\leq T.%
\end{array}%
\right.
\end{equation*}%
From (\ref{w2}), it has
\begin{equation}
\bar{Y}_{s}^{u}\leq \Phi ^{0}\left( \mathbb{P}_{X_{T}^{u}}^{W}\right)
+\int_{t}^{T}f^{0}\Big (\mathbb{P}_{X_{s}^{u}}^{W},\bar{Y}_{s}^{u},\bar{Z}%
_{s}^{u},\mathbb{P}_{\bar{Y}_{s}^{u}},u_{s}\Big )-\int_{t}^{T}\bar{Z}_{s}^{u}%
\mathrm{d}W_{s}\text{.}  \label{w3}
\end{equation}%
Now consider the following BSDE,%
\begin{equation}
Y_{t}^{u}=\Phi ^{0}\left( \mathbb{P}_{X_{T}^{u}}^{W}\right)
+\int_{t}^{T}f^{0}\left( \mathbb{P}_{X_{s}^{u}}^{W},Y_{s}^{u},Z_{s}^{u},%
\mathbb{P}_{Y_{s}^{u}}^{W},u_{s}\right) \mathrm{d}s-\int_{t}^{T}Z_{s}^{u}%
\mathrm{d}W_{s}.  \label{w4}
\end{equation}%
Put%
\begin{equation*}
\tilde{Y}_{s}^{u}=\bar{Y}_{s}^{u}-Y_{s}^{u},\text{ }t\leq s\leq T.
\end{equation*}%
Applying Tanaka-It\^{o}'s formula to formula to $\left\vert \left( \tilde{Y}%
_{s}^{u}\right) ^{+}\right\vert ^{2}$ on $\left[ 0,T\right] $ and repeating
the same argument on Lemma \ref{combsde}, we get
\begin{equation*}
v\left( t,\mu \right) =\bar{Y}_{y}^{u}\leq Y_{s}^{u}=J\left( t,\mu ;u\right)
,\text{ }u\in \mathcal{U}_{ad}^{2}\left[ 0,T\right] .
\end{equation*}%
To prove the equality for some optimal control, we employ the same It\^{o}'s
formula under the feedback control $u^{\ast }\in \mathcal{U}_{ad}^{2}\left[
0,T\right] $ constructed in (\ref{feed}), which immediately leads to the
equality in BSDE (\ref{w3}). Therefore, $v\left( t,\mu \right)
=Y_{s}^{u^{\ast }}$, and consequently, it follows the desired result:
\begin{equation*}
v\left( t,\mu \right) =J\left( t,\mu ;u^{\ast }\right) =V\left( t,\mu
\right) .
\end{equation*}%
The proof is thus complete. \hfill $\Box $

\begin{remark}
Apparently, the verification theorem (Theorem \ref{m2}) provides an analytic
feedback form of the optimal control whenever the HJB equation requires the
value function $\mathcal{V}$ to be a sufficiently smooth in the Wasserstein
space. However, as mentioned before (see Theorem \ref{m1}), the value
function can be shown a unique viscosity solution to HJB equation (\ref{hjb2}%
). Thus, it is more reasonable to investigate the same topic under the
framework of viscosity solution. We refer the similar work in Zhang \cite%
{zhang2012} (see also in \cite{zhangbook}). This issue will be our future
work under consideration.
\end{remark}

\section{LQ stochastic McKean-Vlasov control problem via $g$-expectation}

\label{sec4}

In this section, we study the classical linear-quadratic stochastic
McKean-Vlasov control problem under generalized expectation via BSDEs. In
this framework, we will see, on one hand, the admissible control set needs
more stringent integrability requirements due to the well-posedness for BSDE
regarding the terminal condition. On the other hand, from Proposition 2.1 in
Yong \cite{Y2008}, we know that the control problem under generalized
expectation is not necessarily equivalent to the classical expectation,
which means this topic is interesting in its own way.

Consider the following linear McKean-Vlasov controlled system ($d=1$):%
\begin{equation*}
\left\{
\begin{array}{rcl}
\mathrm{d}X_{s} & = & \left[ AX_{s}+\bar{A}\mathbb{E}X_{s}+Bu_{s}\right]
\mathrm{d}t+\left[ CX_{s}+\bar{C}\mathbb{E}X_{s}+Du_{s}\right] \mathrm{d}%
W_{s}, \\
X_{t} & = & \xi \in L^{2}\left( \mathcal{G};\mathbb{R}^{n}\right) ,\text{ }%
0\leq t\leq s\leq T,%
\end{array}%
\right.
\end{equation*}%
where $A,\bar{A},C,C$ are constant matrices in $\mathbb{R}^{n\times n},$ $B$
and $D$ are constant matrices in $\mathbb{R}^{n\times k},$ and $\xi \in
L^{2}\left( \mathcal{G};\mathbb{R}^{n}\right) .$ Let
\begin{eqnarray*}
\eta \left( t,\xi ,u\left( \cdot \right) \right) &=&\frac{1}{2}\Big [%
\left\langle GX_{T},X_{T}\right\rangle +\left\langle \bar{G}\mathbb{E}X_{T},%
\mathbb{E}X_{T}\right\rangle \\
&&+\int_{t}^{T}\left( X_{s}^{\top }QX_{s}+\mathbb{E}X_{s}^{\top }\bar{Q}%
\mathbb{E}X_{s}+u_{s}^{\top }Ru_{s}\right) \mathrm{d}s\Big ],
\end{eqnarray*}%
where $G,\bar{G},Q$ and $\bar{Q}$ are symmetric matrices in $\mathbb{R}%
^{n\times n}$ while $R$ is a symmetric matrix in $\mathbb{R}^{k\times k},$
and $X\left( \cdot \right) $ is the state process corresponding to $\left(
\xi ,u\left( \cdot \right) \right) \in L^{2}\left( \mathcal{G};\mathbb{R}%
^{n}\right) \times \mathcal{U}_{ad}^{2}\left[ 0,T\right] .$ Now we may
introduce the following which is referred as a \emph{cost functional:}
\begin{eqnarray*}
\mathcal{J}\left( t,\xi ;u\left( \cdot \right) \right) &=&\mathbb{E}\left[
\eta \left( t,\xi ,u\left( \cdot \right) \right) \right] \\
&=&\frac{1}{2}\mathbb{E}\Big [\left\langle GX_{T},X_{T}\right\rangle
+\left\langle \bar{G}\mathbb{E}X_{T},\mathbb{E}X_{T}\right\rangle \\
&&+\int_{t}^{T}\left( X_{s}^{\top }QX_{s}+\mathbb{E}X_{s}^{\top }\bar{Q}%
\mathbb{E}X_{s}+u_{s}^{\top }Ru_{s}\right) \mathrm{d}s\Big ].
\end{eqnarray*}

In reality, the objective expectation does not quite reflect people's
preferences (cf. \cite{All1953, E1961}). An effective way is to employ the
so-called generalized expectation, $g$-expectation, (see Peng \cite{Peng04})
which, in some sense, is subjective.

Consider
\begin{equation}
\left\{
\begin{array}{lll}
\mathrm{d}\zeta (s) & = & -g(\kappa (s))\mathrm{d}s+\kappa (s)\mathrm{d}W(s),
\\
\zeta (T) & = & \chi ,\qquad 0\leq s\leq T,%
\end{array}%
\right.  \label{bsde1}
\end{equation}%
where $g:\mathbb{R}^{d}\rightarrow \mathbb{R}$ is a given map for which we
introduce the following assumption.

\begin{description}
\item[(H1)] The map $g:\mathbb{R}^{d}\rightarrow \mathbb{R}$ satisfies the
following conditions:
\end{description}

\begin{equation}
\left\vert g(\kappa )-g(\kappa ^{\prime })\right\vert \leq C\left\vert
\kappa -\kappa ^{\prime }\right\vert ,\quad \kappa ,\kappa ^{\prime }\in
\mathbb{R}^{d}\text{,}  \label{g1}
\end{equation}%
and
\begin{equation}
g(\kappa )=0\quad \Longleftrightarrow \quad \kappa =0.  \label{g2}
\end{equation}

Assume that (H1) holds, for any $\chi \in L_{\mathcal{F}_{T}}^{2}(0,T;%
\mathbb{R})$, BSDE (\ref{g1}) admits a unique adapted solution $(\zeta
(s),\kappa (\cdot ))\in L_{\mathcal{F}}^{2}(0,T;\mathbb{R})\times L_{%
\mathcal{F}}^{2}(0,T;\mathbb{R}^{d}).$ Further, if (\ref{g2}) holds, we may
define
\begin{equation}
\mathcal{E}_{g}^{t}(\chi )=\zeta (t),\text{ }0\leq t\leq T.  \label{g3}
\end{equation}%
Clearly, the map $\chi \rightarrow \mathcal{E}_{g}^{0}(\chi )$ keeps all the
properties that $\mathbb{E}$ has, except possibly for the linearity. We call
$\mathcal{E}_{g}^{t}(\chi )$ the expectation of $\xi $ associated with $g$.
Whenever $g(\cdot )=0$, $\mathcal{E}_{g}^{t}$ is reduced to the original
(linear) conditional expectation. A typical examples of $g$ satisfying (H1)
as follows:%
\begin{equation}
g(\kappa )=\left\langle \beta ,\kappa \right\rangle ,\text{ }\beta ,\kappa
\in \mathbb{R}^{d}.  \label{eg2}
\end{equation}
Note that (\ref{eg2}) is positively homogeneous\footnote{%
Under this case, one may find the relationship between the Choquet
expectation and $g$-expectation (see \cite{CCD} for more details).}. The $g$%
-expectation, a nonlinear expectation firstly introduced by Peng \cite%
{Peng04} via a nonlinear BSDE above, can be considered as a nonlinear
extension of the well-known Girsanov transformations. The original
motivation for studying $g$-expectation comes from the theory of expected
utility, which is the foundation of modern mathematical economics. Chen and
Epstein \cite{CE} gave an application of $g$-expectation to recursive
utility. The $g$-expectation, apart from their own theoretical values, have
found important applications in various areas especially in finance. For
example, the super- and sub-pricing of contingent claims in an incomplete
market can both be captured by the $g$-probability. Ambiguity in financial
modeling can be described by the $g$-expectation. Indeed, Chen and Epstein
\cite{CE} introduced a $k$-ignorance model involving the $g$-probability to
study ambiguity aversion. The $g$-expectations have also been found to have
intimate connection with the rapidly developed risk measure theory.


With the above-defined $g$-expectation on hand, it is quite natural to
introduce the following cost functional: For any $u(\cdot )\in \mathcal{U}%
_{ad}^{p}\left[ 0,T\right] ,$ $p>2,$%
\begin{eqnarray*}
\mathcal{J}_{g}\left( t,\xi ,u\left( \cdot \right) \right) &=&\mathcal{E}%
_{g}^{t}\left[ \eta \left( t,\xi ,u\left( \cdot \right) \right) \right] \\
&&\frac{1}{2}\mathcal{E}_{g}^{t}\Big \{\left\langle
GX_{T},X_{T}\right\rangle +\left\langle \bar{G}\mathbb{E}X_{T},\mathbb{E}%
X_{T}\right\rangle \\
&&+\int_{t}^{T}\left[ X_{s}^{\top }QX_{s}+\mathbb{E}X_{s}^{\top }\bar{Q}%
\mathbb{E}X_{s}+u_{s}^{\top }Ru_{s}\right] \mathrm{d}t\Big \}.
\end{eqnarray*}%
for some suitable $g$ satisfying (H1).

\begin{remark}
We explain the reason why $p>2$ from a concrete example without mean field ($%
d=1$):
\begin{equation}
\left\{
\begin{array}{rcl}
\mathrm{d}X_{s} & = & u_{s}\mathrm{d}W_{s}, \\
-\mathrm{d}Y_{s} & = & \Big [\left\langle Q_{s}X_{s},X_{s}\right\rangle
+\left\langle R_{s}u_{s},u_{s}\right\rangle +\kappa Z_{s}\Big ]\mathrm{d}%
s-Z_{s}\mathrm{d}W(s), \\
X_{0} & = & x,\quad Y_{T}=\left\langle GX_{T},X_{T}\right\rangle ,%
\end{array}%
\right.  \label{example}
\end{equation}%
where $\kappa \in \mathbb{R}$. In order to get the explicit solution of $%
Y(\cdot ),$ we introduce the following SDE:
\begin{equation}
\left\{
\begin{array}{rcl}
\mathrm{d}\Psi _{s} & = & \kappa \Psi _{s}\mathrm{d}W_{s}, \\
\Psi _{0} & = & 1.%
\end{array}%
\right.  \label{ax}
\end{equation}%
It is easy to obtain the solution to (\ref{ax}) is
\begin{equation*}
\Psi _{s}=\exp \left\{ -\frac{1}{2}\kappa ^{2}s+\kappa W_{s}\right\} .
\end{equation*}%
Applying It\^{o}'s formula to $\Psi (s)Y(s)$, we get
\begin{equation}
Y_{0}=\mathbb{E}\left[ \Psi _{T}\left\langle GX_{T},X_{T}\right\rangle
+\int_{0}^{T}\Psi _{s}(\left\langle Q_{s}X_{s},X_{s}\right\rangle
+\left\langle R_{s}u_{s},u_{s}\right\rangle )\mathrm{d}s\right] .
\label{newf}
\end{equation}%
Define
\begin{eqnarray*}
\tilde{G} &=&\Psi _{T}G, \\
\tilde{Q}_{s} &=&\Psi _{s}Q_{s}, \\
\tilde{R}_{s} &=&\Psi _{s}R_{s},\quad s\in \left[ 0,T\right] .
\end{eqnarray*}%
Observe that (\ref{newf}) is still a quadratic functional. On one hand, the
coefficients $\tilde{Q}$ and $\tilde{R}$ are $\mathcal{F}_{s}$-adapted.
Moreover, these coefficients are unbounded. Therefore, we need to restrict $%
u(\cdot )\in \mathcal{U}_{ad}^{p}\left[ 0,T\right] $ with $p>2$ to guarantee
the finiteness of the cost functional $Y(0)$. On the other hand, to ensure $%
\left\langle \tilde{G}X(T),X(T)\right\rangle \in L_{\mathcal{F}%
_{T}}^{p}(\Omega ),$ $p>1,$ we also need $u(\cdot )\in \mathcal{U}_{ad}^{p}%
\left[ 0,T\right] $ with $p>2.$
\end{remark}

In fact, we can formulate our stochastic LQ problem with generalized
expectation as follows: For any $u(\cdot )\in \mathcal{U}_{ad}^{p}\left[ 0,T%
\right] $, let $(\zeta (\cdot ),\kappa (\cdot ))$ be the adapted solution of
(\ref{bsde1}).

Put%
\begin{equation*}
\left\{
\begin{array}{rcl}
Y_{t} & = & \zeta _{t}-\int_{0}^{t}\left[ X_{s}^{\top }QX_{s}+\mathbb{E}%
X_{s}^{\top }\bar{Q}\mathbb{E}X_{s}+u_{s}^{\top }Ru_{s}\right] \mathrm{d}s,
\\
Z_{t} & = & \kappa _{t},\text{ }t\in \left[ 0,T\right] .%
\end{array}%
\right.
\end{equation*}%
Then $(y(\cdot ),z(\cdot ))$ is the unique adapted solution to the following
BSDE:%
\begin{equation}
\left\{
\begin{array}{rcl}
-\mathrm{d}Y_{s} & = & \left[ \frac{1}{2}\left( X_{s}^{\top }QX_{s}+\mathbb{E%
}X_{s}^{\top }\bar{Q}\mathbb{E}X_{s}+u_{s}^{\top }Ru_{s}\right) +g(Z_{s})%
\right] \mathrm{d}s-Z_{s}\mathrm{d}W_{s}, \\
Y_{T} & = & \frac{1}{2}\left( \left\langle GX_{T},X_{T}\right\rangle
+\left\langle \bar{G}\mathbb{E}X_{T},\mathbb{E}X_{T}\right\rangle \right) ,%
\text{ }0\leq s\leq T.%
\end{array}%
\right.  \label{dbsde}
\end{equation}%
Therefore, the state equation takes the following form which are (decoupled)
FBSDEs (if we consider the state equation in a simple form):%
\begin{equation}
\left\{
\begin{array}{rcl}
\mathrm{d}X_{s} & = & \left[ AX_{s}+\bar{A}\mathbb{E}X_{s}+Bu_{s}\right]
\mathrm{d}s+\left[ CX_{s}+\bar{C}\mathbb{E}X_{s}+Du_{s}\right] \mathrm{d}%
W_{s}, \\
-\mathrm{d}Y_{s} & = & \left[ \frac{1}{2}\left( X_{s}^{\top }QX_{s}+\mathbb{E%
}X_{s}^{\top }\bar{Q}\mathbb{E}X_{s}+u_{s}^{\top }Ru_{s}\right) +g(Z_{s})%
\right] \mathrm{d}s-Z_{s}\mathrm{d}W_{s}, \\
X_{t} & = & \xi ,\text{ }Y_{T}=\frac{1}{2}\left( \left\langle
GX_{T},X_{T}\right\rangle +\left\langle \bar{G}\mathbb{E}X_{T},\mathbb{E}%
X_{T}\right\rangle \right) ,\text{ }0\leq t\leq s\leq T.%
\end{array}%
\right.  \label{fbsde5}
\end{equation}%
The cost functional is defined as%
\begin{eqnarray}
\mathcal{J}_{g}\left( t,\xi ,u\left( \cdot \right) \right) &=&Y_{t}  \notag
\\
&=&\frac{1}{2}\mathbb{E}\Big \{\left\langle GX_{T},X_{T}\right\rangle
+\left\langle \bar{G}\mathbb{E}X_{T},\mathbb{E}X_{T}\right\rangle  \notag \\
&&+\int_{t}^{T}\Big [2g(Z_{s},\mathbb{P}_{Z_{s}})+X_{s}^{\top }QX_{s}+%
\mathbb{E}X_{s}^{\top }\bar{Q}\mathbb{E}X_{s}+u_{s}^{\top }Ru_{s}\Big ]%
\mathrm{d}s\Big \}  \notag \\
&=&\mathcal{J}\left( t,\xi ,u\left( \cdot \right) \right) +\mathbb{E}\left[
\int_{t}^{T}g(Z_{s})\mathrm{d}s\right] .  \label{v3}
\end{eqnarray}

\begin{remark}
From Lemma \ref{vdet}, we point out that $Y_{t}$ obtained in (\ref{fbsde5})
is deterministic. Therefore, $Y_{t}$ can be attained from the last equality
in (\ref{v3}), from which it is possible for Yong (see Theorem 3.3, \cite%
{Y2008}) to investigate the existence of optimal control for classical LQ
issue under $g$-expectation.
\end{remark}

\noindent \textbf{Problem (LQ-}$\mathbf{g}$\textbf{). }Minimize (\ref{v3})
over (\ref{fbsde5}).

The goal is to seek $\bar{u}\left( \cdot \right) \in \mathcal{U}_{ad}^{p}%
\left[ 0,T\right] $ (if it ever exists), such that%
\begin{equation*}
\mathcal{J}_{g}\left( t,\xi ,\bar{u}\left( \cdot \right) \right)
=\inf_{u\left( \cdot \right) \in \mathcal{U}_{ad}^{p}\left[ 0,T\right] }%
\mathcal{J}_{g}\left( t,\xi ,u\left( \cdot \right) \right) .
\end{equation*}

\begin{remark}
From (\ref{v3}), we see that if $g$ is not identically $0$, a minimum point
of $\mathcal{J}\left( \xi ,u\left( \cdot \right) \right) $ is not
necessarily a minimum point of $\mathcal{J}_{g}\left( \xi ,u\left( \cdot
\right) \right) $, and vice versa since the appearance of the term $Z$.
\end{remark}

The similar stochastic control problem has been first studied by Yong \cite%
{Y2008} without mean filed terms. However notice that this class of
McKean-Vlasov stochastic control problems with recursive utility just lands
right on our framework.

Now we suppose that the control set $U$ is the set $U\left( \mathbb{R}^{n};%
\mathbb{R}^{k}\right) $ of \emph{Lipschitz} functions from $\mathbb{R}^{n}$
into $\mathbb{R}^{k}$. In this case, we can introduce the so called
semi-feedback control $u\left( t,x,\omega _{0}\right) $, in the sense that
it is of closed-loop form w.r.t. the state process $X_{s}$, but of open-loop
form w.r.t. the noise $W$. In other words, the random field $\mathcal{F}^{0}$%
-adapted control process $u\left( x\right) $ may be regarded equivalently as
processes valued in some functional space $U$ on $\mathbb{R}^{n}$, for
instance, a closed subset of the Polish space $C(\mathbb{R}^{n},U) $, of
continuous functions from $\mathbb{R}^{n}$ into some Euclidian space $U $.

In order to present the feedback control via Riccati equation, let us look
at the case that $g\left( \cdot \right) $ is given by (\ref{eg2}). To this
end, we proceed our problem by setting%
\begin{eqnarray*}
f^{0}\left( \mu ,y,z,\mathbb{\pi },u\right) &=&\text{Var}\left( \mu \right)
\left( Q\right) +\bar{\mu}^{\top }\left( Q+\bar{Q}\right) \bar{\mu}+%
\overline{u_{\#}\mu }^{2}\left( R\right) +\left\langle b_{1},\kappa
\right\rangle , \\
\Phi ^{0}\left( \mu \right) &=&\text{Var}\left( \mu \right) \left( G\right) +%
\bar{\mu}^{\top }\left( G+\bar{G}\right) \bar{\mu},
\end{eqnarray*}%
where for any matrix $\Pi \in \mathbb{S}^{n}=\left\{ \Lambda \in \mathbb{R}%
^{n\times n}:\Lambda ^{\top }=\Lambda \right\} $, we set
\begin{equation*}
\bar{\mu}^{2}\left( \Pi \right) =\int_{\mathbb{R}^{d}}x^{\top }\Pi x\mu
\left( \mathrm{d}x\right) ,\text{ Var}\left( \mu \right) \left( \Pi \right) =%
\bar{\mu}^{2}\left( \Pi \right) -\bar{\mu}^{\top }\left( \Pi \right) \bar{\mu%
},\text{ }\mu \in \mathcal{P}_{2}\left( \mathbb{R}^{n}\right)
\end{equation*}%
and ($u_{\#}\mu $ denotes as the probability measure on $\mathbb{R}^{k}$
induced by $u$) \
\begin{equation*}
\overline{u_{\#}\mu }=\int_{\mathbb{R}^{d}}u\left( x\right) \mu \left(
\mathrm{d}x\right) ,\text{ }\overline{u_{\#}\mu }^{2}\left( R\right) =\int_{%
\mathbb{R}^{d}}u^{\top }\left( x\right) Ru\left( x\right) \mu \left( \mathrm{%
d}x\right) ,\text{ }\mu \in \mathcal{P}_{2}\left( \mathbb{R}^{n}\right) .
\end{equation*}%
We are going to find a value function $\mathcal{V}\left( t,\mu \right) ,$ $%
\mu =\mathbb{P}_{\xi }$ of the following type:%
\begin{equation*}
\mathcal{V}\left( t,\mu \right) =\text{Var}\left( \mu \right) \left(
P_{1}\left( t\right) \right) +\bar{\mu}^{\top }\left( P_{2}\left( t\right)
\right) \bar{\mu}+\bar{\mu}^{\top }\varphi \left( t\right) +\psi \left(
t\right) ,
\end{equation*}%
where $P_{1},P_{2}\in C^{1}\left( \left[ 0,T\right] ;\mathbb{S}^{n}\right) ,$
$\varphi \in C^{1}\left( \left[ 0,T\right] ;\mathbb{R}^{n}\right) $ and $%
\psi \in C^{1}\left( \left[ 0,T\right] ;\mathbb{R}\right) .$ \newline

It is easy to compute%
\begin{eqnarray*}
\partial _{t}\mathcal{V}\left( t,\mu \right) &=&\text{Var}\left( \mu \right)
\left( P_{1}^{\prime }\left( t\right) \right) +\bar{\mu}^{\top }\left(
P_{2}^{\prime }\left( t\right) \right) \bar{\mu}+\bar{\mu}^{\top }\varphi
^{\prime }\left( t\right) +\psi ^{\prime }\left( t\right) , \\
\partial _{\mu }\mathcal{V}\left( t,\mu \right) \left( x\right)
&=&2P_{1}\left( t\right) \left( x-\bar{\mu}\right) +2P_{2}\left( t\right)
\bar{\mu}+\varphi \left( t\right) , \\
\partial _{x}\partial _{\mu }\mathcal{V}\left( t,\mu \right) \left( x\right)
&=&2P_{1}\left( t\right) , \\
\partial _{\mu }^{2}\mathcal{V}\left( t,\mu \right) \left( x,\tilde{x}%
\right) &=&2\left( P_{2}\left( t\right) -P_{1}\left( t\right) \right) .
\end{eqnarray*}%
Now we are ready to derive the $P_{1},P_{2},\varphi $ and $\psi $ according
to the HJB equation (\ref{hjb2}). First, we compare the terms in
\begin{eqnarray*}
\mathcal{V}\left( T,\mu \right) &=&\text{Var}\left( \mu \right) \left(
P_{1}\left( T\right) \right) +\bar{\mu}^{\top }\left( P_{2}\left( T\right)
\right) \bar{\mu}+\bar{\mu}^{\top }\varphi \left( T\right) +\psi \left(
T\right) \\
&=&\text{Var}\left( \mu \right) \left( G\right) +\bar{\mu}^{\top }\left( G+%
\bar{G}\right) \bar{\mu},
\end{eqnarray*}%
which implies that
\begin{equation*}
P_{1}\left( T\right) =G,\text{ }P_{2}\left( T\right) =G+\bar{G},\text{ }%
\varphi \left( T\right) =0,\text{ }\psi \left( T\right) =0.
\end{equation*}%
Next for any $\mu \in \mathcal{P}_{2}\left( \mathbb{R}^{n}\right) ,$ we have%
\begin{eqnarray}
0 &=&\inf_{u\in U\left( \mathbb{R}^{n};\mathbb{R}^{k}\right) }\Psi _{t}^{\mu
}\left( u\right) +\text{Var}\left( \mu \right) \Big (P_{1}^{\prime }\left(
t\right) +Q+C^{\top }P_{1}\left( t\right) C+P_{1}\left( t\right) \bar{A}+%
\bar{A}^{\top }P_{1}\left( t\right)  \notag \\
&&+\beta \left( P_{1}\left( t\right) C+C^{\top }P_{1}\left( t\right) \right) %
\Big )+\bar{\mu}^{\top }\Big (P_{2}^{\prime }\left( t\right) +Q+\bar{Q}%
+\left( C+\bar{C}\right) ^{\top }P_{2}\left( t\right) \left( C+\bar{C}\right)
\notag \\
&&+P_{2}\left( t\right) \left( A+\bar{A}\right) +\left( A+\bar{A}\right)
^{\top }P_{2}\left( t\right) +\beta P_{2}\left( t\right) \left( C+\bar{C}%
\right) +\beta \left( C+\bar{C}\right) ^{\top }P_{2}\left( t\right) \Big )%
\bar{\mu}  \notag \\
&&+\bar{\mu}^{\top }\Big (\varphi ^{\prime }\left( t\right) +\left[ A+\bar{A}%
+\beta \left( C+\bar{C}\right) \right] ^{\top }\varphi \left( t\right) \Big )%
+\psi ^{\prime }\left( t\right) ,  \label{o0}
\end{eqnarray}%
where
\begin{eqnarray*}
\Psi _{t}^{\mu }\left( u\right) &=&\text{Var}\left( u_{\#}\mu \right) \left(
\Pi _{1}\left( t\right) \right) +\overline{u_{\#}\mu }^{\top }\Pi _{2}\left(
t\right) \overline{u_{\#}\mu } \\
&&+2\int_{\mathbb{R}^{n}}\left( x-\bar{\mu}\right) ^{\top }\Pi _{3}\left(
t\right) u\left( x\right) \mu \left( \mathrm{d}x\right) \\
&&+2\bar{\mu}^{\top }\Pi _{4}\left( t\right) +\Pi _{5}\left( t\right)
\overline{u_{\#}\mu },
\end{eqnarray*}%
with
\begin{eqnarray*}
\Pi _{1}\left( t\right) &=&D^{\top }P_{1}\left( t\right) D+R, \\
\Pi _{2}\left( t\right) &=&D^{\top }P_{2}\left( t\right) D+R, \\
\Pi _{3}\left( t\right) &=&C^{\top }P_{1}\left( t\right) D+P_{1}\left(
t\right) \left( B+\beta D\right) , \\
\Pi _{4}\left( t\right) &=&\left( C+\bar{C}\right) ^{\top }P_{2}\left(
t\right) D+P_{2}\left( t\right) \left( B+\beta D\right), \\
\Pi _{5}\left( t\right) &=&\left( B+\beta D\right) ^{\top }\varphi \left(
t\right) .
\end{eqnarray*}%
To represent the optimal control, we add the following assumption:

\begin{enumerate}
\item[\textbf{(A8)}] Suppose that
\begin{equation}
Q+\bar{Q}\geq 0,Q\geq 0,G+\bar{G}\geq 0,G\geq 0,R\geq \delta I_{k}
\label{A6}
\end{equation}%
for some positive constant $\delta .$
\end{enumerate}

Under (A8), we know that $\Pi _{1}\left( t\right) $ and $\Pi _{2}\left(
t\right) $ are invertible. After some tedious calculation (basing on the
square completion), we are able to show the minimum of $\Psi _{t}^{\mu
}\left( u\right) $ over $U\left( \mathbb{R}^{n};\mathbb{R}^{k}\right) $ is
attained at
\begin{equation}
u^{\ast }\left( t,x,\mu \right) =-\frac{1}{2}\Pi _{1}^{-1}\left( t\right)
\Pi _{3}^{\top }\left( t\right) \left( x-\bar{\mu}\right) -\Pi
_{2}^{-1}\left( t\right) \Pi _{4}\left( t\right) \bar{\mu}-\frac{1}{2}\Pi
_{2}^{-1}\left( t\right) \Pi _{5}\left( t\right) .  \label{o1}
\end{equation}%
The corresponding
\begin{eqnarray}
\Psi _{t}^{\mu }\left( u\right) &=&\text{Var}\left( u-u^{\ast }\left(
t,\cdot ,\mu \right) \right) \left( \Pi _{1}\left( t\right) \right)  \notag
\\
&&+\overline{u-u^{\ast }\left( t,\cdot ,\mu \right) _{\#}\mu }^{\top }\Pi
_{2}\left( t\right) \overline{u-u^{\ast }\left( t,\cdot ,\mu \right)
_{\#}\mu }  \notag \\
&&-\text{Var}\left( \mu \right) \Pi _{3}\left( t\right) \Pi _{1}^{-1}\left(
t\right) \Pi _{3}^{\top }\left( t\right) -\bar{\mu}^{\top }\Pi _{4}\left(
t\right) \Pi _{2}^{-1}\left( t\right) \Pi _{4}^{\top }\left( t\right) \bar{%
\mu}  \notag \\
&&-\Pi _{5}^{\top }\left( t\right) \Pi _{2}^{-1}\left( t\right) \Pi
_{4}^{\top }\left( t\right) \bar{\mu}-\frac{1}{4}\Pi _{5}^{\top }\left(
t\right) \Pi _{2}^{-1}\left( t\right) \Pi _{5}\left( t\right) .  \label{o2}
\end{eqnarray}%
Now substituting (\ref{o1}) into (\ref{o2}) and combing to (\ref{o0}), we
derive the following deterministic differential equations (comparing the
coefficients of Var$\left( \cdot \right) ,$ $\bar{\mu}^{\top }\left( \cdot
\right) \bar{\mu}$ and $\bar{\mu}$)%
\begin{equation}
\left\{
\begin{array}{rcl}
0 & = & P_{1}^{\prime }\left( t\right) +Q+C^{\top }P_{1}\left( t\right)
C+P_{1}\left( t\right) \bar{A}+\bar{A}^{\top }P_{1}\left( t\right) \\
&  & +\beta \left( P_{1}\left( t\right) C+C^{\top }P_{1}\left( t\right)
\right) -\Pi _{3}\left( t\right) \Pi _{1}^{-1}\left( t\right) \Pi _{3}^{\top
}\left( t\right) , \\
P_{1}\left( t\right) & = & G,%
\end{array}%
\right.  \label{Ri1}
\end{equation}%
\begin{equation}
\left\{
\begin{array}{rcl}
0 & = & P_{2}^{\prime }\left( t\right) +Q+\bar{Q}+\left( C+\bar{C}\right)
^{\top }P_{2}\left( t\right) \left( C+\bar{C}\right) \\
&  & +P_{2}\left( t\right) \left( A+\bar{A}\right) +\left( A+\bar{A}\right)
^{\top }P_{2}\left( t\right) +\beta P_{2}\left( t\right) \left( C+\bar{C}%
\right) \\
&  & +\beta \left( C+\bar{C}\right) ^{\top }P_{2}\left( t\right) -\Pi
_{4}\left( t\right) \Pi _{2}^{-1}\left( t\right) \Pi _{4}^{\top }\left(
t\right) \\
P_{2}\left( t\right) & = & G+\bar{G},%
\end{array}%
\right.  \label{Ri2}
\end{equation}%
\begin{equation}
\left\{
\begin{array}{rcl}
0 & = & \varphi ^{\prime }\left( t\right) +\left[ A+\bar{A}+\beta \left( C+%
\bar{C}\right) \right] ^{\top }\varphi \left( t\right) \\
&  & -\Pi _{5}^{\top }\left( t\right) \Pi _{2}^{-1}\left( t\right) \Pi
_{4}^{\top }\left( t\right) , \\
\varphi \left( T\right) & = & 0%
\end{array}%
\right.  \label{Ri3}
\end{equation}%
and%
\begin{equation}
\left\{
\begin{array}{rcl}
0 & = & \psi ^{\prime }\left( t\right) -\frac{1}{4}\Pi _{5}^{\top }\left(
t\right) \Pi _{2}^{-1}\left( t\right) \Pi _{5}\left( t\right) \\
\varphi \left( T\right) & = & 0,%
\end{array}%
\right.  \label{Ri4}
\end{equation}%
Under the previous assumptions (especially (A8)), the Riccati equations (\ref%
{Ri1})-(\ref{Ri4}) have a unique solution (see \cite{won}). Note that all
the coefficients are deterministic and bounded. Consequently, from the
expression of $u^{\ast }(t,x,\mu )$ in (\ref{o1}), we derive that it is a
candidate of optimal control for Problem (LQ-$g$). At this moment, the
closed-loop system reads%
\begin{equation}
\left\{
\begin{array}{rcl}
\mathrm{d}X_{s} & = & \left[ AX_{s}+\bar{A}\mathbb{E}X_{s}+Bu^{\ast }\left(
s,X_{s},\mathbb{P}_{X_{s}}\right) \right] \mathrm{d}t \\
&  & +\left[ CX_{s}+\bar{C}\mathbb{E}X_{s}+Du^{\ast }\left( s,X_{s},\mathbb{P%
}_{X_{s}}\right) \right] \mathrm{d}W_{s}, \\
X_{t} & = & \xi \in L^{2}\left( \mathcal{G};\mathbb{R}^{n}\right) ,\text{ }%
0\leq t\leq s\leq T,%
\end{array}%
\right.  \label{fsde}
\end{equation}%
Therefore, the state equation (\ref{fsde}) $X\in \mathcal{H}^{p}\left[ 0,T%
\right] $, where%
\begin{equation*}
\begin{array}{lll}
\mathcal{H}^{p}\left[ 0,T\right] & \triangleq & \Big \{X:\left[ 0,T\right]
\times \Omega \rightarrow \mathbb{R}^{n}\text{ is }\mathcal{F}\text{-adapted
and continuous process} \\
&  & \text{ satisfying }\mathbb{E}\left[ \sup_{0\leq t\leq T}\left\vert
X_{t}\right\vert ^{p}\right] <\infty \Big \}.%
\end{array}%
\end{equation*}%
for any $p\geq 1$ and thus $u^{\ast }\in \mathcal{U}_{ad}^{p}\left[ 0,T%
\right] $, for any $p\geq 1$. From this, we assert that BSDE in (\ref{fbsde5}%
) admits a unique adapted solution $(Y_{s},Z_{s}(\cdot ))_{t\leq s\leq T}$.
Consequently, $Y_{t}$ as a solution is the optimal value of the cost (\ref%
{v3}).

\begin{remark}
If $\beta \equiv 0,$ as mentioned before, it degenerates the classical
McKean-Vlasov LQ problem (see \cite{PW2017}). However, the
existence/nonexistence of optimal control for Problem (LQ-$g$) with $p>2$,
under (H1) only are still elusive issue.
\end{remark}

\section{Concluding remark}

\label{sec5}

We analyze the optimal control problem, where the state processes are
governed by FBSDEs of McKean-Vlasov type. The value function $\mathcal{V}$
defined by the solution to the backward component of controlled FBSDEs can
be proved a nonlinear extended Feynman-Kac representation in terms of a
class of FBSDEs. We explore the probabilistic representation rigorously to
prove the DPP for value function via BSDE method. As we have well known, the
DPP is playing a crucial part in deriving a characterization of the value
function as a solution of a nonlinear partial differential equation (HJB
equation) on the Wasserstein space of measures basing on the law-invariant
of $\mathcal{V}.$ It is necessary to point out that the usual way of solving
these equations is employing the Pontryagin maximum principle, which
requires some smoothness of coefficients and convexity assumptions due to
the requirement of variation equality and adjoint equations. We should point
it out that the classical object functional for forward optimal control of
McKean-Vlasov type is spelled out by the expectation to the running function
$f$ and terminal cost $\Phi $ w.r.t. state and its distribution, which
undoubtedly is deterministic no matter of the selection of initial
condition. Besides, notice that in such a case, from view point of HJB
equation, $f$ can be only depend on state variable, badly off for the value
function $\mathcal{V}$ and its gradient. In this paper, we prove the
functional is deterministic by virtue of BSDE's property and the viscosity
property together with a uniqueness result for the value function. The LQ
problem of McKean-Vlasov type under nonlinear expectation ($g$-expectation)
has been attacked as a theoretic result. Nevertheless, there are still some
key points should be stated: 1) We obtain the uniqueness of viscosity
solution via comparison principles in the separable Hilbert space. It is
interesting to consider this problem in the Wasserstein space, or more
generally in metric spaces (see \cite{cgkp24}). 2) The Wiener process $W$
accounts for the common random environment in literature of mean field games
in which all the individuals evolve, often called common noise (cf. \cite%
{cdl16}). Therefore, it is quite natural to bring the idiosyncratic noises,
independent of $W$, into the systems. Without doubt, it makes the analysis
rather complicated.

\section{Appendix: Technique Lemmas\label{Techlemmas}}

\subsection{Proof of Lemma \protect\ref{estbsde}}

\paragraph{Proof.}

The first inequality (\ref{best1}) can be obtained directly by (\ref%
{estbdsde2}). We are going to prove (\ref{best2}) and (\ref{best3}). Let%
\begin{eqnarray*}
f^{1}\left( y^{1},z^{1}\right) &=&f\left( x^{1},y^{1},z^{1},\mu
^{1},u^{1}\right) , \\
f^{2}\left( y^{2},z^{2}\right) &=&f\left( x^{2},y^{2},z^{2},\mu
^{2},u^{2}\right) ,
\end{eqnarray*}%
where $\left( x^{i},y^{i},z^{i},\mu ^{i},u^{i}\right) \in \mathbb{R}%
^{n}\times \mathbb{R}\times \mathbb{R}^{d}\times \mathcal{P}_{2}\left(
\mathbb{R}^{n}\times \mathbb{R}\times \mathbb{R}^{d}\right) \times U,$ $%
i=1,2.$ Immediately,
\begin{equation*}
\left\vert f^{1}\left( y^{1},z^{1}\right) -f^{2}\left( y^{1},z^{1}\right)
\right\vert \leq C\big [\left\vert x^{1}-x^{2}\right\vert +\mathcal{W}%
_{2}\left( \mu ^{1},\mu ^{2}\right) +\left\vert u^{1}-u^{2}\right\vert \big ]
\end{equation*}%
For $\xi ,$ $\xi ^{\prime }\in L^{2}\left( \mathcal{G};\mathbb{R}^{n}\right)
,$ set
\begin{eqnarray*}
\hat{X}_{s} &=&X_{s}^{t,\xi ;u}-X_{s}^{t,\xi ^{\prime };u^{\prime }},\hat{Y}%
_{s}=Y_{s}^{t,\xi ;u}-Y_{s}^{t,\xi ^{\prime };u^{\prime }}, \\
\hat{Z}_{s} &=&Z_{s}^{t,\xi ;u}-Z_{s}^{t,\xi ^{\prime };u^{\prime }},\hat{u}%
_{s}=u_{s}-u_{s}^{\prime },\text{ }s\in \left[ t,T\right] .
\end{eqnarray*}%
From Lemma \ref{l1}, we get, for some $\beta \geq 2,$ there exists a
positive constant $C_{\beta }$ (depending on $\beta ,$ $T$ and Lipschitz
constant) such that
\begin{eqnarray}
&&\mathbb{E}\left[ \sup_{0\leq s\leq T}\left\vert \hat{Y}_{s}\right\vert
^{\beta }+\left( \int_{0}^{T}\left\vert \hat{Z}_{s}\right\vert ^{2}\mathrm{d}%
s\right) ^{\frac{\beta }{2}}\right]  \notag \\
&\leq &C_{\beta }\mathbb{E}\Bigg [\left\vert \Phi \left( X_{T}^{t,\xi ;u},%
\mathbb{P}_{X_{T}^{t,\xi ;u}}\right) -\Phi \left( X_{T}^{t,\xi ^{\prime
};u^{\prime }},\mathbb{P}_{X_{T}^{t,\xi ^{\prime };u^{\prime }}}\right)
\right\vert ^{\beta }  \notag \\
&&+\left( \int_{t}^{T}\left\vert f^{1}\left( Y_{s}^{t,\xi ;u},Z_{s}^{t,\xi
;u}\right) -f^{2}\left( Y_{s}^{t,\xi ;u},Z_{s}^{t,\xi ;u}\right) \right\vert
\mathrm{d}s\right) ^{\beta }\Bigg ]  \notag \\
&\leq &C\mathbb{E}\left[ \left\vert \hat{X}_{T}\right\vert ^{\beta }\right]
+C\mathcal{W}_{2}^{\beta }\left( \mathbb{P}_{X_{T}^{t,\xi ;u}},\mathbb{P}%
_{X_{T}^{t,\xi ^{\prime };u^{\prime }}}\right)  \notag \\
&&+C\mathbb{E}\Bigg \{\int_{t}^{T}\Big [\left\vert \hat{X}_{s}\right\vert
^{2}+\left\vert \hat{u}_{s}\right\vert ^{2}+\mathcal{W}_{2}^{2}\left(
\mathbb{P}_{\Theta _{s}^{t,\xi ;u}},\mathbb{P}_{\Theta _{s}^{t,\xi ^{\prime
};u^{\prime }}}\right)  \notag \\
&&+\mathcal{W}_{2}^{2}\left( \mathbb{P}_{u_{s}^{1}},\mathbb{P}%
_{u_{s}^{2}}\right) \Big ]^{\frac{\beta }{2}}\mathrm{d}s\Bigg \}.
\label{aest1}
\end{eqnarray}%
Next, we will estimate the term $\int_{t}^{T}\mathcal{W}_{2}^{2}\left(
\mathbb{P}_{\Theta _{s}^{t,\xi ;u}},\mathbb{P}_{\Theta _{s}^{t,\xi ^{\prime
};u^{\prime }}}\right) \mathrm{d}s.$ To this end, we employ a classical
estimates for the solutions of BSDE (see Lemma A.1 in Li \cite{Li2018}). For
all $\eta >0$, there exists a suitable constant $\gamma $ (depending on
Lipschitz constant and $\eta $) such that
\begin{eqnarray}
&&\left\vert \hat{Y}_{t}\right\vert ^{2}+\frac{1}{2}\mathbb{E}\left[
\int_{t}^{T}e^{\gamma \left( s-t\right) }\left( \left\vert \hat{Y}%
_{s}\right\vert ^{2}+\left\vert \hat{Z}_{s}\right\vert ^{2}\right) \mathrm{d}%
s\right]  \notag \\
&\leq &\mathbb{E}\left[ e^{\gamma \left( T-t\right) }\left\vert \Phi \left(
X_{T}^{t,\xi ;u},\mathbb{P}_{X_{T}^{t,\xi ;u}}\right) -\Phi \left(
X_{T}^{t,\xi ^{\prime };u^{\prime }},\mathbb{P}_{X_{T}^{t,\xi ^{\prime
};u^{\prime }}}\right) \right\vert ^{2}\right]  \notag \\
&&+C\eta \mathbb{E}\left[ \int_{t}^{T}e^{\gamma \left( s-t\right)
}\left\vert f^{1}\left( Y_{s}^{t,\xi ;u},Z_{s}^{t,\xi ;u}\right)
-f^{2}\left( Y_{s}^{t,\xi ;u},Z_{s}^{t,\xi ;u}\right) \right\vert \mathrm{d}s%
\right]  \notag \\
&\leq &C\mathbb{E}\left[ \left\vert \hat{X}_{T}\right\vert ^{2}\right] +C%
\mathcal{W}_{2}^{2}\left( \mathbb{P}_{X_{T}^{t,\xi ;u}},\mathbb{P}%
_{X_{T}^{t,\xi ^{\prime };u^{\prime }}}\right)  \notag \\
&&+C\eta \mathbb{E}\Bigg [\int_{t}^{T}e^{\gamma \left( s-t\right) }\Big (%
\left\vert \hat{X}_{s}\right\vert ^{2}+\left\vert \hat{u}_{s}\right\vert
^{2}+\mathcal{W}_{2}^{2}\left( \mathbb{P}_{\Theta _{s}^{t,\xi ;u}},\mathbb{P}%
_{\Theta _{s}^{t,\xi ^{\prime };u^{\prime }}}\right)  \notag \\
&&+\mathcal{W}_{2}^{2}\left( \mathbb{P}_{u_{s}},\mathbb{P}_{u_{s}^{\prime
}}\right) \Big )\mathrm{d}s\Bigg ],  \label{aest2}
\end{eqnarray}%
where $C$ depends only on the Lipschitz constants of $f,\Phi .$ Notice that,
from Lemma \ref{l0}
\begin{equation}
\mathbb{E}\left[ \left\vert \hat{X}_{T}\right\vert ^{2}\right] +\mathcal{W}%
_{2}^{2}\left( \mathbb{P}_{X_{T}^{t,\xi ;u}},\mathbb{P}_{X_{T}^{t,\xi
^{\prime };u^{\prime }}}\right) \leq C_{2}\left[ \mathcal{W}_{2}^{2}\left(
\mathbb{P}_{\xi },\mathbb{P}_{\xi ^{\prime }}\right) +\Upsilon
_{0,T}^{u,u^{\prime };2}\right] ,  \label{aest3}
\end{equation}%
where $C_{2}$ depends only on the Lipschitz constants of $b,\sigma .$
Therefore, from (\ref{aest2}) and the definition of $2$-Wasserstein metric,
we attain%
\begin{eqnarray}
&&\frac{1}{2}\mathbb{E}\left[ \int_{t}^{T}e^{\gamma \left( s-t\right)
}\left( \left\vert \hat{Y}_{s}\right\vert ^{2}+\left\vert \hat{Z}%
_{s}\right\vert ^{2}\right) \mathrm{d}s\right]  \notag \\
&\leq &C_{2,\gamma ,\eta }\left[ \mathcal{W}_{2}^{2}\left( \mathbb{P}_{\xi },%
\mathbb{P}_{\xi ^{\prime }}\right) +\Upsilon _{0,T}^{u,u^{\prime };2}\right]
\notag \\
&&+\tilde{C}\eta \mathbb{E}\Bigg [\int_{t}^{T}e^{\gamma \left( s-t\right) }%
\Big (\left\vert \hat{u}_{s}\right\vert ^{2}+\mathcal{W}_{2}^{2}\left(
\mathbb{P}_{u_{s}},\mathbb{P}_{u_{s}^{\prime }}\right) \Big )\mathrm{d}s
\notag \\
&&+\tilde{C}\eta \mathbb{E}\Bigg [\int_{t}^{T}e^{\gamma \left( s-t\right)
}\left( \left\vert \hat{Y}_{s}\right\vert ^{2}+\left\vert \hat{Z}%
_{s}\right\vert ^{2}\right) \Bigg ],  \label{aest4}
\end{eqnarray}%
where $C_{2,\gamma ,\eta }$ depends only on the Lipschitz constants of $%
f,\Phi $, $\gamma $ and $\eta ,$ while $\tilde{C}$ depends on Lipschitz
constants of $f,\Phi $. Choosing $\eta >0$ small enough such that $\tilde{C}%
\eta <1/2,$ we get%
\begin{eqnarray*}
&&\mathbb{E}\left[ \int_{t}^{T}e^{\gamma \left( s-t\right) }\left(
\left\vert \hat{Y}_{s}\right\vert ^{2}+\left\vert \hat{Z}_{s}\right\vert
^{2}\right) \mathrm{d}s\right] \\
&\leq &C\left[ \mathcal{W}_{2}^{2}\left( \mathbb{P}_{\xi },\mathbb{P}_{\xi
^{\prime }}\right) +\Upsilon _{0,T}^{u,u^{\prime };2}\right] .
\end{eqnarray*}%
Furthermore, from Lemma \ref{l0} and the properties of $\mathcal{W}_{2}$, we
get%
\begin{eqnarray*}
&&\int_{t}^{T}\mathcal{W}_{2}^{2}\left( \mathbb{P}_{\Theta _{s}^{t,\xi ;u}},%
\mathbb{P}_{\Theta _{s}^{t,\xi ^{\prime };u^{\prime }}}\right) \mathrm{d}s \\
&\leq &\mathbb{E}\left[ \int_{t}^{T}\left( \left\vert \hat{X}_{s}\right\vert
^{2}+\left\vert \hat{Y}_{s}\right\vert ^{2}+\left\vert \hat{Z}%
_{s}\right\vert ^{2}\right) \mathrm{d}s\right] \\
&\leq &C\left[ \mathcal{W}_{2}^{2}\left( \mathbb{P}_{\xi },\mathbb{P}_{\xi
^{\prime }}\right) +\Upsilon _{0,T}^{u,u^{\prime };2}\right] .
\end{eqnarray*}%
The proof is complete. \hfill $\Box $

\subsection{Proof of Lemma \protect\ref{combsde}}

\paragraph{Proof.}

Under (A3)-(A4), we know that BSDE (\ref{cbsde}), for $i=1,2,$ has a unique
solution $\left( Y^{i},Z^{i}\right) $, respectively. Set
\begin{equation*}
\hat{\xi}=\xi ^{2}-\xi ^{1},\text{ }\hat{Y}_{t}=Y_{t}^{2}-Y_{t}^{1},\text{ }%
\hat{Z}_{t}=Z_{t}^{2}-Z_{t}^{1}.
\end{equation*}%
Applying Tanaka-It\^{o}'s formula to formula to $\left\vert \hat{Y}%
_{t}^{+}\right\vert ^{2}$ on $\left[ 0,T\right] ,$ we have%
\begin{eqnarray*}
&&\left\vert \hat{Y}_{t}^{+}\right\vert ^{2}+\mathbb{E}^{\mathcal{F}_{t}}%
\left[ \int_{t}^{T}\sum_{i=1}^{d}\left\vert \hat{Z}_{s}^{i}\right\vert ^{2}%
\mathbf{I}_{\left\{ Y_{s}^{1}<Y_{s}^{2}\right\} }\mathrm{d}s\right] \\
&=&2\mathbb{E}^{\mathcal{F}_{t}}\left[ \int_{t}^{T}\hat{Y}_{s}^{+}\left[
f^{2}\left( Y_{s}^{2},Z_{s}^{2},\mathbb{P}_{Y_{s}^{2}}\right) -f^{1}\left(
Y_{s}^{1},Z_{s}^{1},\mathbb{P}_{Y_{s}^{1}}\right) \right] \mathbf{I}%
_{\left\{ Y_{s}^{1}<Y_{s}^{2}\right\} }\mathrm{d}s\right] \\
&&-\mathbb{E}^{\mathcal{F}_{t}}\left[ \int_{t}^{T}\hat{Y}_{s}^{+}\mathrm{d}%
L_{s}\right] ,
\end{eqnarray*}%
where $L_{s}$ is local time for $Y_{s}^{2}-Y_{s}^{1}$ at $0$ which claims $%
\mathbb{E}^{\mathcal{F}_{t}}\int_{t}^{T}\hat{Y}_{s}^{+}\mathrm{d}L_{s}=0.$
We deal with the first term on the right side of the equality. \newline
Put%
\begin{eqnarray*}
&&f^{2}\left( Y_{s}^{2},Z_{s}^{2},\mathbb{P}_{Y_{s}^{2}}\right) -f^{1}\left(
Y_{s}^{1},Z_{s}^{1},\mathbb{P}_{Y_{s}^{1}}\right) \\
&=&\underset{I_{1}}{\underbrace{f^{2}\left( Y_{s}^{2},Z_{s}^{2},\mathbb{P}%
_{Y_{s}^{2}}\right) -f^{2}\left( Y_{s}^{1},Z_{s}^{1},\mathbb{P}%
_{Y_{s}^{1}}\right) }}+\underset{I_{2}}{\underbrace{f^{2}\left(
Y_{s}^{1},Z_{s}^{1},\mathbb{P}_{Y_{s}^{1}}\right) -f^{1}\left(
Y_{s}^{1},Z_{s}^{1},\mathbb{P}_{Y_{s}^{1}}\right) }}.
\end{eqnarray*}%
Immediately, according to (ii), $I_{2}\leq 0$. From the inequality $2ab\leq
a^{2}+b^{2}$, the Lipschitz condition of $f^{2}$ and (iii), we have the%
\begin{equation*}
\left\vert \left[ f^{2}\left( Y_{s}^{2},Z_{s}^{2},\mathbb{P}%
_{Y_{s}^{2},}\right) -f^{2}\left( Y_{s}^{1},Z_{s}^{1},\mathbb{P}%
_{Y_{s}^{1}}\right) \right] \mathbf{I}_{\left\{ Y_{s}^{1}<Y_{s}^{2}\right\}
}\right\vert \leq C\left( \left\vert \hat{Y}_{s}^{+}\right\vert +\left\vert
\hat{Z}_{s}^{+}\right\vert \right) \mathbf{I}_{\left\{
Y_{s}^{1}<Y_{s}^{2}\right\} }.
\end{equation*}%
Combing the fact that $\hat{Y}_{t}^{+}$ is deterministic, it follows that ($%
2ab\leq \delta a^{2}+\frac{1}{\delta }b^{2}$, $\forall \delta >0$)
\begin{eqnarray*}
&&\left\vert \hat{Y}_{t}^{+}\right\vert ^{2}+\mathbb{E}\left[
\int_{t}^{T}\sum_{i=1}^{d}\left\vert \hat{Z}_{s}^{i}\right\vert ^{2}\mathbf{I%
}_{\left\{ Y_{s}^{1}<Y_{s}^{2}\right\} }\mathrm{d}s\right] \\
&\leq &\mathbb{E}\Bigg \{\int_{t}^{T}\delta \left\vert \hat{Y}%
_{s}^{+}\right\vert ^{2}+\frac{1}{\delta }C\left[ \left\vert \hat{Y}%
_{s}^{+}\right\vert ^{2}+\left\vert \hat{Z}_{s}^{+}\right\vert ^{2}\right]
\mathbf{I}_{\left\{ Y_{s}^{1}<Y_{s}^{2}\right\} }\mathrm{d}s\Bigg \} \\
&\leq &C\mathbb{E}\left\{ \int_{t}^{T}\left[ \left( \delta +\frac{1}{\delta }%
C\right) \left\vert \hat{Y}_{s}^{+}\right\vert ^{2}+\frac{1}{\delta }%
C\sum_{i=1}^{d}\left\vert \hat{Z}_{s}^{i}\right\vert ^{2}\right] \mathbf{I}%
_{\left\{ Y_{s}^{1}<Y_{s}^{2}\right\} }\mathrm{d}s\right\} .
\end{eqnarray*}%
We can take $\delta $ small enough such that $1-\frac{1}{\delta }C>0.$ Then%
\begin{eqnarray*}
&&\left\vert \hat{Y}_{t}^{+}\right\vert ^{2}+\left( 1-\frac{1}{\delta }%
C\right) \mathbb{E}\left[ \int_{t}^{T}\sum_{i=1}^{d}\left\vert \hat{Z}%
_{s}^{i}\right\vert ^{2}\mathbf{I}_{\left\{ Y_{s}^{1}<Y_{s}^{2}\right\} }%
\mathrm{d}s\right] \\
&\leq &\mathbb{E}\left[ \int_{t}^{T}\left( \delta +\frac{1}{\delta }C\right)
\left\vert \hat{Y}_{s}^{+}\right\vert ^{2}\mathbf{I}_{\left\{
Y_{s}^{1}<Y_{s}^{2}\right\} }\mathrm{d}s\right] .
\end{eqnarray*}%
The proof is complete by virtue of Gronwall inequality. \hfill $\Box $

\subsection{Proof of Lemma \protect\ref{rdet}}

For $\xi \in L^{2}\left( \mathcal{F}_{t}^{0};\mathbb{R}^{n}\right) ,$ we
first consider its simple form. In fact, if $\xi =x,$ $x\in \mathbb{R}^{n},$
given $u\in \Lambda _{t,x},$ we can prove that $Y_{t}^{t,x;u}$ is a
deterministic constant. Indeed, for every $\left( t,x\right) \in \left[ 0,T%
\right] \times \mathbb{R}^{n},$ the solution $X_{s}^{t,x;u}$ $\left( \mathbb{%
P}_{X_{s}^{t,x;u}}^{W}=\mathbb{P}_{X_{s}^{t,x;u}}^{0}\right) $ of SDE in (%
\ref{FBSDE}) is $\mathcal{F}_{s}^{t}$ -adapted which also means $%
X_{T}^{t,x;u}$ $\left( \mathbb{P}_{X_{T}^{t,x;u}}^{0}\right) $ is $\mathcal{F%
}_{T}^{t}$-measurable and consequently, the term $f\left( \Theta
_{s}^{t,x,u},\mathbb{P}_{\Theta _{s}^{t,x;u}},u_{s},\mathbb{P}%
_{u_{s}}\right) $ in (\ref{FBSDE}) is also $\mathcal{F}_{s}^{t}$-adapted
according to property of BSDE. Then from the following Lemma \ref{deter},
the solution $\left\{ \Theta _{s}^{t,x,u}\right\} _{t\leq s\leq T}$ is also $%
\mathcal{F}_{s}^{t}$-adapted, from which derives that $Y_{t}^{t,x;u}$ is
deterministic. Afterward, as for the general case, namely, $\xi \in
L^{2}\left( \mathcal{F}_{t}^{0};\mathbb{R}^{n}\right) $ will be constructed
by the approximation of simple random variables in $L^{2}\left( \mathcal{F}%
_{t}^{0};\mathbb{R}^{n}\right) .$

\begin{definition}
For any $t\in \left[ 0,T\right] $, a sequence $\left\{ A_{k}\right\}
_{k=1,\ldots ,N}$ with $1\leq N\leq \infty $ is called a partition of $%
\left( \Omega ,\mathcal{F}_{t}\right) $ if $\bigcup_{k=1}^{N}A_{k}=\Omega $
and $A_{k}\cap A_{l}=\emptyset $, whenever $k\neq l$.
\end{definition}

\begin{remark}
When $\xi \in L^{2}\left( \mathcal{F}_{t}^{0};\mathbb{R}^{n}\right) ,$ $\xi $
is independent of $\Omega ^{1}.$ Then, $\mathbb{P}_{X_{s}^{t,\xi ,u}}^{W}=%
\mathbb{P}_{X_{s}^{t,\xi ,u}\left( \cdot \right) }^{0},$ $s\in \left[ t,T%
\right] .$
\end{remark}

\paragraph{Proof.}

Consider a simple random variable $\xi $ of the form: $\xi
=\sum_{i=1}^{N}I_{A_{i}}x^{i}$, and a control process $u$ of the form $%
u_{\left( \cdot \right) }=\sum_{i=1}^{N}I_{A_{i}}u_{\cdot }^{i}$, where $%
\left\{ A_{i}\right\} _{i=1,\ldots ,N}$ is a finite partition of $\left(
\Omega ,\mathcal{F}^{0}\right) ,$ $x^{i}\in \mathbb{R}^{n}$ and $u_{\cdot
}^{i}$ is $\left( \mathcal{F}^{0}\right) _{\cdot }^{t}$-adapted, for $1\leq
i\leq N.$ For each $i,$ set $\left( X_{s}^{i},Y_{s}^{i},Z_{s}^{i}\right)
=\left( X_{s}^{i},\Theta _{s}^{i}\right) =\left(
X_{s}^{t,x^{i},u^{i}},Y_{s}^{t,x^{i},u^{i}},Z_{s}^{t,x^{i},u^{i}}\right) ,$
for $s\in \left[ t,T\right] $. Then $\left( X^{i},Y^{i},Z^{i}\right) $ is
the solution to the FBSDEs%
\begin{equation}
\left\{
\begin{array}{lll}
\mathrm{d}X_{s}^{i} & = & b\left( X_{s}^{i},\mathbb{P}%
_{X_{s}^{i}}^{W},u_{s}^{i}\right) \mathrm{d}s+\sigma \left( X_{s}^{i},%
\mathbb{P}_{X_{s}^{i}}^{W},u_{s}^{i}\right) \mathrm{d}W_{s}, \\
\mathrm{d}Y_{s}^{i} & = & -f^{0}\left( \mathbb{P}_{X_{s}^{i}}^{W},\Theta
_{s}^{i},\mathbb{P}_{\Theta _{s}^{i}},u_{s}^{i}\right) \mathrm{d}s+Z_{s}^{i}%
\mathrm{d}W_{s}, \\
X_{t}^{i} & = & x^{i},\text{ }Y_{T}^{i}=\Phi ^{0}\left( X_{T}^{i},\mathbb{P}%
_{X_{T}^{i}}^{W}\right) .%
\end{array}%
\right.  \label{fbsde2}
\end{equation}%
Immediately, we have%
\begin{eqnarray*}
\sum_{i=1}^{N}I_{A_{i}}X_{r}^{i}
&=&\sum_{i=1}^{N}I_{A_{i}}x_{i}+\int_{t}^{r}b\left(
\sum_{i=1}^{N}I_{A_{i}}X_{s}^{i},\mathbb{P}_{%
\sum_{i=1}^{N}I_{A_{i}}X_{s}^{i}}^{W},\sum_{i=1}^{N}I_{A_{i}}u_{s}^{i}%
\right) \mathrm{d}s \\
&&+\int_{t}^{r}\sigma \left( \sum_{i=1}^{N}I_{A_{i}}X_{s}^{i},\mathbb{P}%
_{\sum_{i=1}^{N}I_{A_{i}}X_{s}^{i}}^{W},\sum_{i=1}^{N}I_{A_{i}}u_{s}^{i}%
\right) \mathrm{d}W_{s}, \\
\sum_{i=1}^{N}I_{A_{i}}Y_{t}^{i} &=&\Phi ^{0}\left(
\sum_{i=1}^{N}I_{A_{i}}X_{T}^{i},\mathbb{P}_{%
\sum_{i=1}^{N}I_{A_{i}}X_{T}^{i}}^{W}\right) \\
&&+\int_{t}^{T}f^{0}\left( \mathbb{P}_{%
\sum_{i=1}^{N}I_{A_{i}}X_{T}^{i}}^{W},\Theta _{s}^{i},\mathbb{P}_{\Theta
_{s}^{i}},\sum_{i=1}^{N}I_{A_{i}}u_{s}^{i}\right) \mathrm{d}s \\
&&-\int_{t}^{T}\sum_{i=1}^{N}I_{A_{i}}Z_{s}^{i}\mathrm{d}W_{s},
\end{eqnarray*}%
where we have taken into account the fact $\sum_{i}\phi \left( x_{i}\right)
I_{A_{i}}=\phi \left( \sum_{i}x_{i}I_{A_{i}}\right) .$

Therefore,
\begin{equation}
Y_{t}^{t,\xi ;u}=\sum_{i=1}^{N}I_{A_{i}}Y_{t}^{i}.  \label{yed}
\end{equation}%
Clearly, from Lemma \ref{deter}, each $Y_{t}^{i}$ is deterministic. Now
taking the conditional expectation $\mathbb{E}^{\mathcal{F}_{t}^{0}}$ on
both sides of (\ref{yed}) yields
\begin{eqnarray*}
\mathbb{E}^{\mathcal{F}_{t}^{0}}\left[ Y_{t}^{t,\xi ;u}\right] &=&\mathbb{E}%
^{\mathcal{F}_{t}^{0}}\left[ \sum_{i=1}^{N}I_{A_{i}}Y_{t}^{i}\right] \\
&=&\sum_{i=1}^{N}Y_{t}^{i}\mathbb{E}^{\mathcal{F}_{t}^{0}}\left[ I_{A_{i}}%
\right] \\
&=&\sum_{i=1}^{N}Y_{t}^{i}I_{A_{i}}.
\end{eqnarray*}%
Therefore, for simple random variables, we have the desired result. For a
general $\xi \in L^{2}\left( \mathcal{F}_{t}^{0};\mathbb{R}^{n}\right) ,$ we
can take a sequence of simple random variables $\left\{ \xi _{i}\right\}
_{i\in N}$ which converge to $\xi $ in $L^{2}\left( \mathcal{F}_{t}^{0};%
\mathbb{R}^{n}\right) $. Consequently, from the estimates (\ref{best2}), we
have%
\begin{equation*}
\mathbb{E}\left[ \left\vert Y_{t}^{t,\xi _{i};u}-Y_{t}^{t,\xi ;u}\right\vert
^{2}\right] \leq C\mathbb{E}\left[ \left\vert \xi _{i}-\xi \right\vert ^{2}%
\right] ,\text{ as }i\rightarrow \infty .
\end{equation*}%
Thus the proof is complete. \hfill $\Box $

\subsection{Proof of Lemma \protect\ref{estvalue}}

\paragraph{Proof.}

Applying It\^{o}'s formula to $\left\vert Y_{t}^{t,\xi ;u}-Y_{t}^{t,\xi
^{\prime };u}\right\vert ^{2}$ and repeating the classical method, we have,
for some positive constant $C_{0},$
\begin{eqnarray}
\left\vert Y_{t}^{t,\xi ;u}-Y_{t}^{t,\xi ^{\prime };u}\right\vert &\leq
&C_{0}\mathbb{E}\left[ \left\vert \xi -\xi ^{\prime }\right\vert \right] ,
\label{ves1} \\
\left\vert Y_{t}^{t,\xi ;u}\right\vert &\leq &C_{0}\mathbb{E}\left[ \left(
1+\left\vert \xi \right\vert \right) \right] .  \label{ves2}
\end{eqnarray}%
On the other hand, for any arbitrary $\epsilon >0,$ there exists $u$ and $%
u^{\prime }\in \mathcal{U}_{ad}^{2}\left[ 0,T\right] $ such that%
\begin{equation*}
Y_{t}^{t,\xi ;u}-\epsilon \leq V\left( t,\xi \right) \leq Y_{t}^{t,\xi ;u},%
\text{ a.s,}
\end{equation*}%
and%
\begin{equation*}
Y_{t}^{t,\xi ^{\prime };u^{\prime }}-\epsilon \leq V\left( t,\xi ^{\prime
}\right) \leq Y_{t}^{t,\xi ^{\prime };u^{\prime }}\text{ a.s,}
\end{equation*}%
respectively. Combined with the estimates (\ref{ves1}) and (\ref{ves2}), we
have
\begin{equation*}
-C_{0}\mathbb{E}\left[ \left( 1+\left\vert \xi \right\vert \right) \right]
-\epsilon \leq Y_{t}^{t,\xi ;u}-\epsilon \leq V\left( t,\xi \right) \leq
Y_{t}^{t,\xi ;u}\leq C_{0}\mathbb{E}\left[ \left( 1+\left\vert \xi
\right\vert \right) \right] .
\end{equation*}%
Due to the arbitrariness of $\epsilon ,$ we get (\ref{est3}). Next, we deal
with%
\begin{equation*}
Y_{t}^{t,\xi ;u}-Y_{t}^{t,\xi ^{\prime };u}-\epsilon \leq V\left( t,\xi
\right) -V\left( t,\xi ^{\prime }\right) \leq Y_{t}^{t,\xi ;u^{\prime
}}-Y_{t}^{t,\xi ^{\prime };u^{\prime }}+\epsilon ,\text{ a.s.}
\end{equation*}%
Thus
\begin{equation*}
-C_{0}\mathbb{E}\left[ \left\vert \xi -\xi ^{\prime }\right\vert \right]
-\epsilon \leq V\left( t,\xi \right) -V\left( t,\xi ^{\prime }\right) \leq
C_{0}\mathbb{E}\left[ \left\vert \xi -\xi ^{\prime }\right\vert \right]
+\epsilon ,\text{ a.s.,}
\end{equation*}%
from which we immediately get (\ref{est1}). We are going to prove (\ref{est2}%
). From the definition of value function (\ref{d1}), for any given $%
\varepsilon >0,$ there exists $u^{\varepsilon }\in \mathcal{U}_{ad}^{2}\left[
0,T\right] $ such\ that
\begin{equation}
\mathbb{G}_{t,t+\delta }^{t,\xi ;u^{\varepsilon }}\left[ V\left( t+\delta
,X_{t+\delta }^{t,\xi ;u^{\varepsilon }}\right) \right] \geq V\left( t,\xi
\right) \geq \mathbb{G}_{t,t+\delta }^{t,\xi ;u^{\varepsilon }}\left[
V\left( t+\delta ,X_{t+\delta }^{t,\xi ;u^{\varepsilon }}\right) \right]
-\varepsilon .  \label{vbg}
\end{equation}%
We proceed the following inequalities by means of the above first (second)
inequality:%
\begin{equation*}
V\left( t,\xi \right) -V\left( t+\delta ,\xi \right) \geq -\tilde{C}\delta ^{%
\frac{1}{2}}\text{ (res., }\leq \tilde{C}\delta ^{\frac{1}{2}}\text{).}
\end{equation*}%
We prove the first inequality of (\ref{vbg}), the other one is analogous. It
follows from (\ref{estv1}) that, for an arbitrarily small $\varepsilon >0,$%
\begin{equation*}
V\left( t,\xi \right) -V\left( t+\delta ,\xi \right) \geq \mathcal{X}_{1}+%
\mathcal{X}_{2}-\varepsilon ,
\end{equation*}%
where
\begin{eqnarray*}
\mathcal{X}_{1} &\triangleq &\mathbb{G}_{t,t+\delta }^{t,\xi ;u^{\varepsilon
}}\left[ V\left( t+\delta ,X_{t+\delta }^{t,\xi ;u^{\varepsilon }}\right) %
\right] -\mathbb{G}_{t,t+\delta }^{t,\xi ;u^{\varepsilon }}\left[ V\left(
t+\delta ,\xi \right) \right] , \\
\mathcal{X}_{2} &\triangleq &\mathbb{G}_{t,t+\delta }^{t,\xi ;u^{\varepsilon
}}\left[ V\left( t+\delta ,\xi \right) \right] -V\left( t+\delta ,\xi
\right) .
\end{eqnarray*}%
From (\ref{est1}), there exists a constant $C$ independent of $%
u^{\varepsilon },$ such that%
\begin{eqnarray*}
\left\vert \mathcal{X}_{1}\right\vert &\leq &C\left[ \mathbb{E}\left[
\left\vert V\left( t+\delta ,X_{t+\delta }^{t,\xi ;u^{\varepsilon }}\right)
-V\left( t+\delta ,\xi \right) \right\vert ^{2}\right] \right] ^{\frac{1}{2}}
\\
&\leq &C\left[ \mathbb{E}\left[ \left\vert X_{t+\delta }^{t,\xi
;u^{\varepsilon }}-\xi \right\vert \right] \right] ^{\frac{1}{2}}.
\end{eqnarray*}%
From \cite{bbp1997}, applying Proposition 1.1 to c\`{a}dl\`{a}g process $%
X^{t,\xi ,u^{\varepsilon }}$, we obtain
\begin{equation*}
\mathbb{E}\left[ \left\vert X_{t+\delta }^{t,\xi ,u^{\varepsilon }}-\xi
\right\vert ^{2}\right] \leq C\delta \mathbb{E}\left[ \left( 1+\left\vert
\xi \right\vert ^{2}\right) \right] .
\end{equation*}%
Instantly, we get%
\begin{equation}
\mathcal{X}_{1}\geq -C\delta ^{\frac{1}{2}}\mathbb{E}\left[ \left(
1+\left\vert \xi \right\vert ^{2}\right) \right] .  \label{ee1}
\end{equation}%
Meanwhile$,$
\begin{eqnarray*}
\mathcal{X}_{2} &=&\mathbb{G}_{t,t+\delta }^{t,\xi ;u^{\varepsilon }}\left[
V\left( t+\delta ,\xi \right) \right] -V\left( t+\delta ,\xi \right) \\
&=&\mathbb{E}^{\mathcal{F}_{t}^{0}}\Bigg [V\left( t+\delta ,\mu \right)
+\int_{t}^{t+\delta }f^{0}\left( \rho _{s}^{t,\mu ;u^{\varepsilon }},\tilde{%
\Theta}_{s}^{t,\xi ,u^{\varepsilon }},\mathbb{P}_{\tilde{\Theta}_{s}^{t,\xi
,u^{\varepsilon }}},u_{s}^{\varepsilon }\right) \mathrm{d}s \\
&&-\int_{t}^{t+\delta }\tilde{Z}_{s}^{t,\xi ;u^{\varepsilon }}\mathrm{d}W_{s}%
\Bigg ]-V\left( t+\delta ,\mu \right) \\
&=&\mathbb{E}^{\mathcal{F}_{t}^{0}}\Bigg [\int_{t}^{t+\delta }f^{0}\left(
\rho _{s}^{t,\mu ;u^{\varepsilon }},\tilde{\Theta}_{s}^{t,\xi
,u^{\varepsilon }},\mathbb{P}_{\tilde{\Theta}_{s}^{t,\xi ,u^{\varepsilon
}}},u_{s}^{\varepsilon }\right) \mathrm{d}s\Bigg ],
\end{eqnarray*}%
where $\left( \tilde{Y}_{s}^{t,\xi ;u^{\varepsilon }},\tilde{Z}_{s}^{t,\xi
;u^{\varepsilon }}\right) _{t\leq s\leq t+\delta }$ is the solution to BSDE (%
\ref{bsg}) with the terminal condition $\eta =V\left( t+\delta ,\xi \right)
. $ Employing the classical method (Schwartz inequality), we have that, for
certain constant $C$ independent of $t$ and $\delta ,$%
\begin{eqnarray}
\left\vert \mathcal{X}_{2}\right\vert ^{2} &\leq &\delta ^{\frac{1}{2}}%
\mathbb{E}^{\mathcal{F}_{t}^{0}}\Bigg [\int_{t}^{t+\delta }\left\vert
f^{0}\left( \rho _{s}^{t,\mu ;u^{\varepsilon }},\tilde{\Theta}_{s}^{t,\xi
,u^{\varepsilon }},\mathbb{P}_{\tilde{\Theta}_{s}^{t,\xi ,u^{\varepsilon
}}},u_{s}^{\varepsilon }\right) \right\vert \mathrm{d}s\Bigg ]  \notag \\
&\leq &C\delta ^{\frac{1}{2}}.  \label{ee2}
\end{eqnarray}%
From the inequalities, (\ref{ee1}) and (\ref{ee2}), it follows%
\begin{eqnarray*}
V\left( t,\xi \right) -V\left( t+\delta ,\xi \right) &\geq &\mathcal{X}_{1}+%
\mathcal{X}_{2}-\varepsilon \\
&\geq &-C\delta ^{\frac{1}{2}}\mathbb{E}\left[ \left( 1+\left\vert \xi
\right\vert ^{2}\right) \right] -\varepsilon .
\end{eqnarray*}%
The proof ends from letting $\varepsilon \rightarrow 0.$ \hfill $\Box $

\subsection{Proof of Lemma \protect\ref{estv}}

\paragraph{Proof.}

The idea is similar to\ Lemma \ref{estvalue}. From Lemma \ref{p2}, for any
given $\varepsilon >0,$ there exists $u^{\varepsilon }\in \mathcal{U}%
_{ad}^{2}\left[ 0,T\right] $ such\ that
\begin{equation}
\mathbb{G}_{t,t+\delta }^{t,\mu ;u^{\varepsilon }}\left[ \mathcal{V}\left(
t+\delta ,\rho _{t+\delta }^{t,\mu ;u^{\varepsilon }}\right) \right] \geq
\mathcal{V}\left( t,\mu \right) \geq \mathbb{G}_{t,t+\delta }^{t,\mu
;u^{\varepsilon }}\left[ \mathcal{V}\left( t+\delta ,\rho _{t+\delta
}^{t,\mu ;u^{\varepsilon }}\right) \right] -\varepsilon .  \label{estv1}
\end{equation}

In the following, we prove the following inequalities with the help of the
above first (second) inequality:%
\begin{equation}
\mathcal{V}\left( t,\mu \right) -\mathcal{V}\left( t+\delta ,\mu \right)
\geq -\tilde{C}\delta ^{\frac{1}{2}}\text{ (res., }\leq \tilde{C}\delta ^{%
\frac{1}{2}}\text{).}  \label{estv2}
\end{equation}%
We prove the first inequality of (\ref{estv2}), the other one is analogous.
It follows from (\ref{estv1}) that, for an arbitrarily small $\varepsilon
>0, $%
\begin{equation*}
\mathcal{V}\left( t,\mu \right) -\mathcal{V}\left( t+\delta ,\mu \right)
\geq \mathcal{I}_{1}+\mathcal{I}_{2}-\varepsilon ,
\end{equation*}%
where
\begin{eqnarray*}
\mathcal{I}_{1} &\triangleq &\mathbb{G}_{t,\tau }^{t,\mu ;u^{\varepsilon }}%
\left[ \mathcal{V}\left( t+\delta ,\rho _{t+\delta }^{t,\mu ;u^{\varepsilon
}}\right) \right] -\mathbb{G}_{t,\tau }^{t,\mu ;u^{\varepsilon }}\left[
\mathcal{V}\left( t+\delta ,\mu \right) \right] , \\
\mathcal{I}_{2} &\triangleq &\mathbb{G}_{t,\tau }^{t,\mu ;u^{\varepsilon }}%
\left[ \mathcal{V}\left( t+\delta ,\mu \right) \right] -\mathcal{V}\left(
t+\delta ,\mu \right) .
\end{eqnarray*}%
From Lemma \ref{p1} and \ref{estbsde}, we get that, for some constant $C$
independent of $u^{\varepsilon },$%
\begin{eqnarray*}
\left\vert \mathcal{I}_{1}\right\vert &\leq &C\left[ \mathbb{E}\left[
\left\vert \mathcal{V}\left( t+\delta ,\rho _{t+\delta }^{t,\mu
;u^{\varepsilon }}\right) -\mathcal{V}\left( t+\delta ,\mu \right)
\right\vert ^{2}\right] \right] ^{\frac{1}{2}} \\
&\leq &C\left[ \mathbb{E}\left[ \mathcal{W}_{2}^{2}\left( \rho _{t+\delta
}^{t,\mu ;u^{\varepsilon }},\mu \right) \right] \right] ^{\frac{1}{2}}.
\end{eqnarray*}%
From \cite{bbp1997}, applying Proposition 1.1 to c\`{a}dl\`{a}g process $%
X^{t,\xi ,u^{\varepsilon }}$, we obtain
\begin{equation*}
\mathbb{E}\left[ \left\vert X_{t+\delta }^{t,\xi ,u^{\varepsilon }}-\xi
\right\vert ^{2}\right] \leq C\delta \mathbb{E}\left[ \left( 1+\left\vert
\xi \right\vert ^{2}\right) \right] .
\end{equation*}%
Meanwhile, from the definition of $\mathcal{W}_{2},$ it yields%
\begin{equation*}
\mathbb{E}\left[ \mathcal{W}_{2}^{2}\left( \rho _{t+\delta }^{t,\mu
;u^{\varepsilon }},\mu \right) \right] \leq C\delta \mathbb{E}\left[ \left(
1+\left\vert \xi \right\vert ^{2}\right) \right] .
\end{equation*}%
Now taking into account the arbitrariness of $\xi \in L^{2}\left( \mathcal{G}%
_{t};\mathbb{R}^{n}\right) ,$ we get%
\begin{equation*}
\mathbb{E}\left[ \mathcal{W}_{2}^{2}\left( \rho _{t+\delta }^{t,\mu
;u^{\varepsilon }},\mu \right) \right] \leq C\delta \left( 1+\mathcal{W}%
_{2}^{2}\left( \mu ,\delta _{0}\right) \right) ,
\end{equation*}%
from which we immediately derive that%
\begin{equation}
\mathcal{I}_{1}\geq -C\delta ^{\frac{1}{2}}\left( 1+\mathcal{W}_{2}\left(
\mu ,\delta _{0}\right) \right) .  \label{estv3}
\end{equation}%
We now proceed the second term $\mathcal{I}_{2}.$ According to the
definition of $\mathbb{G}_{t,t+\delta }^{t,\mu ;u^{\varepsilon }}\left[
\mathcal{\cdot }\right] ,$ we deduce that $\mathcal{I}_{2}$ can be expressed
as
\begin{eqnarray*}
\mathcal{I}_{2} &=&\mathbb{G}_{t,t+\delta }^{t,\mu ;u^{\varepsilon }}\left[
\mathcal{V}\left( t+\delta ,\mu \right) \right] -\mathcal{V}\left( t+\delta
,\mu \right) \\
&=&\mathbb{E}^{\mathcal{F}_{t}^{0}}\Bigg [\mathcal{V}\left( t+\delta ,\mu
\right) +\int_{t}^{t+\delta }f^{0}\left( \rho _{s}^{t,\mu ;u^{\varepsilon }},%
\tilde{\Theta}_{s}^{t,\mu ,u^{\varepsilon }},\mathbb{P}_{\tilde{\Theta}%
_{s}^{t,\mu ,u^{\varepsilon }}},u_{s}^{\varepsilon }\right) \mathrm{d}s \\
&&-\int_{t}^{t+\delta }\tilde{Z}_{s}^{t,\mu ;u^{\varepsilon }}\mathrm{d}W_{s}%
\Bigg ]-\mathcal{V}\left( t+\delta ,\mu \right) \\
&=&\mathbb{E}^{\mathcal{F}_{t}^{0}}\Bigg [\int_{t}^{t+\delta }f^{0}\left(
\rho _{s}^{t,\mu ;u^{\varepsilon }},\tilde{\Theta}_{s}^{t,\mu
,u^{\varepsilon }},\mathbb{P}_{\tilde{\Theta}_{s}^{t,\mu ,u^{\varepsilon
}}},u_{s}^{\varepsilon }\right) \mathrm{d}s\Bigg ],
\end{eqnarray*}%
where $\left( \tilde{Y}_{s}^{t,\mu ;u^{\varepsilon }},\tilde{Z}_{s}^{t,\mu
;u^{\varepsilon }}\right) _{t\leq s\leq t+\delta }$ is the solution of BSDE (%
\ref{bsg}) with the terminal condition $\eta =\mathcal{V}\left( t+\delta
,\mu \right) .$ Employing the classical method (Schwartz inequality), from
Lemma \ref{estbsde}, we have that, for certain constant $C$ independent of $%
t $ and $\delta ,$%
\begin{eqnarray*}
\left\vert \mathcal{I}_{2}\right\vert ^{2} &\leq &\delta ^{\frac{1}{2}}%
\mathbb{E}^{\mathcal{F}_{t}^{0}}\Bigg [\int_{t}^{t+\delta }\left\vert
f^{0}\left( \rho _{s}^{t,\mu ;u^{\varepsilon }},\tilde{\Theta}_{s}^{t,\xi
,u^{\varepsilon }},\mathbb{P}_{\tilde{\Theta}_{s}^{t,\xi ,u^{\varepsilon
}}},u_{s}^{\varepsilon }\right) \right\vert \mathrm{d}s\Bigg ] \\
&\leq &C\delta ^{\frac{1}{2}}.
\end{eqnarray*}%
Consequently, combining (\ref{estv3}), we have%
\begin{eqnarray*}
\mathcal{V}\left( t,\mu \right) -\mathcal{V}\left( t+\delta ,\mu \right)
&\geq &\mathcal{I}_{1}+\mathcal{I}_{2}-\varepsilon \\
&\geq &-C\delta ^{\frac{1}{2}}\left( 1+\mathcal{W}_{2}\left( \mu ,\delta
_{0}\right) \right) -\varepsilon .
\end{eqnarray*}%
The desired result follows from letting $\varepsilon \rightarrow 0.$ The
proof is then complete. \hfill $\Box $

\subsection{Proof of Lemma \protect\ref{estvis}}

\paragraph{Proof.}

From \cite{bbp1997}, applying Proposition 1.1 to c\`{a}dl\`{a}g process $%
X^{t,\xi ,u}$, we have, for $p\geq 2$ the existence of constant $C_{p}>0$
such that%
\begin{equation*}
\mathbb{E}\left[ \sup_{t\leq s\leq t+\delta }\left\vert X^{t,\xi ,u}-\xi
\right\vert ^{p}\right] \leq C_{p}\delta \left( 1+\left\vert \xi \right\vert
^{p}\right) ,\text{ uniformly in }u\in \mathcal{U}_{ad}^{2}\left[ 0,T\right]
.
\end{equation*}%
Therefore, whenever $\delta \downarrow 0,$ the random variable $\kappa
^{\delta }=\sup_{t\leq s\leq t+\delta }\left\vert X^{t,\xi ,u}-\xi
\right\vert $ converges monotone $0.$ Denote $\widetilde{\mathscr{Y}}%
_{s}^{u}=\mathscr{Y}_{s}^{1,u}-\mathscr{Y}_{s}^{2,u}$ and $\widetilde{%
\mathscr{Z}}_{s}^{u}=\mathscr{Z}_{s}^{1,u}-\mathscr{Z}_{s}^{2,u}$for $t\leq
s\leq t+\delta .$
\begin{equation}
\left\{
\begin{array}{rcl}
-\mathrm{d}\widetilde{\mathscr{Y}}_{s}^{u} & = & \Big [F\left( s,\mathscr{X}%
_{s}^{t,\xi ;u},\mathscr{Y}_{s}^{1,u},\mathscr{Z}_{s}^{1,u},u_{s}\right)
-F\left( s,\mathscr{X}_{s}^{t,\xi ;u},\mathscr{Y}_{s}^{2,u},\mathscr{Z}%
_{s}^{2,u},u_{s}\right) \\
&  & -F\left( s,\xi ,\mathscr{Y}_{s}^{2,u},\mathscr{Z}_{s}^{2,u},u_{s}%
\right) +F\left( s,\mathscr{X}_{s}^{t,\xi ;u},\mathscr{Y}_{s}^{2,u},%
\mathscr{Z}_{s}^{2,u},u_{s}\right) \Big ]\mathrm{d}s \\
&  & -\widetilde{\mathscr{Z}}_{s}^{u}\mathrm{d}W_{s},\text{ } \\
\widetilde{\mathscr{Y}}_{t+\delta }^{u} & = & 0.%
\end{array}%
\right.  \label{bs4}
\end{equation}%
Applying It\^{o}'s formula to $\left\vert \widetilde{\mathscr{Y}}%
_{s}^{u}\right\vert ^{2}e^{\beta \left( s-t\right) }$ on the interval $\left[
t,T\right] ,$ we get the following estimate with same technique used in
Lemma \ref{estvalue}
\begin{eqnarray*}
\mathbb{E}\left[ \int_{t}^{t+\delta }\left( \left\vert \widetilde{\mathscr{Y}%
}_{s}^{u}\right\vert ^{2}+\left\vert \widetilde{\mathscr{Z}}%
_{s}^{u}\right\vert ^{2}\right) \mathrm{d}s\right] &\leq &C\mathbb{E}\left[
\int_{t}^{t+\delta }\rho _{0}\left( \left\vert \mathscr{X}_{s}^{t,\xi
;u}-\xi \right\vert ^{2}\right) \mathrm{d}s\right] \\
&\leq &C\delta \rho _{0}\left( \kappa ^{\delta }\right) ^{2}.
\end{eqnarray*}%
For some $\rho _{0}:\mathbb{R}\rightarrow \mathbb{R}$ with $\rho _{0}\left(
\varepsilon \right) \downarrow 0$ as $\varepsilon \downarrow 0.$ On the
other hand, from Lemma \ref{vdet}, we know that $\mathscr{Y}_{t}^{1,u}$ and $%
\mathscr{Y}_{t}^{2,u}$ are deterministic whenever $u\in \mathcal{U}_{ad}^{2}%
\left[ 0,T\right] $ is $\mathcal{F}_{s}^{t}$-adapted. It yields%
\begin{eqnarray*}
\left\vert \widetilde{\mathscr{Y}}_{t}^{u}\right\vert &=&\left\vert \mathbb{E%
}\left[ \widetilde{\mathscr{Y}}_{t}^{u}\right] \right\vert \\
&=&\Bigg |\mathbb{E}\Bigg [\int_{t}^{t+\delta }F\left( s,\mathscr{X}%
_{s}^{t,\xi ;u},\mathscr{Y}_{s}^{1,u},\mathscr{Z}_{s}^{1,u},u_{s}\right)
-F\left( s,\xi ,\mathscr{Y}_{s}^{2,u},\mathscr{Z}_{s}^{2,u},u_{s}\right)
\mathrm{d}s\Bigg ]\Bigg | \\
&\leq &C\mathbb{E}\Bigg [\int_{t}^{t+\delta }\Big [\rho _{0}\left(
\left\vert \mathscr{X}_{s}^{t,\xi ;u}-\xi \right\vert \right) +\left\vert %
\mathscr{Y}_{s}^{1,u}-\mathscr{Y}_{s}^{2,u}\right\vert
\end{eqnarray*}%
\begin{eqnarray*}
&&+\left\vert \mathscr{Z}_{s}^{1,u}-\mathscr{Z}_{s}^{2,u}\right\vert +%
\mathcal{W}_{2}\left( \mathbb{P}_{\left( \mathscr{Y}_{s}^{1,u},\mathscr{Z}%
_{s}^{1,u}\right) },\mathbb{P}_{\left( \mathscr{Y}_{s}^{2,u},\mathscr{Z}%
_{s}^{2,u}\right) }\right) \Big ]\mathrm{d}s\Bigg ] \\
&\leq &C\delta \mathbb{E}\left[ \rho _{0}\left( \kappa ^{\delta }\right) %
\right] \\
&&+C\delta ^{\frac{1}{2}}\Bigg \{\mathbb{E}\int_{t}^{t+\delta }\left(
\left\vert \mathscr{Y}_{s}^{1,u}-\mathscr{Y}_{s}^{2,u}\right\vert
^{2}+\left\vert \mathscr{Z}_{s}^{1,u}-\mathscr{Z}_{s}^{2,u}\right\vert
^{2}\right) \mathrm{d}s\Bigg \}^{\frac{1}{2}},
\end{eqnarray*}
which implies that
\begin{equation*}
\left\vert \widetilde{\mathscr{Y}}_{t}^{u}\right\vert \leq C\delta \left\{
\mathbb{E}\left[ \rho _{0}\left( \kappa ^{\delta }\right) \right] +\left[
\mathbb{E}\left[ \rho _{0}^{2}\left( \kappa ^{\delta }\right) \right] \right]
^{\frac{1}{2}}\right\} .
\end{equation*}%
Observe that, for each $\delta >0,$ $\rho _{0}\left( \kappa ^{\delta
}\right) $ is square integrable. Thus we set $\rho _{1}\left( \delta \right)
=\mathbb{E}\left[ \rho _{0}\left( \kappa ^{\delta }\right) \right] +\left[
\mathbb{E}\left[ \rho _{0}^{2}\left( \kappa ^{\delta }\right) \right] \right]
^{\frac{1}{2}}.$ For the general $u\in \mathcal{U}_{ad}^{2}\left[ 0,T\right]
$, since $\mathscr{Y}_{s}^{1,u}$ and $\mathscr{Y}_{s}^{2,u}$ are continuous
in $u,$ the proof is completed via the same technique employed before.
\hfill $\Box $

\bigskip

\noindent \textbf{Acknowledgement.} The author thanks the Handling Editor,
Prof. Enrique Zuazua and anonymous reviewers for their many good comments
and suggestions. The author also thanks Dr. Xiaoli Wei for her conversations
and suggestions.

\end{document}